\documentclass[twoside]{irmaems}
\usepackage{amsmath, amssymb, amscd,amsopn} 
\usepackage{latexsym}
\usepackage{graphicx} 
\usepackage{makeidx} 

\renewcommand\theequation{\arabic{section}.\arabic{equation}}

\theoremstyle{definition} %%% for statements in roman typeface

 \newtheorem{definition}{Definition}[section]

\newtheorem*{note}{Note}
\theoremstyle{plain}      %%% for statements in italic typeface

 \newtheorem{proposition}[definition]{Proposition}
 \newtheorem{theorem}[definition]{Theorem}
 \newtheorem{corollary}[definition]{Corollary}
 \newtheorem{lemma}[definition]{Lemma}

\newcommand{\GL}{\operatorname{GL}}
\newcommand{\End}{\operatorname{End}}

\begin{document}

\def\R{{\mathbb R}} \def\Sym{{\rm Sym}} \def\proj{{\rm proj}}
\def\BB{{\mathcal B}}\def\Id{{\rm Id}} \def\PP{{\mathcal P}}
\def\Z{{\mathbb Z}} \def\C{{\mathbb C}} \def\Diag{{\rm Diag}}
\def\z{{\bf z}} \def\Sylv{{\rm Sylv}} \def\Res{{\rm Res}}
\def\Disc{{\rm Disc}} \def\DD{{\mathcal D}} \def\u{{\bf u}}
\def\Homeo{{\rm Homeo}} \def\MM{{\mathcal M}}
\def\x{{\bf x}} \def\D{{\mathbb D}} \def\Aut{{\rm Aut}}
\def\Inn{{\rm Inn}} \def\Out{{\rm Out}} \def\N{{\mathbb N}}
\def\pprod{{\rm prod}} \def\Cub{{\rm Cub}} \def\sgn{{\rm sgn}}
\def\EE{{\mathcal E}} \def\RR{{\mathcal R}} \def\Im{{\rm Im}}
\def\SS{{\mathcal S}} \def\LL{{\mathcal L}} \def\AA{{\mathcal A}}
\def\F{{\mathbb F}} \def\Ker{{\rm Ker}} \def\cd{{\rm cd}}
\def\Q{{\mathbb Q}} \def\FF{{\mathcal F}} \def\CC{{\mathcal C}}
\def\XX{{\mathcal X}} \def\Sal{{\rm Sal}} \def\B{{\mathbb B}}
\def\codim{{\rm codim}} \def\min{{\rm min}} \def\Cox{{\rm Cox}}
\def\Min{{\rm Min}} \def\S{{\mathbb S}} \def\K{{\mathbb K}}
\def\VV{{\mathcal V}} \def\End{{\rm End}} \def\GL {{\rm GL}}
\def\PV{{\rm PV}} \def\St{{\rm St}} \def\pos{{\rm pos}}
\def\An{{\rm An}} \def\CP{{\rm CP}} \def\A{{\mathbb A}}

\title{Braid groups and Artin groups}

\author{Luis Paris}

\address{
Institut de Math\'ematiques de Bourgogne
-- UMR 5584 du CNRS\\
Universit\'e de Bourgogne,
BP 47870\\
21078 Dijon cedex, France\\
email:\,\tt{lparis@u-bourgogne.fr}
}

\maketitle

\tableofcontents   
%%%%%%%%%%%%%%%%%%%%%%%%%%%%%%%%%%%%%

\section{Introduction}

The braids go back to several centuries and were universally used for ornamental purposes or even 
practical ones, for example in the fashioning of ropes. Now, they are described by means of abstract 
models known under the name of ``theory of braids''. The theory of braids studies the concept of braids 
(such as we imagine them) as well as various generalizations arising from various branches of the 
mathematics. The idea is that the braids form a group. The number of strands must be fixed so that the 
operation is well-defined. So, we have a braid group on two strands, a braid group on three strands, 
and so on. The braid group on one strand is trivial because a string cannot be braided (although it can 
be knotted). 

We generally make the mathematical study of braids go back to an article of Emile Artin \cite{Artin1} 
dated from 1925, in which is described the notion of braids under various aspects, one being that 
obvious, like a ``series of tended and interlaced strings'', and others more conceptual but equally 
deep, such as a presentation by generators and relations, or a presentation as the mapping class group 
of a punctured disk.

Since the 30s, a strong link between braids and links (and knots) were established by people such as 
Alexander and Markov (see \cite{Birma1}). This link is at the origin in the 80s of a deep revival in 
the theory of knots with the work of Jones and his invariant defined from the theory of braids (see 
\cite{Jones1}, \cite{Jones2}, \cite{FYHLMO1}, and \cite{PrzTra1}).

Later, interesting relations with the algebraic geometry and the theory of finite groups generated 
by reflections were established, in particular by Arnol'd \cite{Arnol2}, \cite{Arnol3}, \cite{Arnol4} 
and Brieskorn \cite{Bries3}, \cite{Bries2}. These relations become particularly interesting when we 
extend the notion of braid groups to that of Artin groups of spherical type, also called generalized 
braid groups. Although the Artin groups were introduced by Tits \cite{Tits3} as extensions of Coxeter 
groups, their study really began in the seventies with the works of Brieskorn \cite{Bries1}, 
\cite{Bries2}, Saito \cite{BriSai1} and Deligne \cite{Delig1}, where different aspects of these groups 
are studied, such as their combinatorics, as well as their link with the hyperplane arrangements and 
the singularities.

Some problems in group theory, often very close to the algorithmics, such as the word and conjugacy 
problems, have a renewal of interest not only through their applications in the other domains, but also 
because the notion of mathematical demonstration is changing. Indeed, we distinguish now the notion of
demonstration from the notion of effective demonstration, the one which builds up the solution. Such a 
demonstration gives rise to an algorithm, and its complexity (calculation time) is of importance. The 
algorithmics in the braid groups is especially active. Problems of decision such as 
the conjugacy problem were solved by Garside \cite{Garsi1} in 1969 with methods which are now the 
source of numerous works on the braid groups. In \cite{DehPar1} is introduced a more formal and more 
general framework to study algorithmic problems on the braid groups: the Garside groups. The idea is to 
isolate certain combinatorial properties of the braid groups, in particular these emphasized by Garside 
\cite{Garsi1}. It is a less restrictive model which uses tools from the language theory (monoids, 
rewriting systems) and the combinatorics (ordered sets), tools that are especially adapted to treat 
algorithmic problems. Now, the major part of the algorithmic problems on the braid groups are studied 
within the framework of the Garside groups. Also, let us indicate that the Artin groups of spherical 
type are Garside groups.

This survey is written from these viewpoints but also maintaining two other objectives: (1) to make a survey 
understandable by non-specialists; (2) to make as often as possible the link with the mapping class 
groups.

The first section is about the ``classical'' theory of braid groups. Various aspects as well as some 
properties of them are presented. The second section is an introduction to the Artin groups, and the 
third is an introduction to the Garside groups. There, the reader will find algorithms to solve some 
decision problems such as the conjugacy one for the braid groups (and Garside groups).

The fifth section is about the cohomology of Artin groups, although the exposition goes beyond by 
explaining the Salvetti complexes. These are tools originally from the theory of hyperplane 
arrangements that turn out to be useful in the context of the braid groups.

The sixth section is about the linear representations of the braid groups studied by Bigelow 
\cite{Bigel1} and Krammer \cite{Kramm4}, \cite{Kramm3}, as well as about its various generalizations (to 
the Artin groups). Both, the algebraic aspect and the topological aspect of these representations, are 
explained. Other linear representations of the braid groups have been studied and are also interesting 
but, for lack of place and for reason of coherence, these will not be treated in this text. We refer to 
\cite{BirBre1} for a survey on the other linear representations.

The seventh section is about the geometric representations of the Artin groups. (By a geometric 
representation we simply mean a homomorphism in a mapping class group.) This subject is less popular 
than the previous ones but I strongly believe in its future. In particular, Subsection 7.3, where are 
explained the results of Castel \cite{Caste2}, shows all the power of such a study.

Finally, I would like to indicate two aspects of the braid groups which are not in this survey and 
which ``should be in any survey on the braid groups''.

The first aspect is the link of the braids with links and knots. This is very important in the theory 
but amply explained in all the books and almost all the surveys on the subject. So, I voluntarily 
ignore this aspect in order to be able to treat in a more detailed way the other ones. The reader will 
find in \cite{Birma1}, \cite{Hanse1}, \cite{MurKur1}, \cite{KasTur1} detailed expositions on this 
aspect and on the braid groups in general.

I would have wanted to make an eighth section to explain the second aspect: the orders in the braid 
groups. But, unfortunately, this article is long enough and there is no more room for another section. 
Inspired by problems of set theory, Dehornoy \cite{Dehor3} founded an explicit construction of a total 
ordering invariant by left multiplication in the braid group. The fact that the braid group is 
orderable is not maybe completely new, in the sense that it results from Nielsen theory \cite{Niels1}, 
but Dehornoy's ordering is interesting in itself. In my opinion, it is an important tool to understand 
the braid groups, and I augur numerous developments in this direction. The Artin groups of type $B_n$ 
and $\tilde A_n$ embed into braid groups (see Section 3) thus they are also orderable. The Artin groups 
of type $D_n$ embed into mapping class groups of surfaces with boundary (see Section 7), and, by 
\cite{RouWie1}, such a group is orderable. We do not know whether the other Artin groups are orderable 
or not. We encourage the reader to consult \cite{DeDyRoWi1} for a detailed discussion on this subject.

%%%%%%%%%%%%%%%%%%%%%%%%%%%%%%%%%%%%%%%%%%%%%%%%%%%%%%%%%%%%%%%%%%%%%%%%

\section{Braid groups}

\subsection{Braids}

Let $n\ge 1$ be an integer, and let $P_1, \dots, P_n$ be $n$ distinct points in the plane $\R^2$ 
(except mention of the contrary, we will always assume $P_k= (k,0)$ for all $1 \le k\le n$). Define a {\it 
braid on $n$ strands} \index{Braid} to be a $n$-tuple $\beta = (b_1, \dots, b_n)$ of paths, $b_k : [0,1] 
\to \R^2$, such that
\begin{itemize}
\item
$b_k(0)=P_k$ for all $1 \le k\le n$;
\item
there exists a permutation $\chi = \theta (\beta) \in \Sym_n$ such that $b_k(1)= P_{\chi (k)}$ for all 
$1 \le k \le n$;
\item
$b_k(t) \neq b_l(t)$ for all $k \neq l$ and all $t \in [0,1]$.
\end{itemize}
Two braids $\alpha$ and $\beta$ are said to be {\it homotopic} if there exists a continuous family $\{ 
\gamma_s \}_{s \in [0,1]}$ of braids such that $\gamma_0 = \alpha$ and $\gamma_1 = \beta$. Note that 
$\theta (\alpha) = \theta (\beta)$ if $\alpha$ and $\beta$ are homotopic.

We represent graphically a homotopy class of braids as follows. Let $I_k$ be a copy of the interval 
$[0,1]$. Take a braid $\beta= (b_1, \dots, b_n)$ and define the {\it geometric braid} \index{Geometric braid}
\[
\beta^g : I_1 \sqcup \cdots \sqcup I_n \to \R \times [0,1]
\]
by $\beta^g(t) = (b_k(t),t)$ for all $t \in I_k$ and all $1 \le k\le n$. Let $\proj: \R^2 \times [0,1] 
\to \R \times [0,1]$ be the projection defined by
\[
\proj (x,y,t)= (x,t)\,.
\]
Up to homotopy, we can assume that $\proj \circ \beta^g$ is a smooth immersion with only transversal 
double points that we call {\it crossings}. In each crossing we indicate graphically like in Figure~2.1 
which strand goes over the other. Such a representation of $\beta$ is called a {\it braid diagram} 
\index{Braid diagram} of 
$\beta$. An example is illustrated in Figure 2.2.

%%%%%%%%%%%
\begin{figure}[htb]
\centerline{
\setlength{\unitlength}{.4cm}
\begin{picture}(9,5)
\put(0,2){\includegraphics[width=3.6cm]{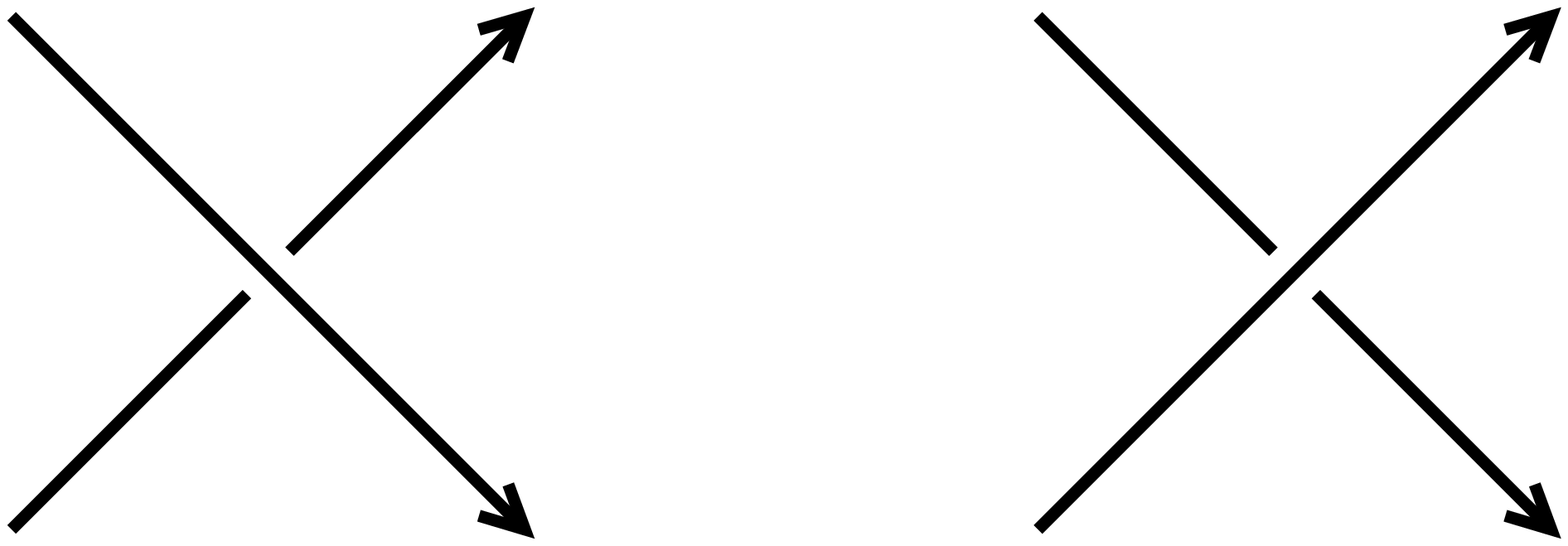}}
\put(0,1){\small positive}
\put(0,0){\small crossing}
\put(6,1){\small negative}
\put(6,0){\small crossing}
\end{picture}}
\medskip
\centerline{{\bf Figure 2.1.} Crossings in a braid diagram.}
\end{figure}
%%%%%%%%%

%%%%%%%
\begin{figure}[htb]
\centerline{
\setlength{\unitlength}{.4cm}
\begin{picture}(16,10)
\put(0,1){\includegraphics[width=6cm]{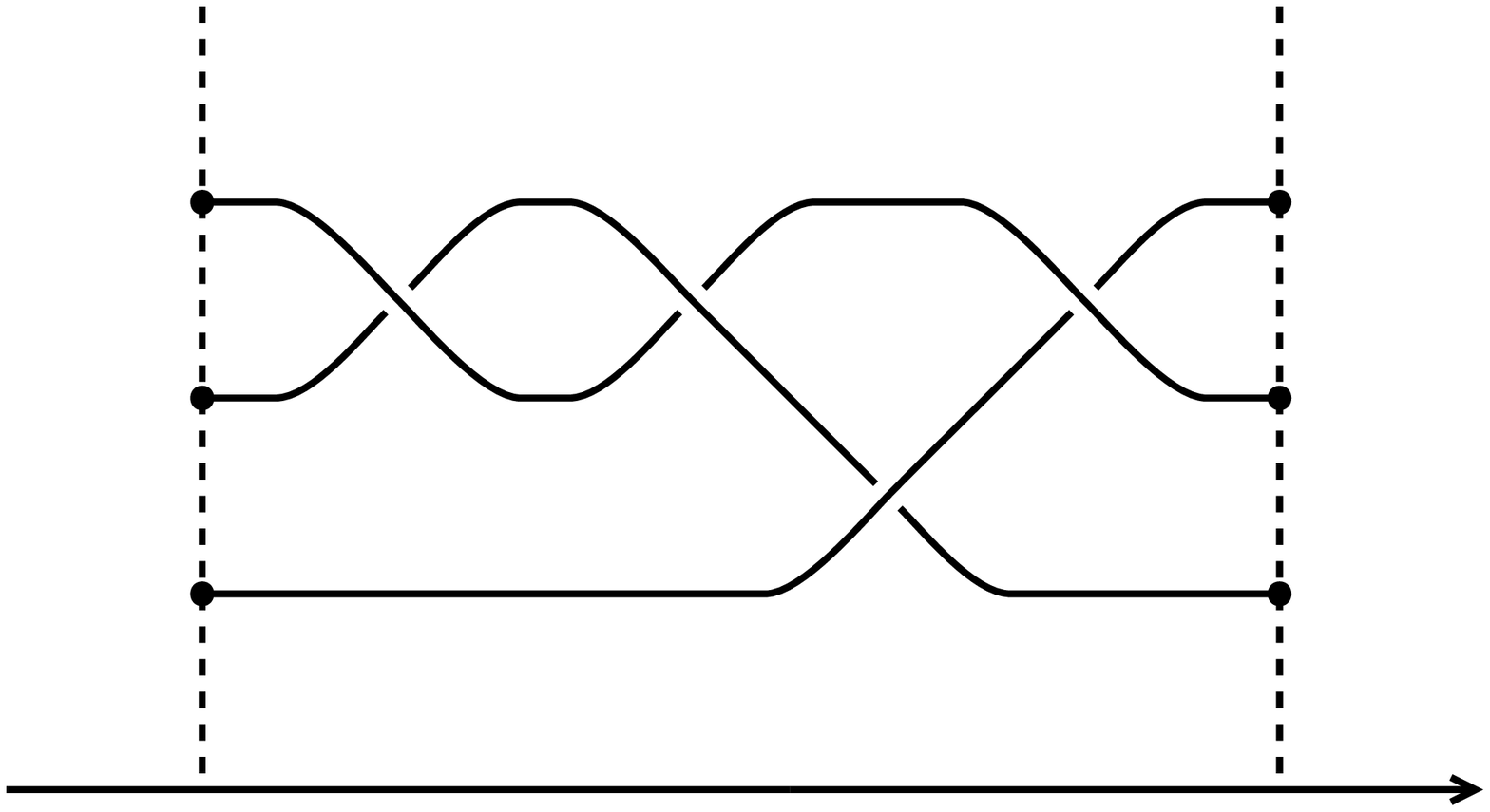}}
\put(1.9,0.3){\small $0$}
\put(12.9,0.3){\small $1$}
\put(1,2.9){\small $P_1$}
\put(1,4.9){\small $P_2$}
\put(1,6.9){\small $P_3$}
\put(13.3,2.9){\small $P_1$}
\put(13.3,4.9){\small $P_2$}
\put(13.3,6.9){\small $P_3$}
\put(15.5,1){\small $t$}
\put(1,9.5){\small $\R^2 \times \{0\}$}
\put(12,9.5){\small $\R^2 \times \{1\}$}
\end{picture}}
\medskip
\centerline{{\bf Figure 2.2.} A braid diagram.}
\end{figure}
%%%%%%%%%%

The {\it product} of two braids $\alpha = (a_1, \dots, a_n)$ and $\beta = (b_1, \dots, b_n)$ is defined 
to be the braid
\[
\alpha \cdot \beta = (a_1 b_{\chi(1)}, \dots, a_nb_{\chi(n)})\,,
\]
where $\chi= \theta (\alpha)$. An example is illustrated in Figure 2.3.

%%%%%%%%%%%
\begin{figure}[htb]
\centerline{
\setlength{\unitlength}{.4cm}
\begin{picture}(21,6)
\put(2,0){\includegraphics[width=7.6cm]{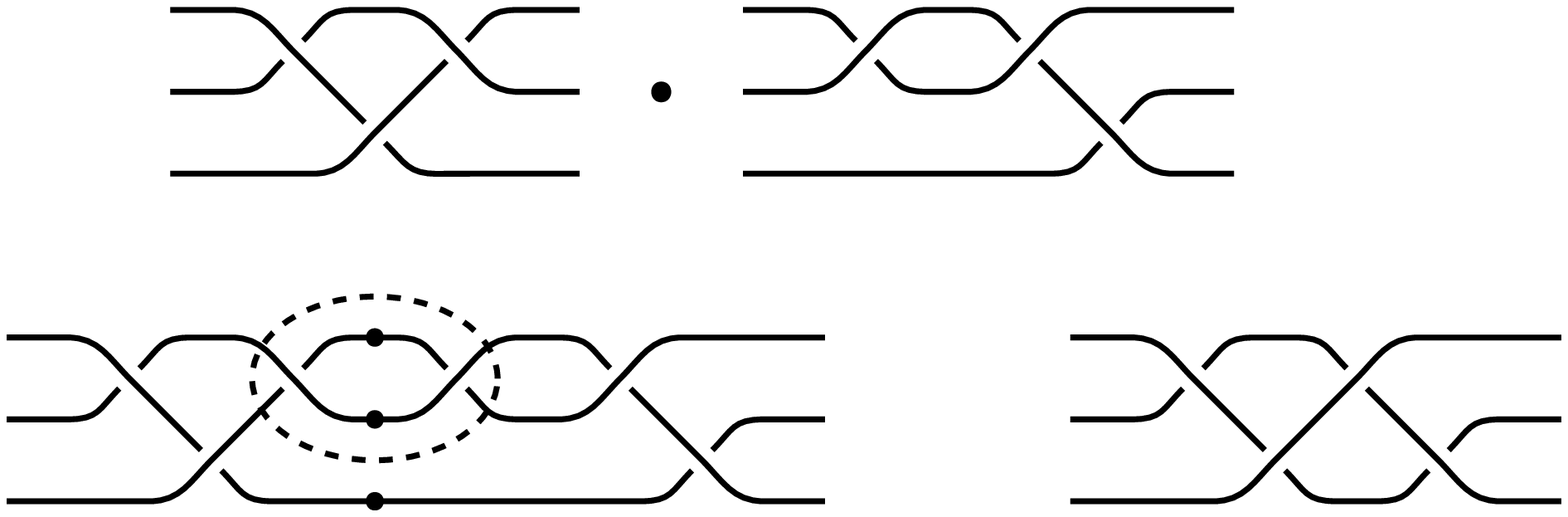}}
\put(0.5,1){$=$}
\put(13.2,1){$=$}
\end{picture}}
\medskip
\centerline{{\bf Figure 2.3.} Product of two braids.}
\end{figure}
%%%%%%%%%%

Let $\BB_n$ denote the set of homotopy classes of braids on $n$ strands. It is easily seen that the 
above defined multiplication of braids induces an operation on $\BB_n$. Moreover, we have the 
following.

\begin{proposition}
The set $\BB_n$ endowed with this operation is a group.
\end{proposition}

From now on, except mention of the contrary, by a braid we will mean a homotopy class of braids. The 
group $\BB_n$ of Proposition 2.1 is called the {\it braid group on $n$ strands}\index{Braid group}.  
The identity is the {\it constant braid}\index{Constant braid} $\Id= (\Id_1, \dots, 
\Id_n)$, where, for $1 \le k\le n$, $\Id_k$ denotes the constant path on $P_k$. The inverse of a braid 
$\beta$ is its {\it mirror} as illustrated in Figure 2.4.

%%%%%%%%%%%%%%%
\begin{figure}[htb]
\centerline{
\setlength{\unitlength}{.4cm}
\begin{picture}(16,3)
\put(0,1){\includegraphics[width=6.4cm]{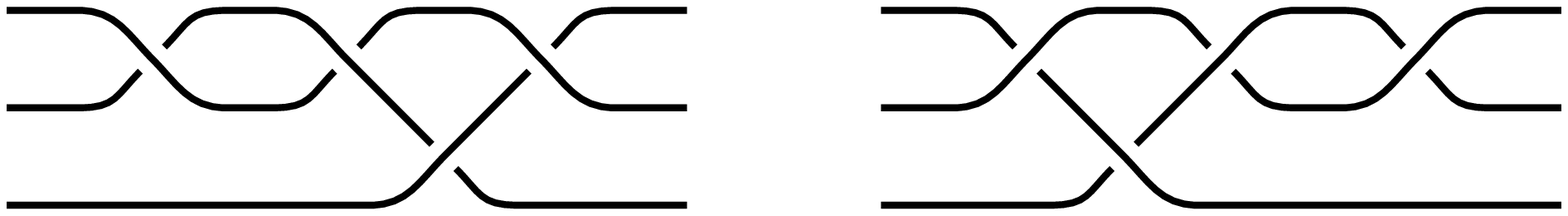}}
\put(3,0){\small $\beta$}
\put(12,0){\small $\beta^{-1}$}
\end{picture}}
\medskip
\centerline{{\bf Figure 2.4.} Inverse of a braid.}
\end{figure}
%%%%%%%%%%%%%%%

Recall that, if two braids $\alpha, \alpha'$ are homotopic, then $\theta( \alpha) = \theta(\alpha')$. 
Hence, the map $\theta$ from the set of braids on $n$ strands to $\Sym_n$ induces a map $\theta: \BB_n 
\to \Sym_n$. It is easily checked that this map is an epimorphism. Its kernel is called the {\it 
pure braid group on $n$ strands}\index{Pure braid group} and is denoted by $\PP\BB_n$. It plays an important role in the 
theory.

Let $\sigma_k$ be the braid illustrated in Figure 2.5. One can easily verify that $\sigma_1, \dots, 
\sigma_{n-1}$ generate the braid group $\BB_n$ and satisfy the relations 
\[\begin{array}{cr}
\sigma_k \sigma_l = \sigma_l \sigma_l &\quad \text{if } |k-l| \ge 2\,,\\
\sigma_k \sigma_l \sigma_k = \sigma_l \sigma_k \sigma_l &\quad \text{if } |k-l|=1\,.
\end{array}\]
(See Figure 2.6.) These relations suffice to define the braid group, namely:

%%%%%%%%%%%%%%%%%
\begin{figure}[htb]
\centerline{
\setlength{\unitlength}{.4cm}
\begin{picture}(6,7)
\put(2,0){\includegraphics[width=1.6cm]{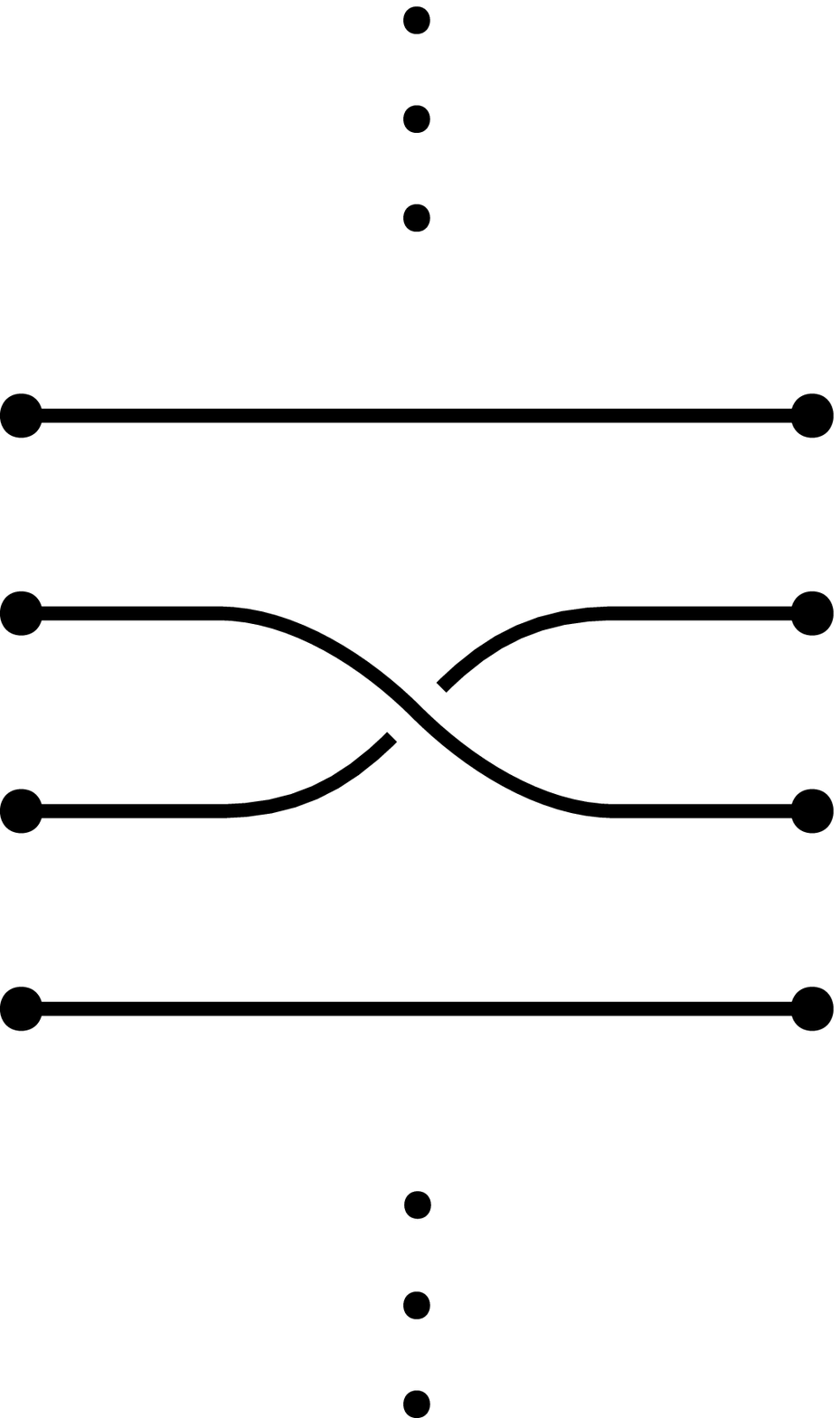}}
\put(0.6,2.7){\small $P_k$}
\put(0.1,3.7){\small $P_{k+1}$}
\end{picture}}
\medskip
\centerline{{\bf Figure 2.5.} The braid $\sigma_k$.}
\end{figure}
%%%%%%%%%%%%%%

%%%%%%%%%%%%
\begin{figure}[htb]
\centerline{
\setlength{\unitlength}{.4cm}
\begin{picture}(13,10)
\put(0,1){\includegraphics[width=5.2cm]{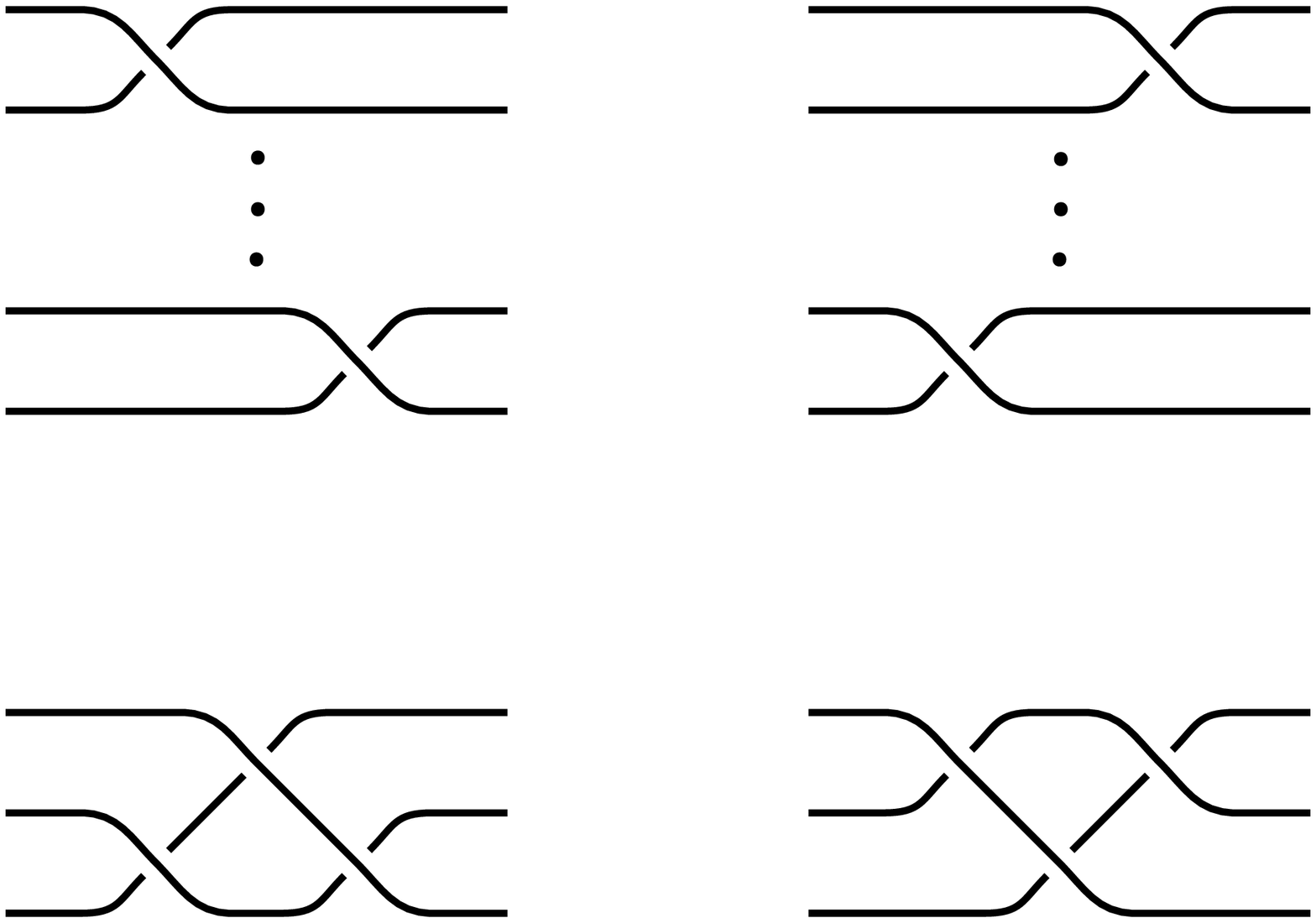}}
\put(1,0){\small $\sigma_k \sigma_{k+1} \sigma_k$}
\put(8.8,0){\small $\sigma_{k+1} \sigma_k \sigma_{k+1}$}
\put(1.3,5){\small $\sigma_k$}
\put(3.3,5){\small $\sigma_l$}
\put(9.3,5){\small $\sigma_l$}
\put(11.3,5){\small $\sigma_k$}
\end{picture}}
\medskip
\centerline{{\bf Figure 2.6.} Relations in $\BB_n$.}
\end{figure}
%%%%%%%%%

\begin{theorem}[Artin \cite{Artin1}, \cite{Artin2}, Magnus \cite{Magnu1}]
The group $\BB_n$ has 
a presentation with generators $\sigma_1, \dots, \sigma_{n-1}$ and relations
\[\begin{array}{cr}
\sigma_k \sigma_l = \sigma_l \sigma_k &\quad \text{if } |k-l| \ge 2\,,\\
\sigma_k \sigma_l \sigma_k = \sigma_l \sigma_k \sigma_l &\quad \text{if } |k-l|=1\,.
\end{array}\]
\end{theorem}

\begin{theorem}[Burau \cite{Burau1}, Markov \cite{Marko1}]
For $1 \le k<l\le n$, let
\[
\delta_{k\,l} = \sigma_{l-1} \cdots \sigma_{k+1} \sigma_k^2 \sigma_{k+1}^{-1} \cdots \sigma_{l-1}^{-1} \,.
\]
Then the pure braid group $\PP\BB_n$ has a presentation with generators
\[
\delta_{k\,l}\,, \quad 1 \le k<l\le n\,,
\]
and relations
\[\begin{array}{cl}
\delta_{r\,s} \delta_{k\,l} \delta_{r\,s}^{-1} = \delta_{k\,l} &\ \text{if } 1 \le r<s<k<l \le n \\
&\ \text{ or } 1 \le k<r<s<l \le n\,,\\
\delta_{r\,k} \delta_{k\,l} \delta_{r\,k}^{-1} = \delta_{k\,l}^{-1} \delta_{r\,l}^{-1} \delta_{k\,l} 
\delta_{r\,l} \delta_{k\,l} &\ \text{if } 1 \le r < k < l \le n\,,\\
\delta_{r\,k} \delta_{r\,l} \delta_{r\,k}^{-1} = \delta_{k\,l}^{-1} \delta_{r\,l} \delta_{k\,l} 
&\ \text{if } 1 \le r <k<l \le n\,,\\
\delta_{r\,s} \delta_{k\,l} \delta_{r\,s}^{-1} = \delta_{s\,l}^{-1} \delta_{r\,l}^{-1} \delta_{s\,l} 
\delta_{r\,l} \delta_{k\,l} \delta_{r\,l}^{-1} \delta_{s\,l}^{-1} \delta_{r\,l} \delta_{s\,l} &\ 
\text{if } 1 \le r<k<s<l\le n\,.
\end{array}\]\end{theorem}

\begin{note}
Most of the proofs of Theorems 2.2 and 2.3 that can be found in the literature proceed as follows. Given an 
exact sequence
\[
1 \to K \longrightarrow G \longrightarrow H \to 1\,,
\]
there is a machinery to compute a presentation of $G$ from presentations of $K$ and $H$. We start with 
the observation that $\PP\BB_2 \simeq \Z$ and with the exact sequence
\begin{equation}\label{R21}
1 \to F_n \longrightarrow \PP\BB_{n+1} \longrightarrow \PP\BB_n \to 1\,,
\end{equation}
where $F_n$ is a free group of rank $n$, to prove Theorem 2.3 by induction on $n$. (The exact sequence 
\eqref{R21} will be explained in Subsection 2.2.) Then we use the exact sequence
\[
1 \to \PP\BB_n \longrightarrow \BB_n \longrightarrow \Sym_n \to 1
\]
to prove Theorem 2.2 from Theorem 2.3. Another proof which, as far as I know, is not written in the 
literature but is known to experts, consists on extracting the presentation of Theorem 2.2 from the 
Salvetti complex of $\BB_n$. This is a cellular complex which is a $K(\BB_n,1)$ (see Section 5). 
\end{note}

\subsection{Configuration spaces}

We identify $\R^2$ with $\C$ and $P_k$ with $k \in \C$ for all $1 \le k\le n$. For $1 \le k<l\le n$ we 
denote by $H_{k\,l}$ the linear hyperplane of $\C^n$ defined by the equation $z_k=z_l$. The {\it big 
diagonal} of $\C^n$ is defined to be
\[
\Diag_n= \bigcup_{1 \le k<l\le n} H_{k\,l}\,.
\]
The {\it space of ordered configurations of $n$ points}\index{Configuration space} in $\C$ is defined to be
\[
M_n= \C^n \setminus \Diag_n\,.
\]
This is the space of $n$-tuples $\z= (z_1, \dots, z_n)$ of complex numbers such that $z_k \neq z_l$ for 
$k \neq l$. The symmetric group $\Sym_n$ acts freely on $M_n$. The quotient
\[
N_n= M_n/\Sym_n
\]
is called the {\it space of configurations of $n$ points}\index{Configuration space}
 in $\C$. This is the space of unordered $n$-tuples 
$\z = \{ z_1, \dots, z_n\}$ of complex numbers such that $z_k \neq z_l$ for $k \neq l$.

\begin{proposition}
Let $P_0= (1,2, \dots, n) \in M_n$. Then $\pi_1 (M_n, P_0) = \PP\BB_n$.
\end{proposition}

\begin{proof}
For a pure braid $\beta = (b_1, \dots, b_n)$ we set
\[\begin{array}{rccc}
\varphi(\beta): &[0,1] &\to& M_n\\
&t&\mapsto&(b_1(t), \dots, b_n(t))\,.
\end{array}\]Clearly, $\varphi(\beta)$ is a loop based at $P_0$. Moreover, two pure braids $\alpha$ and 
$\alpha'$ are homotopic if and only if $\varphi(\alpha)$ and $\varphi(\alpha')$ are homotopic. Thus 
$\varphi$ induces a bijection $\varphi_\ast : \PP\BB_n \to \pi_1 (M_n,P_0)$ which turns out to be a 
homomorphism.
\end{proof} 

For $\z \in M_n$, we denote by $[\z]$ the element of $N_n = M_n/ \Sym_n$ represented by $\z$.

\begin{proposition}
$\pi_1 (N_n, [P_0]) = \BB_n$.
\end{proposition}

\begin{proof}
For a braid $\beta = (b_1, \dots, b_n)$ we set
\[\begin{array}{rccc}
\hat\varphi (\beta): &[0,1] &\to&N_n\\
&t& \mapsto &[b_1(t), \dots, b_n(t)]\,.
\end{array}\]
Clearly, $\hat \varphi (\beta)$ is a loop based at $[P_0]$. It is easily checked that $\hat\varphi$ 
induces a homomorphism $\hat\varphi_\ast: \BB_n \to \pi_1(N_n, [P_0])$, and that the following diagram 
commutes
\[
\begin{CD}
1 @>>> \PP\BB_n @>>> \BB_n @>>> \Sym_n @>>> 1 \\
&& @V\varphi_\ast V\simeq V @VV\hat \varphi_\ast V @VV\Id V \\
1 @>>> \pi_1 (M_n, P_0) @>>> \pi_1(N_n, [ P_0]) @>>> \Sym_n @>>> 1
\end{CD}
\]
The first row is exact by definition, and the second one is associated to the regular covering $M_n 
\to N_n = M_n/\Sym_n$, so it is exact, too. We conclude by the five lemma that $\hat \varphi_\ast$ is 
an isomorphism.
\end{proof}

Let $f,g \in \C [x]$ be two non-constant polynomials. Set
\[\begin{array}{rcl}
f&=&a_0x^m +a_1 x^{m-1} + \cdots +a_m\,, \quad a_0 \neq 0\,,\\
g&=&b_0x^n+b_1x^{n-1} + \cdots + b_n\,, \quad b_0 \neq 0\,.
\end{array}\]
The {\it Sylvester matrix}\index{Sylvester matrix} of $f$ and $g$ is defined to be
\[
\Sylv (f,g) = \left(
\underbrace{
\begin{array}{cccc}
a_0&0&\cdots&0\\
a_1&a_0&\ddots&\vdots\\
\vdots&a_1&\ddots&0\\
a_m&\vdots&\ddots&a_0\\
0&a_m&&a_1\\
\vdots&\ddots&\ddots&\vdots\\
0&\cdots&0&a_m
\end{array}}_{n \text{ columns}}
\underbrace{
\begin{array}{cccc}
b_0&0&\cdots&0\\
b_1&b_0&\ddots&\vdots\\
\vdots&b_1&\ddots&0\\
b_n&\vdots&\ddots&b_0\\
0&b_n&&b_1\\
\vdots&\ddots&\ddots&\vdots\\
0&\cdots&0&b_n
\end{array}}_{m \text{ columns}}
\right)
\]
The {\it resultant}\index{Resultant} of $f$ and $g$ is defined to be
\[
\Res (f,g)= \det (\Sylv (f,g))\,.
\]
The following is classical in algebraic geometry (see \cite{CoLiOS1}, for example).

\begin{theorem}
Let $f,g \in \C[x]$ be two non-constant polynomials. Then $f$ and $g$ have a 
common root if and only if $\Res(f,g)=0$.
\end{theorem}

\begin{corollary}
Let $f \in \C[x]$ be a polynomial of degree $d \ge 2$. Then $f$ has a 
multiple root if and only if $\Res(f,f') =0$.
\end{corollary}

The number $\Res(f,f')$ is called the {\it discriminant}\index{Discriminant} of $f$ and is denoted by $\Disc(f)$. For 
instance, if $f=ax^2 +bx +c$, then $\Disc(f) = b^2-4ac$.

Let $n \ge 2$ and let $\C_n [x]$ be the set of monic polynomials of degree $n$. In particular, $\C_n[x]$ is 
isomorphic to $\C^n$. The map $\Disc: \C_n[x] \to \C$ is clearly a polynomial function, thus
\[
\DD = \{ f \in \C_n[x]; f \text{ has a multiple root}\} = \{ f\in \C_n[x]; \Disc (f)=0\}
\]
is an algebraic hypersurface called the {\it $n$-th discriminant}\index{Discriminant}. 
It is related to the braid group by 
the following.

\begin{proposition}
$N_n= \C_n[x] \setminus \DD$.
\end{proposition}

\begin{proof}
Let $\Phi: M_n \to \C_n [x] \setminus \DD$ be the map defined by
\[
\Phi (z_1, \dots, z_n) = (x-z_1) \cdots (x-z_n)\,.
\]
Then $\Phi$ is surjective and we have $\Phi(\u) = \Phi({\bf v})$ if and only if there exists $\chi \in 
\Sym_n$ such that ${\bf v} = \chi (\u)$. Thus $\C_n[x] \setminus \DD \simeq M_n/\Sym_n = N_n$.
\end{proof}

Now, recall the homotopy long exact sequence\index{Homotopy long exact sequence}
of a fiber bundle (see \cite{Hu1}, for example).

\begin{theorem}
Let $p: M \to B$ be a locally trivial fiber bundle. Let $b_0 \in B$, let 
$F=p^{-1}(b_0)$, and let $P_0 \in F$. Assume that $F$ is connected. Then there is a long exact sequence
of homotopy groups
\begin{multline*}
\cdots \to \pi_{k+1}(B,b_0) \to \pi_k (F,P_0) \to \pi_k(M,P_0) \to \pi_k(B,b_0) \to \cdots\\ 
\cdots\to \pi_2(B,b_0) \to \pi_1(F,P_0) \to \pi_1(M,P_0) \to \pi_1(B,b_0) \to 1\,.
\end{multline*}
\end{theorem}

There are two cases where this long exact sequence becomes a short exact sequence: when $\pi_2(B,b_0) = \{0\}$, 
and when $p$ admits a cross-section $\kappa: B \to M$. In the latter case the short exact sequence 
splits. It turns out that both situations hold in the study of $M_n$.

\begin{theorem}[Fadell, Neuwirth \cite{FadNeu1}]
Let $p: M_{n+1} \to M_n$ be defined by
\[
p(z_1, \dots, z_n, z_{n+1}) = (z_1, \dots, z_n)\,.
\]
Then $p$ is a locally trivial fiber bundle which admits a cross-section $\kappa: M_n \to M_{n+1}$.
\end{theorem}

Let $b_0 =(1,2, \dots, n)$. Then the fiber $p^{-1} (b_0)$ is naturally homeomorphic to $\C \setminus \{ 
1,2,\dots, n\}$ whose fundamental group is the free group $F_n$ of rank $n$. A cross-section of $p$ is 
the map $\kappa: M_n \to M_{n+1}$ defined by
\[
\kappa (z_1, \dots, z_n) = (z_1, \dots, z_n, |z_1| + \cdots + |z_n| +1)\,.
\]

\begin{corollary}
Let $n \ge 2$. Then there is a split exact sequence
\[
\begin{CD}
1 @>>> F_n @>>> \PP\BB_{n+1} @> p_\ast >> \PP\BB_n @>>> 1\,.\\
\noalign{\vskip-16pt}
 &&&& @<<\kappa_\ast< \\
\end{CD}
\]
\end{corollary}

A connected CW-complex $X$ is called $K(\pi,1)$\index{$\K(\pi,1)$ space} if its universal cover is contractible. Equivalently, 
$X$ is $K(\pi,1)$ if $\pi_k(X) = \{0\}$ for all $k \ge 2$. In particular, a space $X$ is $K(\pi,1)$ if 
an only if some of its connected cover $Y$ is $K(\pi,1)$. The notion of $K(\pi,1)$ spaces is of 
importance in the calculation of the (co)homology of groups. We refer to \cite{Brown1} for detailed 
explanations on the subject.

It is easily seen that $\C \setminus \{1, \dots, n\}$ is $K(\pi,1)$, thus, from Theorems 2.9 and 
2.10 follows:

\begin{corollary}
The spaces $M_n$ and $N_n$ are $K(\pi,1)$.
\end{corollary}

It is also known that the fundamental group of a finite dimensional $K(\pi,1)$ space is torsion free 
(see \cite{Brown1}), thus:

\begin{corollary}
$\BB_n = \pi_1(N_n)$ is torsion free.
\end{corollary}

\subsection{Mapping class groups}

Let $\Sigma$ be an oriented compact surface, possibly with boundary, and let $\PP = \{ P_1, \dots, 
P_n\}$ be a collection of $n$ punctures in the interior of $\Sigma$. Let 
\linebreak
$\Homeo^+ (\Sigma,\PP)$ denote 
the group of homeomorphisms $h: \Sigma \to \Sigma$ which preserve the orientation, which pointwise 
fix the boundary of $\Sigma$, and such that $h(\PP) = \PP$. Let $\Homeo^+_0 (\Sigma, \PP)$ denote 
the connected component of the identity in $\Homeo (\Sigma, \PP)$. 
The {\it mapping class group}\index{Mapping class group} of 
the pair $(\Sigma,\PP)$ is defined to be
\[
\MM (\Sigma,\PP) = \pi_0 (\Homeo^+(\Sigma,\PP)) = \Homeo^+(\Sigma,\PP)/\Homeo^+_0 (\Sigma,\PP)\,.
\]

A {\it braid}\index{Braid} of $\Sigma$ based at $\PP$ is defined to be a $n$-tuple $\beta = (b_1, \dots, b_n)$ of 
paths, $b_k: [0,1] \to \Sigma$, such that
\begin{itemize}
\item
$b_k(0)=P_k$ for all $1 \le k\le n$;
\item
there exists a permutation $\chi = \theta (\beta) \in \Sym_n$ such that $b_k(1) = P_{\chi(k)}$ for all 
$1 \le k\le n$;
\item
$b_k(t) \neq b_l(t)$ for all $k \neq l$ and all $t \in [0,1]$.
\end{itemize}
The homotopy classes of braids based at $\PP$ form a group denote by $\BB_n( \Sigma, \PP)$ and called 
the {\it braid group of $\Sigma$ on $n$ strands based at $\PP$}\index{Braid group}. It does not depend up to isomorphism on 
the choice of $\PP$ but only on the cardinality $n= |\PP|$. So, we may often write $\BB_n(\Sigma)$ in 
place of $\BB_n(\Sigma, \PP)$. If $\Sigma = \D$ is a disk, then $\BB_n(\Sigma)$ is naturally 
isomorphic to the braid group $\BB_n$.

For $1 \le k<l\le n$, we denote by $H_{k\,l}(\Sigma)$ the set of $n$-tuples $\x= (x_1, \dots, x_n) \in 
\Sigma^n$ such that $x_k=x_l$. The {\it big diagonal} of $\Sigma^n$ is defined to be
\[
\Diag_n (\Sigma) = \bigcup_{1 \le k<l\le n} H_{k\,l} (\Sigma)\,.
\]
The {\it space of ordered configurations of $n$ points in $\Sigma$}\index{Configuration space}
 is defined to be
\[
M_n(\Sigma) = \Sigma^n \setminus \Diag_n (\Sigma)\,.
\]
This is the space of $n$-tuples $\x = (x_1, \dots, x_n)$ in $\Sigma^n$ such that $x_k \neq x_l$ for all 
$1 \le k \neq l \le n$. The symmetric group $\Sym_n$ acts freely on $M_n (\Sigma)$, and the quotient
\[
N_n(\Sigma) = M_n(\Sigma) / \Sym_n
\]
is called the {\it space of configurations of $n$ points in $\Sigma$}\index{Configuration space}. 
This is the space of unordered 
$n$-tuples $\x = \{ x_1, \dots, x_n\}$ of elements of $\Sigma$ such that $x_k \neq x_l$ for all $1 \le 
k \neq l \le n$.

Set ${\bf P}_0=(P_1, \dots, P_n) \in M_n(\Sigma)$. For $\x \in M_n(\Sigma)$, we denote by $[\x]$ the element 
of $N_n(\Sigma)$ represented by $\x$. The following can be proved in the same way as Proposition 2.5.

\begin{proposition}
$\pi_1 (N_n(\Sigma), [{\bf P}_0]) \simeq \BB_n (\Sigma)$.
\end{proposition}

Now, the surface braid groups and the mapping class groups are related by the following exact sequence.

\begin{theorem}[Birman \cite{Birma2}]
Suppose $\Sigma$ is neither a sphere, nor a torus. Then 
we have the exact sequence
\[
1 \to \BB_n (\Sigma,\PP) \longrightarrow \MM (\Sigma, \PP) \longrightarrow \MM (\Sigma) \to 1 \,.
\]
\end{theorem}

\begin{note}
Let 
\[\begin{array}{rccc}
\Phi: &\Homeo^+ (\Sigma) &\to& N_n(\Sigma)\\
&\varphi &\mapsto &\{ \varphi(P_1), \dots, \varphi(P_n)\}\,.
\end{array}\]
Then $\Phi$ is a locally trivial fiber bundle, and the fiber of $\Phi$ over $\PP=[{\bf P}_0]$ is 
$\Homeo^+(\Sigma,\PP)$. Furthermore, it is known that $\pi_1 (\Homeo^+ (\Sigma)) = \{ 1\}$ (see 
\cite{Hamst1}), thus, by the homotopy long exact sequence of a fiber bundle (see \cite{Hu1}), we have 
the short exact sequence
\[
1\to \pi_1 (N_n (\Sigma), \PP) \longrightarrow \pi_0 (\Homeo^+ (\Sigma, \PP)) \longrightarrow \pi_0 
(\Homeo^+ (\Sigma)) \to 1\,,
\]
which is the same as the exact sequence of Theorem 2.15.
\end{note}

It is known that $\MM (\D) = \{1\}$ (see \cite{Alexa1}), thus, by Theorem 2.15:

\begin{theorem}[Artin \cite{Artin1}, \cite{Artin2}]
Let $\PP = \{ P_1, \dots, P_n\}$ be a 
collection of $n$ punctures in the interior of the disk $\D$. Then $\MM (\D, \PP) \simeq \BB_n$.
\end{theorem}

The isomorphism $\Phi: \MM (\D, \PP) \to \BB_n$ can be easily described as follows. Let $\varphi \in 
\Homeo^+ (\D,\PP)$. We know by \cite{Alexa1} that $\pi_0( \Homeo^+ (\D)) = \{1\}$, thus there exists a 
continuous path $\{ \varphi_t\}_{t \in [0,1]}$ in $\Homeo^+ (\D)$ such that $\varphi_0 = \Id$ and 
$\varphi_1 = \varphi$. Let $\beta = (b_1, \dots, b_n)$ be the braid defined by
\[
b_k(t) = \varphi_t (P_k)\,, \quad 1 \le k\le n \text{ and } t \in [0,1]\,.
\]
Then $\Phi (\varphi)$ is the homotopy class of $\beta$.

The reverse isomorphism $\Phi^{-1} : \BB_n \to \MM (\D,\PP)$ is more complicated to describe, but the 
images of the standard generators can be easily defined in terms of braid twists as follows.

We come back to the situation where $\Sigma$ is an oriented compact surface and $\PP = \{P_1, \dots,
P_n\}$ is a collection of $n$ punctures in the interior of $\Sigma$. Let $P_k,P_l \in \PP$, $k \neq l$. 
An {\it essential arc}\index{Essential arc} joining $P_k$ to $P_l$ is defined to be an embedding 
$a : [0,1] \to \Sigma$ such that $a(0) 
= P_k$, $a(1)= P_l$, $a((0,1)) \cap \PP = \emptyset$, and $a([0,1]) \cap \partial \Sigma = \emptyset$. 
Two essential arcs $a$ and $a'$ are said to be {\it isotopic} if there is a continuous family $\{a_t\}_{t \in 
[0,1]}$ of essential arcs such that $a_0=a$ and $a_1=a'$. Isotopy of essential arcs is an equivalence relation that we 
denote by $a \sim a'$.

Let $a$ be an essential arc joining $P_k$ to $P_l$. Let $\D = \{ z \in \C; |z| \le 1\}$ be the standard disk, and 
let $A : \D \to \Sigma$ be an embedding such that
\begin{itemize}
\item
$a(t) = A(t-\frac{1}{2})$ for all $t \in [0,1]$;
\item
$A( \D) \cap \PP = \{ P_k,P_l \}$.
\end{itemize}
Let $T \in \Homeo^+ (\Sigma, \PP)$ be defined by
\[
(T \circ A) (z) = A (e^{2i\pi |z|} z)\,, \quad z \in \D\,,
\]
and $T$ is the identity outside the image of A (see Figure 2.7). The {\it braid twist along $a$}\index{Braid twist} is 
defined to be the element $\tau_a \in \MM (\Sigma,\PP)$ represented by $T$, that is, the isotopy class 
of $T$. Note that:
\begin{itemize}
\item
the definition of $\tau_a$ does not depend on the choice of $A: \D \to \Sigma$;
\item
if $a$ is isotopic to $a'$, then $\tau_a = \tau_{a'}$.
\end{itemize}

%%%%%%%%%%%%%%%%%%%ù
\begin{figure}[htb]
\centerline{
\setlength{\unitlength}{.4cm}
\begin{picture}(22,8)
\put(0,0){\includegraphics[width=8.8cm]{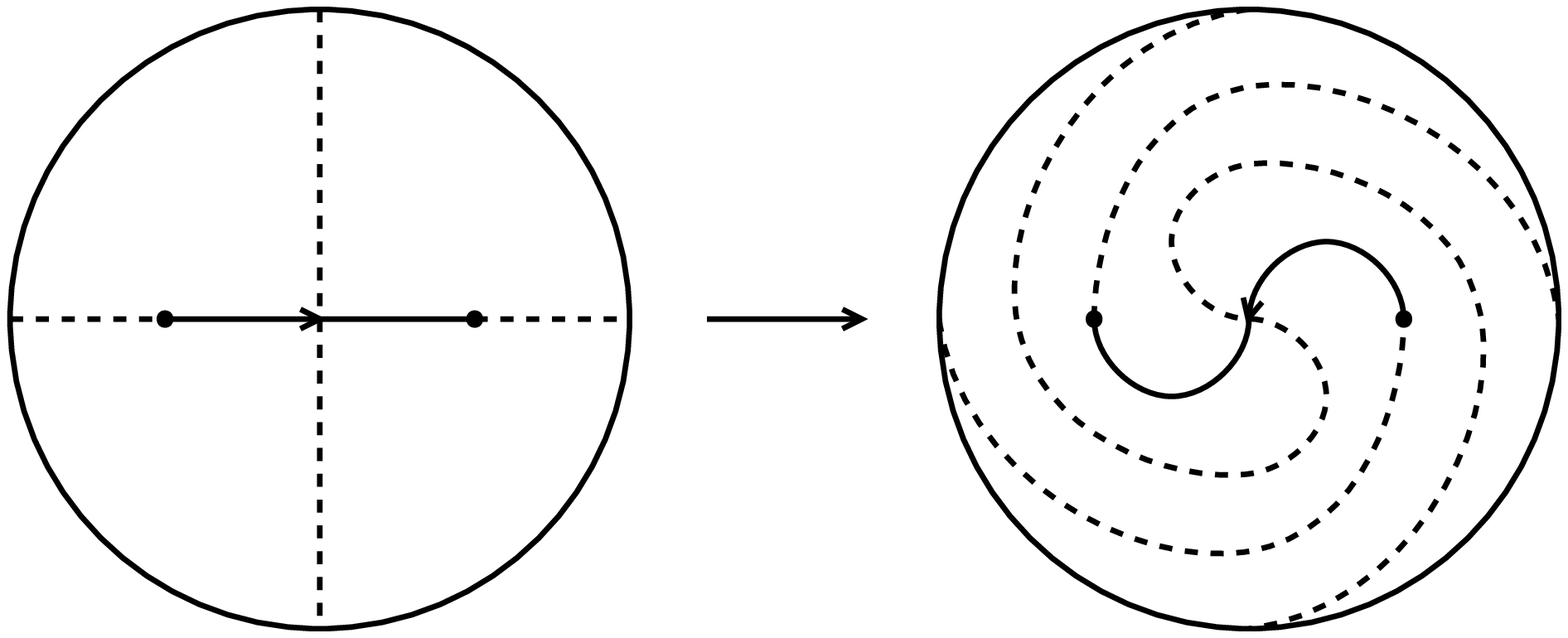}}
\put(2.2,5){\small $P_k$}
\put(6.5,5){\small $P_l$}
\put(4.8,4.8){\small $a$}
\put(10.8,4.8){\small $T$}
\end{picture}}
\medskip
\centerline{{\bf Figure 2.7.} Braid twist.}
\end{figure}
%%%%%%%%%%%%%%%%%%%%

Now, we view the disk $\D$ as the disk in $\C$ of radius $\frac{n+1}{2}$ centered at 
$\frac{n+1}{2}$, and we set $P_k=k$ for $1 \le k\le n$. Let $a_k: [0,1] \to \D$ be the arc defined by
\[
a_k(t)= k+t\,, \quad t \in [0,1]\,.
\]
(See Figure 2.8.) Then:

%%%%%%%%%%%%%
\begin{figure}[htb]
\centerline{
\setlength{\unitlength}{.4cm}
\begin{picture}(14,8)
\put(0,0){\includegraphics[width=5.6cm]{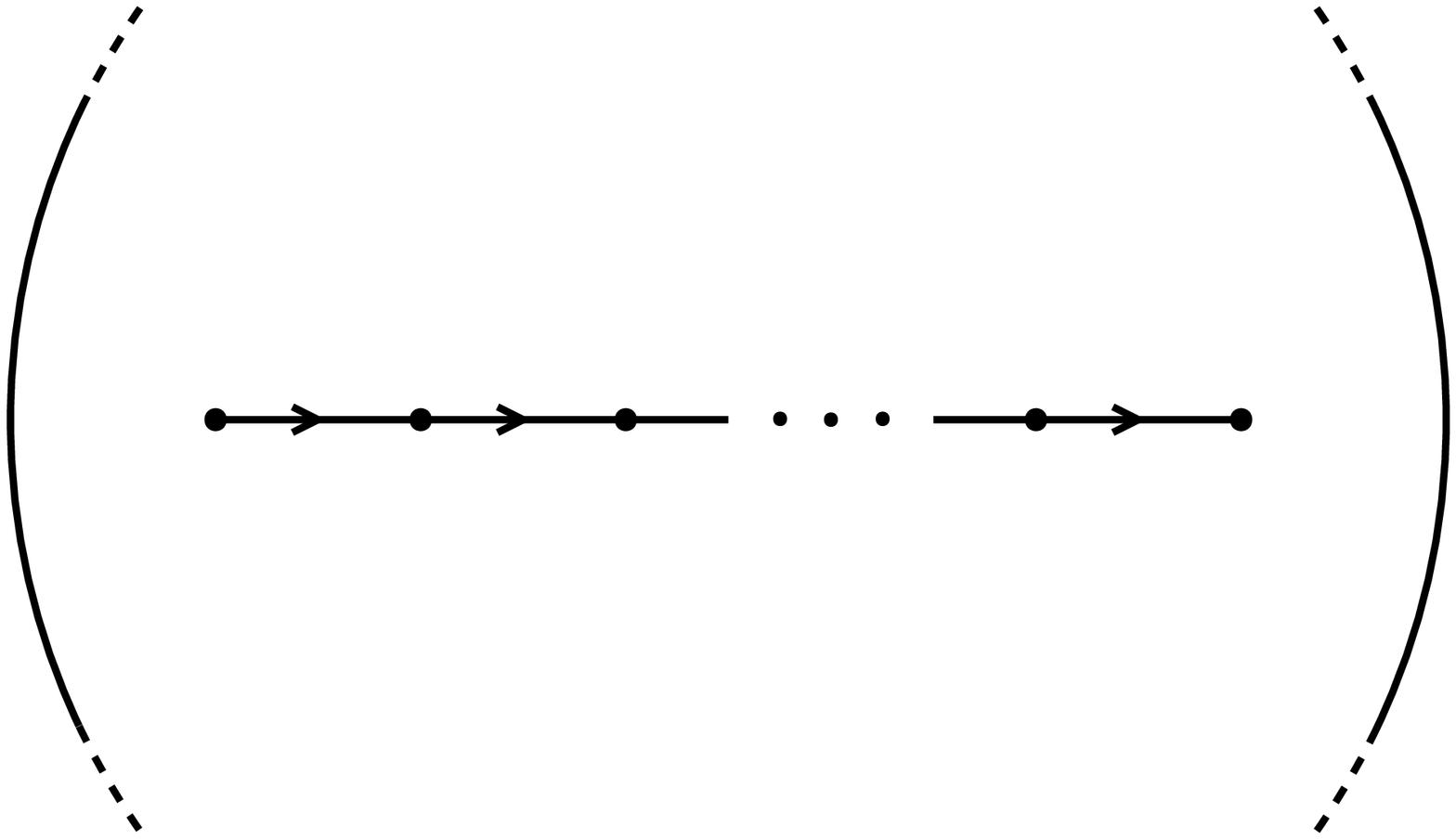}}
\put(1.8,3.2){\footnotesize $1$}
\put(3.8,3.2){\footnotesize $2$}
\put(5.8,3.2){\footnotesize $3$}
\put(11.8,3.2){\footnotesize $n$}
\put(2.5,4.5){\small $a_1$}
\put(4.5,4.5){\small $a_2$}
\put(10.5,4.5){\small $a_{n-1}$}
\end{picture}}
\medskip
\centerline{{\bf Figure 2.8.} The standard generators of $\MM (\D, \PP) = \BB_n$.}
\end{figure}
%%%%%%%%%%%ù

\begin{lemma}
The reverse isomorphism $\Phi^{-1} : \BB_n \to \MM (\D, \PP)$ is defined by
\[
\Phi^{-1} (\sigma_k) = \tau_{a_k}\,, \quad 1 \le k\le n-1\,.
\]
\end{lemma}

\subsection{Automorphisms of free groups}

For a group $G$, we denote by $\Aut(G)$ the group of automorphisms of $G$, by $\Inn (G)$ the group of 
inner automorphisms of $G$, and by $\Out(G) = 
\linebreak
\Aut(G) / \Inn (G)$ the group of outer automorphisms of 
$G$.

Let $F_n= F(x_1, \dots, x_n)$ be the free group of rank $n$. For $1 \le k\le n-1$, let $\tau_k: F_n \to 
F_n$ be the automorphism defined by
\[
\tau_k: \left\{ \begin{array}{cccl}
x_k&\mapsto&x_k^{-1} x_{k+1} x_k\\
x_{k+1}& \mapsto& x_k\\
x_l& \mapsto& x_l&\quad \text{if } l \neq k,k+1
\end{array}\right.
\]
One can easily show the following.

\begin{proposition}
The mapping $\sigma_k \mapsto \tau_k$, $1 \le k\le n-1$, determines a 
representation $\rho: \BB_n \to \Aut (F_n)$.
\end{proposition}

The above representation $\rho: \BB_n \to \Aut(F_n)$ is called the 
{\it Artin representation}\index{Artin representation}. It is 
faithful, more precisely:

\begin{theorem}[Artin \cite{Artin1}, \cite{Artin2}]
\begin{enumerate}
\item
The Artin representation $\rho: \BB_n \to 
\linebreak
\Aut(F_n)$ is faithful.
\item
An automorphism $\alpha \in \Aut(F_n)$ belongs to $\Im\, \rho$ if and only if
\linebreak 
$\alpha (x_n \cdots x_2 
x_1) = x_n \cdots x_2 x_1$ and
there exists a permutation 
$\chi \in \Sym_n$ such that $\alpha (x_k)$ is conjugate to $x_{\chi(k)}$ for all $1 \le k\le n$.
\end{enumerate}
\end{theorem}

In particular, $\BB_n$ can be viewed as a subgroup of $\Aut (F_n)$. This has some consequences on 
$\BB_n$ itself such as the two properties defined below.

A group $G$ is called {\it residually finite}\index{Residually finite group} 
if for all $g \in G \setminus \{ 1\}$ there exists a 
homomorphism $\varphi: G \to H$ such that $H$ is finite and $\varphi(g) \neq 1$. A group $G$ is called 
{\it Hopfian}\index{Hopfian group} if every epimorphism $\varphi: G \to G$ is an isomorphism. 
It is known that the subgroups 
of $\Aut(F_n)$ are both, residually finite and Hopfian (see \cite{MaKaSo1}), thus, by Theorem 2.19:

\begin{corollary}
The braid group $\BB_n$ is residually finite and Hopfian.
\end{corollary}

There are several ways to describe geometrically the Artin representation.
The first way is using the Fadell-Neuwirth fiber bundle $p: M_{n+1} \to M_n$ of Theorem 2.10. Let 
$\Sym_n$ act on $M_n$ and on $M_{n+1}$. The second action is on the first $n$ coordinates, that is,
\[
\chi (z_1, \dots, z_n,z_{n+1}) = (z_{\chi^{-1} (1)}, \dots, , z_{\chi^{-1}(n)}, z_{n+1})\,, \quad \text{for } \chi 
\in \Sym_n\,.
\]
The map $p: M_{n+1} \to M_n$ induces a map $\bar p: M_{n+1} / \Sym_n \to M_n/ \Sym_n = N_n$ which turns out 
to be a locally trivial fiber bundle. The fiber is again homeomorphic to $\C \setminus \{1,2, \dots, 
n\}$, and $\bar p : M_{n+1}/ \Sym_n \to N_n$ has also a cross-section $\bar \kappa : N_n \to 
M_{n+1}/\Sym_n$. So, from the homotopy long exact sequence of a fiber bundle (see Theorem 2.9) we 
obtain the following split exact sequence
\begin{multline}\label{R22}
1  \to F_n = \pi_1(\C \setminus \{1, \dots, n\}) \longrightarrow \pi_1 (M_{n+1} /\Sym_n)\\
\begin{array}{c}
{\scriptstyle \bar p_\ast}\\
\noalign{\vskip-6pt}
\longrightarrow\\
\noalign{\vskip-8pt}
\longleftarrow\\
\noalign{\vskip-6pt}
{\scriptstyle\bar\kappa_\ast}
\end{array}
\pi_1 (N_n)= \BB_n \to 1\,.
\end{multline}
The action of $\BB_n = \pi_1 (N_n)$ on $F_n = \pi_1 (\C \setminus \{1, \dots, n\})$ derived from the 
above split exact sequence is exactly the Artin representation.

Another way to represent the Artin representation is using the isomorphism $\BB_n \simeq 
\MM (\D, 
\{ P_1, \dots, P_n\})$. Fix a basepoint $P_0 \in \partial \D$. Then it is easily shown that $\MM (\D, \{ 
P_1, \dots, P_n \})$ acts on $\pi_1( \D \setminus \{P_1, \dots, P_n\}, P_0) = F_n$, and that this action is the 
Artin representation.

The latter point of view of the Artin representations can be extended to all the mapping class groups. In 
this setting, it is known as the Dehn-Nielsen-Baer theorem.
Here is a version of this theorem.

\begin{theorem}[Dehn, Nielsen \cite{Niels1}, Baer \cite{Baer1}, Magnus \cite{Magnu1}]
Let 
$\Sigma$ be a closed oriented surface, and let $\PP = \{ P_1, \dots, P_n\}$ be a collection of $n$ 
punctures in $\Sigma$. Then the natural homomorphism $\rho: \MM (\Sigma, \PP) \to \Out (\pi_1 (\Sigma 
\setminus \PP))$ is injective. Moreover, if $\PP = \emptyset$, then the image of $\rho$ is an index 
2 subgroup of $\Out (\pi_1 (\Sigma))$.
\end{theorem}

We refer to 
\cite{Ivano1} for a detailed exposition on the Dehn-Nielsen-Baer theorem which include other versions of it.

\begin{note}
There are some variants of the Artin representations introduced in 
\cite{Wada1} and  \cite{CriPar1} that lead to invariants of links.
\end{note}

%%%%%%%%%%%%%%%%%%%%%%%%%%%%%%%%%

\section{Artin groups}

\subsection{Definitions and examples}

Let $S$ be a finite set. A {\it Coxeter matrix}\index{Coxeter matrix} over $S$ is a square matrix $M=(m_{s\,t})_{s,t \in S}$ 
indexed by the elements of $S$ such that
\begin{itemize}
\item
$m_{s\,s} =1$ for all $s \in S$;
\item
$m_{s\,t} = m_{t\,s} \in \{ 2,3,4, \dots, +\infty\}$ for all $s,t \in S$, $s \neq t$.
\end{itemize}
A Coxeter matrix $M=(m_{s\,t})_{s,t\in S}$ is usually represented by its {\it Coxeter graph}\index{Coxeter graph}, $\Gamma = 
\Gamma (M)$. This is a labeled graph defined by the following data.
\begin{itemize}
\item
$S$ is the set of vertices of $\Gamma$.
\item
Two vertices $s,t \in S$, $s \neq t$, are joined by an edge if $m_{s\,t} \ge 3$. This edge is labeled by 
$m_{s\,t}$ if $m_{s\,t} \ge 4$.
\end{itemize}
Let $\Gamma$ be a Coxeter graph. Define the {\it Coxeter system of type $\Gamma$}\index{Coxeter system}  to be the pair 
$(W,S)$, where $W=W_\Gamma$ is the group presented by the generating set $S$ and the relations
\[
\begin{array}{cl}
s^2=1&\quad\text{for all } s\in S\,,\\
(st)^{m_{s\,t}}=1&\quad\text{for all } s,t \in S,\ s\neq t,\text{ and } m_{s\,t} \neq +\infty\,,
\end{array}\]
where $M=(m_{s\,t})_{s,t\in S}$ is the Coxeter matrix of $\Gamma$. The group $W=W_\Gamma$ is called the 
{\it Coxeter group}\index{Coxeter group} of type $\Gamma$.

If $a,b$ are two letters and $m \in \N$, then $\pprod (a,b:m)$ denotes the word
\[
\pprod (a,b:m) = \left\{
\begin{array}{ll}
(ab)^{\frac{m}{2}}&\quad\text{if }m \text{ is even}\,,\\
(ab)^{\frac{m-1}{2}}a&\quad\text{if }m \text{ is odd}\,.
\end{array}\right.
\]
Let $\Sigma= \{ \sigma_s; s \in S\}$ be an abstract set in one-to-one correspondence with $S$. Define 
the {\it Artin system of type $\Gamma$}\index{Artin system} to be the pair $(G,\Sigma)$, where $G=G_\Gamma$ is the group 
presented by the generating set $\Sigma$ and the relations
\[
\pprod (\sigma_s, \sigma_t:m_{s\,t}) = \pprod (\sigma_t, \sigma_s: m_{s\,t}) \quad \text{for }s,t \in 
S,\ s\neq t, \text{ and } m_{s\,t} \neq +\infty\,.
\]
The group $G$ is called the {\it Artin group of type $\Gamma$}\index{Artin group}.

It is easily checked that the group $W_\Gamma$ is also presented by the generating set $S$ and the 
relations
\[\begin{array}{cl}
s^2=1&\ \text{for all } s \in S\,,\\
\pprod(s,t: m_{s\,t}) = \pprod (t,s: m_{s\,t})& \ \text{for all } s,t \in S,\ s \neq t\, \text{ and 
} m_{s\,t} \neq +\infty\,.
\end{array}\]
This shows that the mapping $\Sigma \to S$, $\sigma_s \mapsto s$, induces a {\it canonical 
epimorphism} $\theta: G_\Gamma \to W_\Gamma$.

If $m_{s\,t}=2$, then
\[
\sigma_s \sigma_t = \pprod (\sigma_s, \sigma_t: m_{s\,t}) = \pprod (\sigma_t, \sigma_s: m_{s\,t}) = 
\sigma_t \sigma_s\,,
\]
that is, $\sigma_s$ and $\sigma_t$ commute. So, if $\Gamma_1, \dots, \Gamma_l$ are the connected 
components of $\Gamma$, then
\[
G_\Gamma= G_{\Gamma_1} \times G_{\Gamma_2} \times \cdots \times G_{\Gamma_l}\,.
\]
Similarly, we have
\[
W_\Gamma= W_{\Gamma_1} \times W_{\Gamma_2} \times \cdots \times W_{\Gamma_l}\,.
\]
We say that $G_\Gamma$ (or $W_\Gamma$) is {\it irreducible}\index{Irreducible Coxeter group}
\index{Irreducible Artin group}  if $\Gamma$ is connected. We say that 
$\Gamma$ (or $G_\Gamma$) is {\it of spherical type}\index{Coxeter graph of spherical type}
\index{Artin group of spherical type} if $W_\Gamma$ is finite.

\bigskip\noindent
{\bf Example 1.} Suppose that $\Gamma$ is the graph $A_n$ of Figure 3.1. Then $W_\Gamma= \Sym_{n+1}$ is 
the symmetric group of $\{1, \dots, n,n+1\}$, and the Coxeter generators are the transpositions 
$s_1=(1,2), s_2=(2,3), \dots, s_n=(n,n+1)$. The Artin group $G_\Gamma$ is the braid group $\BB_{n+1}$ 
on $n+1$ strands, and the Artin generators are the standard generators of $\BB_{n+1}$ given in Theorem 
2.2. The canonical epimorphism coincides with the epimorphism described in Subsection 2.1.

%%%%%%%%%%%%%%%%%%%%%%%%%%%%%%%%%%%%%%%%%%%%%%%%%
\begin{figure}[htb]
\centerline{
\setlength{\unitlength}{.4cm}
\begin{picture}(18,16)
\put(3,1){\includegraphics[width=4cm]{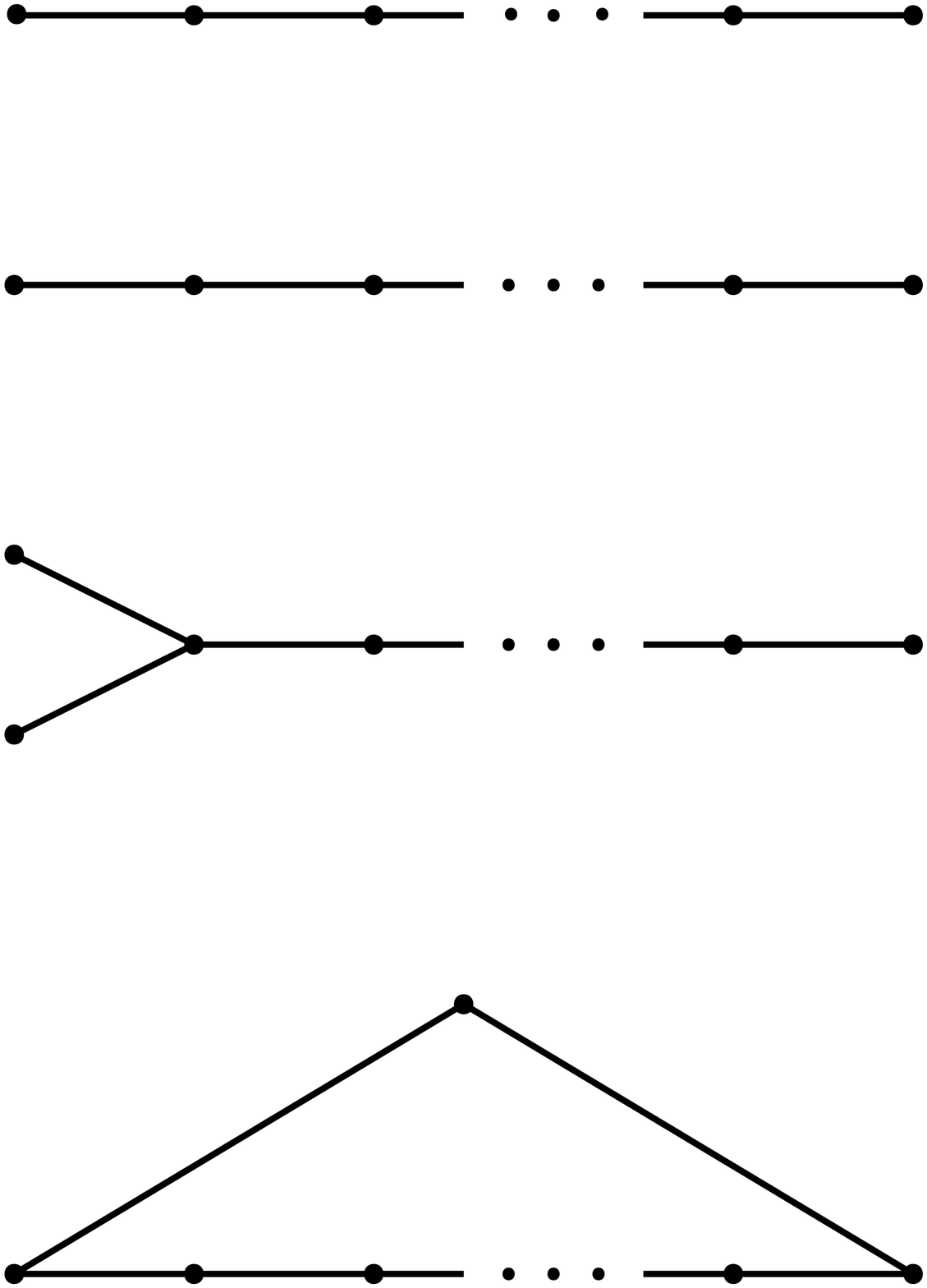}}
\put(3,14){\footnotesize $1$}
\put(5,14){\footnotesize $2$}
\put(6.9,14){\footnotesize $3$}
\put(12.8,14){\footnotesize $n$}
\put(0,14.5){\small $A_n$}
\put(16,14.5){\small $n \ge 1$}
\put(3,11){\footnotesize $1$}
\put(5,11){\footnotesize $2$}
\put(6.9,11){\footnotesize $3$}
\put(12.8,11){\footnotesize $n$}
\put(0,11.5){\small $B_n$}
\put(3.8,12.2){\small $4$}
\put(16,11.5){\small $n \ge 2$}
\put(3,8){\footnotesize $1$}
\put(3,6){\footnotesize $2$}
\put(5,7){\footnotesize $3$}
\put(6.9,7){\footnotesize $4$}
\put(12.8,7){\footnotesize $n$}
\put(0,7.5){\small $D_n$}
\put(16,7.5){\small $n \ge 4$}
\put(3,0){\footnotesize $1$}
\put(5,0){\footnotesize $2$}
\put(6.9,0){\footnotesize $3$}
\put(12.8,0){\footnotesize $n$}
\put(8.5,4){\footnotesize $n+1$}
\put(0,1){\small $\tilde A_n$}
\put(16,1){\small $n \ge 2$}
\end{picture}}
\medskip
\centerline{{\bf Figure 3.1.} The Coxeter graphs $A_n$, $B_n$, $D_n$, and $\tilde A_n$.}
\end{figure}
%%%%%%%%%%%%%%%

\bigskip\noindent
{\bf Example 2.} Suppose that $\Gamma$ is the Coxeter graph $B_n$ of Figure 3.1. Let $C_2 = \{ \pm 
1\}$ denote the cyclic group of order 2. Set $\Cub_n = C_2^n \rtimes \Sym_n$, where $\Sym_n$ acts on 
$C_2^n$ by permutation of the coordinates. This is the group of isometries of a regular $n$-cube (see 
\cite{Harpe1}, for example). The group $\Cub_n$ is the Coxeter group of type $B_n$, and the Coxeter 
generators are
\[
s_1=(-1,1, \dots, 1) \in C_2^n,\quad s_i=(i-1,i)\in \Sym_n\text{ for } 2 \le i\le n\,.
\]
Recall the Artin representation $\rho: \BB_n \to \Aut(F_n)$ defined in Subsection 2.4. Set $G= F_n 
\rtimes_\rho \BB_n$. Recall also the action of $\Sym_n$ on $M_{n+1}$ defined in Subsection 2.4. It 
follows from the exact sequence \eqref{R22} that $G= \pi_1(M_{n+1} /\Sym_n)$. In particular, $G$ is an 
index $n+1$ subgroup of $\pi_1(M_{n+1}/ \Sym_{n+1}) = \pi_1 (N_{n+1}) = \BB_{n+1} = G_{A_n}$. Now, $G$ 
is the Artin group of type $B_n$, and the Artin generators are
\[
\tau_1 = x_1 \in F_n, \quad \tau_i = \sigma_{i-1} \in \BB_{n} \text{ for } 2 \le i\le n\,.
\]
(See \cite{CriPar2}).

\bigskip\noindent
{\bf Example 3.} Suppose that $\Gamma$ is the Coxeter graph $D_n$ of Figure 3.1, where $n \ge 4$. Let 
$\sgn: C_2^n \to C_2$ be the homomorphism defined by
\[
\sgn (\varepsilon_1, \dots, \varepsilon_n) = \prod_{i=1}^n \varepsilon_i\,,
\]
and let $K$ be the kernel of $\sgn$. The subgroup $K$ is invariant under the action of $\Sym_n$, thus 
one can consider the subgroup $W= K \rtimes \Sym_n$ of $\Cub_n = C_2^n \rtimes \Sym_n$. This is the 
Coxeter group of type $D_n$, and the Coxeter generators are
\[
s_1 = (-1,-1,1, \dots, 1) \cdot (1,2), \quad s_i = (1,1,1, \dots, 1) \cdot (i-1,i) \text{ for }2 \le 
i\le n\,.
\]
(See \cite{Harpe1}, for example).

Let $F_{n-1} = F(y_1, \dots, y_{n-1})$ be a free group of rank $n-1$. Let $\rho_{D,1} : F_{n-1} \to 
F_{n-1}$ be the automorphism defined by
\[
\rho_{D,1}: \left\{\begin{array}{rcll}
y_1 &\mapsto& y_1\,,\\
y_j &\mapsto& y_1^{-1} y_j &\quad\text{if } j\ge 2\,.
\end{array}\right.
\]
For $2 \le i\le n-1$, let $\rho_{D,i}: F_{n-1} \to F_{n-1}$ be the automorphism defined by
\[
\rho_{D,i} : \left\{\begin{array}{rcll}
y_{i-1} &\mapsto& y_i\,,\\
y_i &\mapsto& y_i y_{i-1}^{-1} y_i\,,\\
y_j &\mapsto& y_j &\quad\text{if } j \neq i-1, i\,.
\end{array}\right.
\]
One can easily show the following.

\begin{lemma}
The mapping $\sigma_i \mapsto \rho_{D,i}$, $1 \le i\le n-1$, determines a 
representation $\rho_D: \BB_n \to \Aut (F_{n-1})$.
\end{lemma}

The following is implicit in \cite{PerVan1} and explicit in \cite{CriPar2}.

\begin{theorem}[Perron, Vannier \cite{PerVan1}]
The representation $\rho_D: \BB_n \to 
\linebreak
\Aut(F_{n-1})$ is faithful, and the semidirect product $F_{n-1} \rtimes_{\rho_D} \BB_n$ is isomorphic to the 
Artin group $G_{D_n}$ of type $D_n$.
\end{theorem}

\begin{note}
It was shown by Allcock \cite{Allco1} that the Artin group $G_{D_n}$ of type $D_n$ can 
be also presented as an index 2 subgroup of the $n$-strand braid group of a plane with a single 
orbifold point of degree 2.
\end{note}

\noindent
{\bf Example 4.} Suppose that $\Gamma$ is the graph $\tilde A_n$ of Figure 3.1. Let $\Sym_{n+1}$ act on 
$\Z^{n+1}$ by permutations of the coordinates. Then $\Z^{n+1} \rtimes \Sym_{n+1}$ is the Coxeter group 
of type $\Gamma$ (see \cite{Bourb1}).

Let $\Phi: G_{B_{n+1}} \to \Z$ be the homomorphism defined by
\[
\Phi( \sigma_1) = 1\,, \quad \Phi( \sigma_i) = 0 \text{ for } 2 \le i\le n\,.
\]
It was observed by several authors \cite{Allco1}, \cite{ChaPei1}, \cite{Dieck1}, \cite{KenPei1}, 
that the kernel of $\Phi$ is isomorphic to the Artin group $G_{\tilde A_n}$ of type $\tilde A_n$. In 
particular, $G_{\tilde A_n}$ is a subgroup of $\BB_{n+2}$.

The Artin generators of $G_{\tilde A_n}$, viewed as a subgroup of $\BB_{n+2} = 
\linebreak
\MM (\D, \{P_1,P_2, 
\dots, P_{n+2}\})$, can be described in terms of braid twists as follows. We place $P_1, \dots, P_{n+2}$ 
in the interior of $\D$ like in Figure 3.2. For $1 \le i\le n+1$, let $\tau_i$ denote the braid twist 
along the arc $a_i$. Then $\tau_1, \dots, \tau_{n+1}$ are the Artin generators of $G_{\tilde A_n}$.

%%%%%%%%%%%%%%%
\begin{figure}[htb]
\centerline{
\setlength{\unitlength}{.4cm}
\begin{picture}(12,12)
\put(0,0){\includegraphics[width=4.8cm]{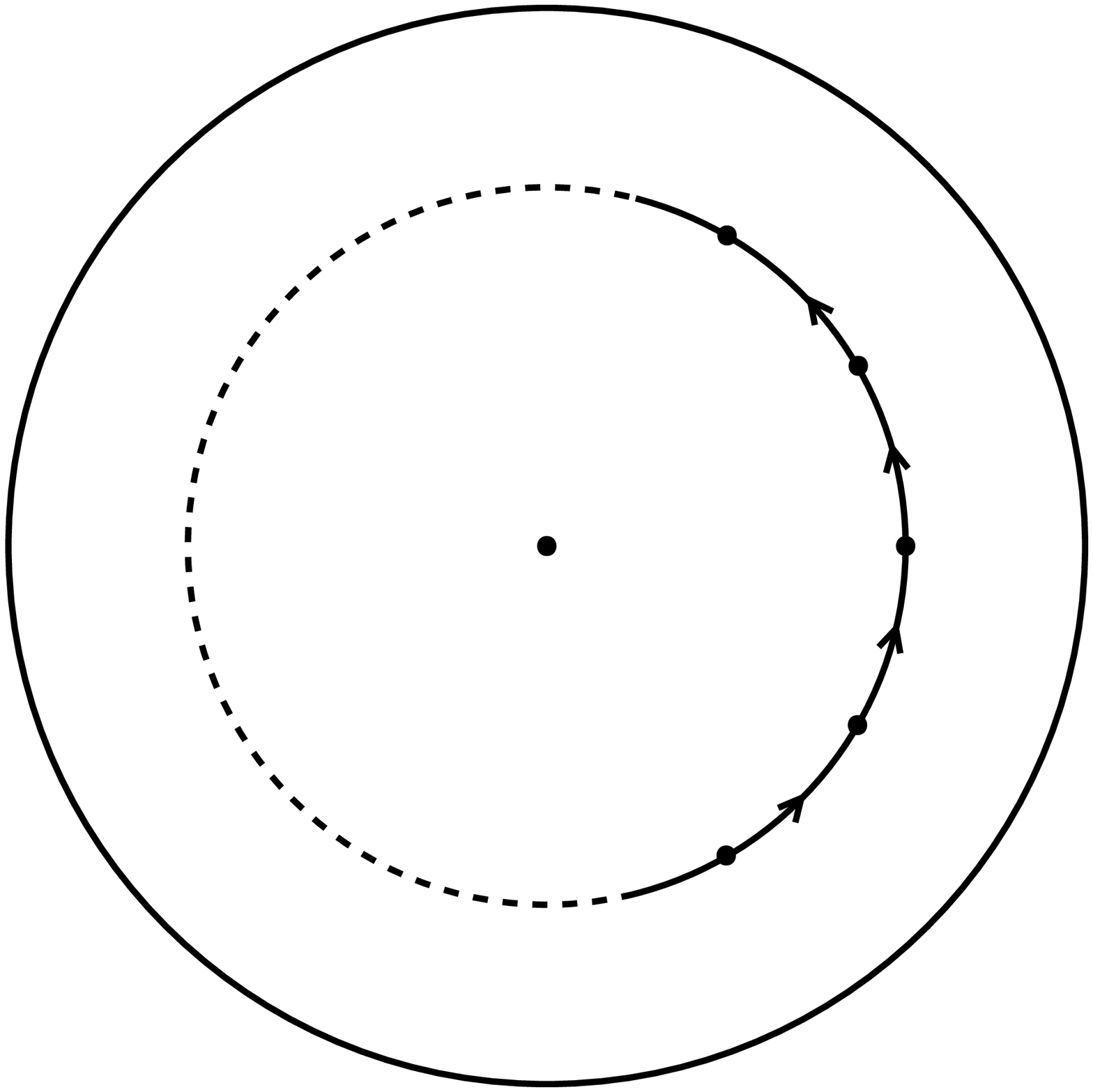}}
\put(10.3,6){\footnotesize $P_1$}
\put(9.7,7.8){\footnotesize $P_2$}
\put(8.1,9.8){\footnotesize $P_3$}
\put(8,1.7){\footnotesize $P_n$}
\put(9.8,3.5){\footnotesize $P_{n+1}$}
\put(5.4,5.4){\footnotesize $P_{n+2}$}
\put(9,6.5){\footnotesize $a_1$}
\put(8,8.5){\footnotesize $a_2$}
\put(8,3.5){\footnotesize $a_n$}
\put(8.3,5){\footnotesize $a_{n+1}$}
\end{picture}}
\medskip
\centerline{{\bf Figure 3.2.} Standard generators of $G_{\tilde A_n}$.}
\end{figure}
%%%%%%%%%%%%%%%%
 
\begin{note}
For a group $G$ we denote by $Z(G)$ the center of $G$. If $\Gamma = A_n, B_n$, or $\tilde 
A_n$, then $G_\Gamma/ Z(G_\Gamma)$ can be viewed as a finite index subgroup of the mapping class group 
of a punctured sphere. This has been cleverly exploited to study the group $G_\Gamma$ itself, in 
particular, to compute the group of automorphisms of $G_\Gamma$ (see \cite{ChaCri1}, \cite{BelMar1}). 
Note that the center of $G_{A_n}$ and $G_{B_n}$ is an infinite cyclic group (see \cite{Delig1}, 
\cite{BriSai1}), and the center of $G_{\tilde A_n}$ is trivial (see \cite{JohAlb1}).
\end{note}

\subsection{Coxeter groups}

The Coxeter groups were introduced by Tits \cite{Tits1} in a manuscript which was recently published, 
and whose results appeared in the seminal Bourbaki's book \cite{Bourb1}. The present 
subsection is a brief survey on these groups with a special emphasis on the results that are needed to study 
the Artin groups. Standard references for the subject are \cite{Bourb1}, \cite{Humph1}.

Let $\Gamma$ be a Coxeter graph, let $M= (m_{s\,t})_{s,t\in S}$ be its associated Coxeter matrix, and 
let $(W,S)$ be the Coxeter system of type $\Gamma$.

Let $\Pi= \{ e_s; s \in S\}$ be an abstract set in one-to-one correspondence with $S$, whose elements are 
called {\it simple roots}\index{Simple root}. We denote by $V$ the real vector space having $\Pi$ as a basis, and by 
$\langle, \rangle: V \times V \to \R$ the symmetric bilinear form defined by
\[
\langle e_s, e_t \rangle = \left\{ \begin{array}{cl}
- \cos (\frac{\pi}{m_{s\,t}}) &\quad \text{if } m_{s\,t} \neq +\infty\,,\\
-1 &\quad \text{if } m_{s\,t}= +\infty\,.
\end{array}\right.
\]
For $s \in S$ we define the reflection $r_s: V \to V$ by
\[
r_s (\x) = \x -2 \langle \x, e_s \rangle e_s\,, \quad \x \in V\,.
\]

\begin{theorem}[Tits \cite{Tits1}]
The mapping $s \mapsto r_s$, $s \in S$, determines a 
faithful linear representation $\rho: W \to GL(V)$.
\end{theorem}

The above linear representation is called the {\it canonical representation}
\index{Canonical representation} of $(W,S)$. Note that 
the bilinear form $\langle, \rangle$ is invariant under the action of $W$.

The {\it root system}\index{Root system} $\Phi$ of $(W,S)$ is defined to be the orbit of $\Pi$ under 
the action of $W$, that is,
\[
\Phi= \{ w \cdot e_s\, ;\, w \in W,\ s \in S\}\,.
\]
Let $f \in \Phi$. Write $f= \sum_{s \in S} \lambda_s e_s$, where $\lambda_s \in \R$ for all $s \in S$. 
We say that $f$ is a {\it positive root} \index{Positive root} (resp. a {\it negative root}) 
\index{Negative root}if $\lambda_s \ge 0$ (resp. 
$\lambda_s \le 0$) for all $s \in S$. The set of positive roots (resp. negative roots) is denoted by 
$\Phi_+$ (resp.  by $\Phi_-$). The following is proved in \cite{Bourb1} for the finite root systems, 
but the same proof works in general (see also \cite{Humph1}, \cite{Deodh1}).

\begin{proposition}
We have the disjoint union $\Phi = \Phi_+ \sqcup \Phi_-$.
\end{proposition}

Let $A$ be a finite set that we call an {\it alphabet}. Let $A^\ast$ denote the set of finite sequences 
of elements of $A$ that we call {\it words on $A$}. We define an operation on $A^\ast$ by
\[
(a_1, \dots, a_p) \cdot (b_1, \dots, b_q) = (a_1, \dots, a_p, b_1, \dots, b_q)\,.
\]
Clearly, $A^\ast$ endowed with this operation is a monoid which is called the {\it free monoid}\index{Free monoid} on $A$. 
The unit in $A^\ast$ is the empty word $\epsilon= ()$.

\bigskip\noindent
Each element $w$ in the Coxeter group $W$ can be written in the form $w=s_1s_2 \cdots s_l$, where $s_1, 
s_2, \dots, s_l \in S$. If $l$ is as small as possible, then $l$ is called the {\it word length}\index{Word length} of $w$ 
and is denoted by $l=\lg_S(w)$. If $w= s_1s_2 \cdots s_l$, then the word $\omega= (s_1, s_2, \dots, 
s_l)$ is called an {\it expression} of $w$. If in addition $l=\lg_S(w)$, then $\omega$ is called a 
{\it reduced expression}\index{Reduced expression} of $w$.

For $w \in W$ we set
\[
\Phi_w = \{ f \in \Phi_+\, ;\, w^{-1} f \in \Phi_- \}\,.
\]
Then the word length and the root systems are related by the following.

\begin{proposition}[Bourbaki \cite{Bourb1}]
We have $| \Phi_w| = \lg_S(w)$ for all $w \in W$.
\end{proposition}

Let $G$ be a group. A subset $S \subset G$ is called a {\it positive generating set} of $G$ if it 
generates $G$ as a monoid. Let $S$ be a positive generating set of $G$. 
for $\omega \in S^\ast$, we denote by $\bar \omega$ the element of $G$ represented by $\omega$.
A {\it solution to the word 
problem}\index{Word problem} for $G$ is an algorithm which, given $\omega \in S^\ast$, decides whether $\bar \omega$ is 
trivial or not.

We turn now to describe Tits' solution to the word problem for Coxeter groups.

Let $\omega, \omega' \in S^\ast$. We say that {\it $\omega$ is transformable to $\omega'$ by an 
$M$-operation of type I} if there exist $\omega_1, \omega_2 \in S^\ast$ and $s \in S$ such that
\[
\omega = \omega_1 \cdot (s,s) \cdot \omega_2 \quad \text{and} \quad \omega' = \omega_1 \cdot 
\omega_2\,.
\]
We say that {\it $\omega$ is transformable to $\omega'$by an $M$-operation of type II} if there exist 
$\omega_1, \omega_2 \in S^\ast$ and $s,t \in S$ such that $s \neq t$, $ùm_{s\,t} \neq +\infty$,
\[
\omega= \omega_1 \cdot \pprod(s,t: m_{s\,t}) \cdot \omega_2 \quad \text{and} \quad \omega'= \omega_1 
\cdot \pprod(t,s: m_{s\,t}) \cdot \omega_2\,.
\]
Note that an $M$-operation of type I shortens the length of the word, but not an $M$-operation of type 
II. An $M$-operation of type II is reversible, but not an $M$-operation of type I. If $\omega$ is 
transformable to $\omega'$ by an $M$-operation, then $\bar \omega = \bar \omega'$.

A word $\omega$ is called {\it $M$-reduced} if its length cannot be reduced by means of $M$-operations.

\begin{theorem}[Tits \cite{Tits2}]
\begin{enumerate}
\item
A word $\omega \in S^\ast$ is reduced if and only if it is $M$-reduced.
\item
Let $\omega, \omega' \in S^\ast$ be two reduced words. We have $\bar \omega = \bar \omega'$ if and only 
if one can pass from $\omega$ to $\omega'$ with a finite sequence of $M$-operations of type II.
\end{enumerate}
\end{theorem}

Now, we introduce a partial order on the Coxeter group $W$ whose role is of importance in the study 
of the associated Artin group and monoid.

For $u,v \in W$, we set $u\le_L v$ if there exists $w \in W$ such that $v=uw$ and $\lg_S (v) = \lg_S 
(u) + \lg_S (w)$.

\begin{proposition}[Bourbaki \cite{Bourb1}]
\begin{enumerate}
\item
Let $u,v \in W$. There exists a unique $w^o \in W$ such that $w^o \le_L u$, $w^o \le_L v$, and $w \le_L 
w^o$ whenever $w\le_L u$ and $w \le_L v$.
\item
Suppose that $W$ is finite. Let $u,v \in W$. There exists a unique $w_o \in W$ such that $u\le_L w_o$, 
$v\le_L w_o$, and $w_o \le_L w$ whenever $u\le_L w$ and $v\le_L w$.
\end{enumerate}
\end{proposition}

The element $w^o$ of Proposition 3.7 is denoted by $w^o = u\wedge_L v$, and the element $w_o$ is denoted by 
$w_o=u\vee_L v$ (if it exists). Note that, by the above, $(W, \le_L)$ is a lattice if $W$ is finite. In 
that case, $W$ has a greatest element which is often denoted by $w_0$. 

We finish the subsection with the classification of the spherical type Coxeter graphs.

Recall that, if $\Gamma_1, \dots, \Gamma_l$ are the connected components of a Coxeter graph $\Gamma$, 
then
\[
W_\Gamma = W_{\Gamma_1} \times W_{\Gamma_2} \times \cdots \times W_{\Gamma_l}\,.
\]
In particular, $\Gamma$ is of spherical type if and only if all the components $\Gamma_1, \dots, 
\Gamma_l$ are of spherical type. So, we only need to classify the connected Coxeter graphs of spherical 
type.

\begin{theorem}[Coxeter \cite{Coxet1}, \cite{Coxet2}]
\begin{enumerate}
\item
A Coxeter graph $\Gamma$ is of spherical type if and only if the canonical bilinear form $\langle, 
\rangle : V \times V \to \R$ is positive definite.
\item
The connected spherical type Coxeter graphs are the Coxeter graphs listed in Figure 3.3.
\end{enumerate}
\end{theorem}

%%%%%%%%%%%%%%%
\begin{figure}[htb]
\centerline{
\setlength{\unitlength}{.4cm}
\begin{picture}(25,25)
\put(2,0){\includegraphics[width=9.2cm]{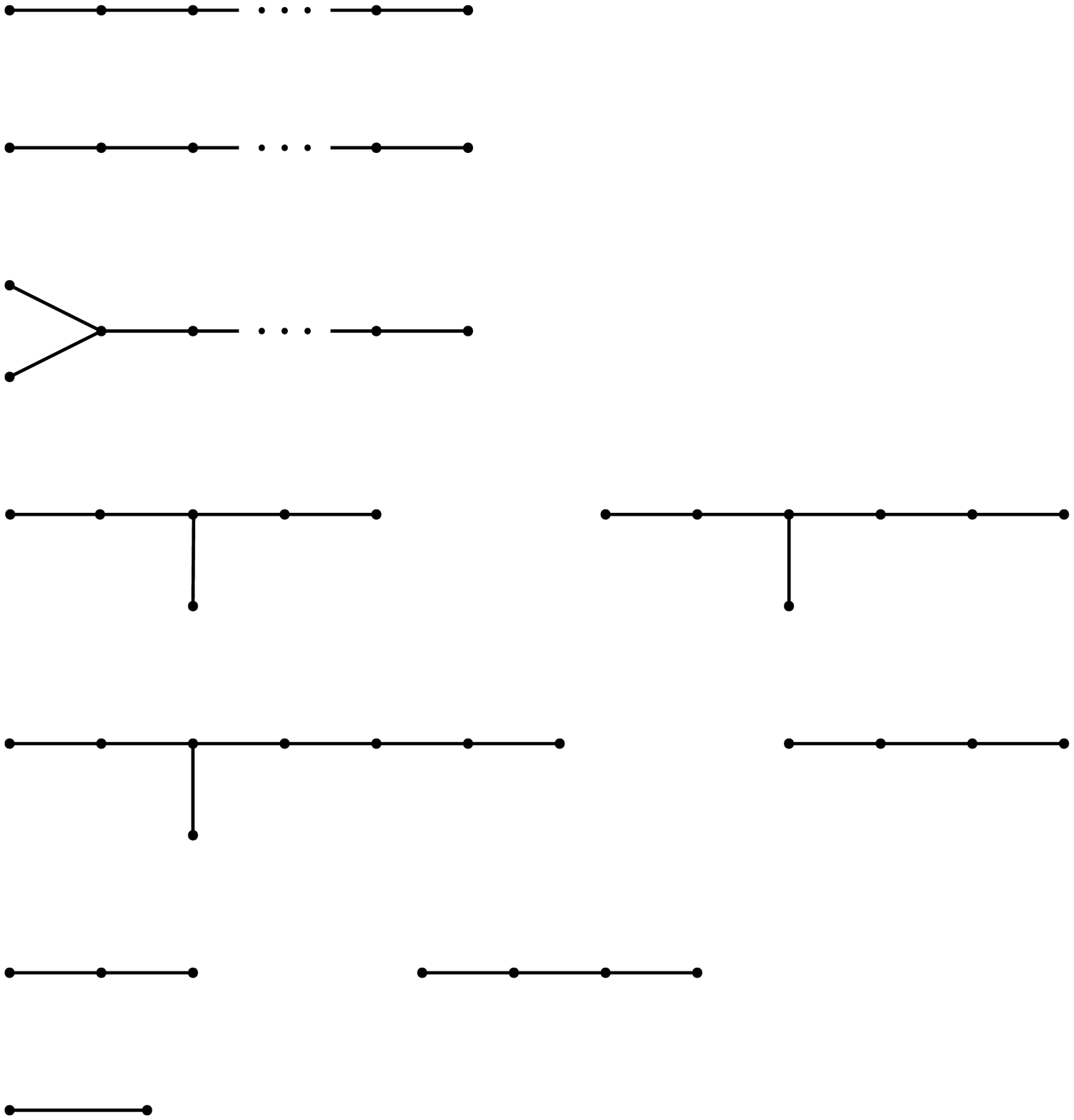}}
\put(0,23.6){\small $A_n$}
\put(15,23.6){\small $n\ge 1$}
\put(0,20.6){\small $B_n$}
\put(15,20.6){\small $n\ge 2$}
\put(0,16.7){\small $D_n$}
\put(15,16.7){\small $n\ge 4$}
\put(0,12.7){\small $E_6$}
\put(13,12.7){\small $E_7$}
\put(0,7.8){\small $E_8$}
\put(17,7.8){\small $F_4$}
\put(0,3){\small $H_3$}
\put(9,3){\small $H_4$}
\put(0,0){\small $I_2(p)$}
\put(8,0){\small $p \ge 5$}
\put(2.8,21.2){\small $4$}
\put(21.5,8.4){\small $4$}
\put(2.8,3.4){\small $5$}
\put(11.7,3.4){\small $5$}
\put(3.5,0.5){\small $p$}
\end{picture}}
\medskip
\centerline{{\bf Figure 3.3.} The connected spherical type Coxeter graphs.}
\end{figure}
%%%%%%%%%%

\subsection{Artin monoids}

Let $\Gamma$ be a Coxeter graph, let $(W,S)$ be the Coxeter system of type $\Gamma$, and let $(G, 
\Sigma)$ be the Artin system of type $\Gamma$. Define the {\it Artin monoid of type $\Gamma$}\index{Artin monoid} to be the 
monoid $G^+=G_\Gamma^+$ presented as a monoid by the generating set $\Sigma = \{ \sigma_s; s \in S\}$ 
and the relations
\begin{multline*}
\pprod (\sigma_s, \sigma_t: m_{s\,t}) =  \pprod (\sigma_t, \sigma_s: m_{s\,t}) \quad \text{for all } 
s,t \in S,\\
\ s \neq t\, \text{ and } m_{s\,t} \neq +\infty\,.
\end{multline*}

\begin{theorem}[Paris \cite{Paris1}]
The natural homomorphism $G_\Gamma^+ \to G_\Gamma$ is 
injective.
\end{theorem}

Recall the homomorphism $\theta: G_\Gamma \to W_\Gamma$, $\sigma_s \mapsto s$. We denote by $\theta^+: 
G_\Gamma^+ \to W_\Gamma$ the restriction of $\theta$ to $G_\Gamma^+$. 
we define a set-section $\kappa : W_\Gamma \to G_\Gamma^+$ of $\theta^+$ as follows.
Let $w \in W$, and let $\omega= (s_1, s_2, \dots, s_l)$ be a reduced expression of $w$. 
Then
\[
\kappa (w) = \sigma_{s_1} \sigma_{s_2} \cdots \sigma_{s_l}\,.
\]
By Theorem 3.6, the definition of $\kappa (w)$ does not depend on the choice of the reduced expression of $w$.

Observe also that the defining relations of $G_\Gamma^+$ are homogeneous, thus $G_\Gamma^+$ has a 
well-defined word length $\lg: G_\Gamma^+ \to \N$, $\sigma_{s_1} \cdots \sigma_{s_l} \mapsto l$. This word 
length satisfies the following properties:
\begin{itemize}
\item
$\lg (\alpha) = 0$ if and only if $\alpha =1$;
\item
$\lg(\alpha \beta) = \lg (\alpha) + \lg (\beta)$ for all $\alpha, \beta \in G_\Gamma^+$.
\end{itemize}
We define partial orders $\le_L$ and $\le_R$ on $G_\Gamma^+$ by
\begin{itemize}
\item
$\alpha \le_L \beta$ if there exists $\gamma \in G_\Gamma^+$ such that $\alpha \gamma = \beta$;
\item
$\alpha \le_R \beta$ if there exists $\gamma \in G_\Gamma^+$ such that $\gamma \alpha = \beta$.
\end{itemize}
The following is again a direct consequence of Theorem 3.6.

\begin{lemma}
Let $u,v \in W$. We have $u \le_L v$ if and only if $\kappa (u) \le_L \kappa 
(v)$.
\end{lemma}

The set $\SS = \{ \kappa (w); w \in W\}$ is called the set of {\it simple elements}\index{Simple element} of $G_\Gamma^+$. If 
$W$ is finite and $w_0$ is the greatest element of $W$, then $\kappa (w_0)$ is called the {\it Garside 
element}\index{Garside element} of $G_\Gamma^+$ and is denoted by $\Delta = \kappa (w_0)$.

The following theorems 3.11 and 3.12 are key results in the study of Artin monoids and groups. They are 
implicit in the work of Brieskorn and Saito \cite{BriSai1}, and explicit for the spherical type Artin groups in the 
work of Deligne \cite{Delig1}. Complete and detailed proofs of them can be found in \cite{Miche1}.

\begin{theorem}
Let $\alpha \in G_\Gamma^+$. Set
\[
E(\alpha) = \{ a \in \SS\,;\, a \le_L \alpha\}\,.
\]
Then $E(\alpha)$ has a greatest element. That is, there exists $a_0 \in E(\alpha)$ such that $E(\alpha) 
= \{ a \in \SS; a \le_L a_0\}$.
\end{theorem}

For $\alpha \in G_\Gamma^+$ we denote by $\delta (\alpha)$ the greatest element of $E(\alpha)$.

\begin{theorem}
Let $\alpha, \beta \in G_\Gamma^+$. Then $\delta( \alpha \beta) = \delta 
(\alpha \delta(\beta))$.
\end{theorem}

Theorems 3.11 and 3.12 have the following consequences whose significance will become clear in the next 
section.

\begin{theorem}
\begin{enumerate}
\item
Let $\alpha, \beta \in G_\Gamma^+$. There exists a unique $\gamma^o \in G_\Gamma^+$ such that 
$\gamma^o \le_L \alpha$, $\gamma^o \le_L \beta$, and $\gamma \le_L \gamma^o$ whenever $\gamma \le_L 
\alpha$ and $\gamma \le_L \beta$.
\item
Suppose that $\Gamma$ is of spherical type. Let $\alpha, \beta \in G_\Gamma^+$. There exists a 
unique $\gamma_o \in G_\Gamma^+$ such that $\alpha \le_L \gamma_o$, $\beta \le_L \gamma_o$, and 
$\gamma_o \le_L \gamma$ whenever $\alpha \le_L \gamma$ and $\beta \le_L \gamma$.
\end{enumerate}
\end{theorem}

The element $\gamma^o$ of Theorem 3.13 is denoted by $\gamma^o = \alpha \wedge_L \beta$, and the element $\gamma_o$ 
is denoted by $\gamma_o=\alpha \vee_L \beta$ (if it exists). Note that the same result is valid if we replace 
$\le_L$ by $\le_R$.

\begin{proof}
We prove (1) by induction on $\lg (\alpha) + \lg (\beta)$. By Proposition 3.7 and by Lemma~3.10, 
$\alpha \wedge_L \beta$ exists if $\alpha, \beta \in \SS$.

Let $\alpha, \beta \in G_\Gamma^+$. Set $a = \delta(a) \wedge_L \delta (\beta)$ ($a$ exists by the 
above observation). If $a=1$, then we must have $\gamma^o = \alpha \wedge_L \beta =1$. Suppose $a 
\neq 1$. Let $\alpha', \beta' \in G_\Gamma^+$ such that $\alpha = a \alpha'$ and $\beta = a \beta'$. 
The element $\alpha' \wedge_L \beta'$ exists by induction. Then $\gamma^o = a \cdot (\alpha' \wedge_L 
\beta')$ (the proof of this equality is left to the reader).

Now, we assume that $\Gamma$ is of spherical type and turn to prove (2). Let $w_0$ be the greatest 
element of $W$, and let $\Delta = \kappa (w_0)$ be the Garside element of $G_\Gamma^+$. It is shown in 
\cite{Bourb1} that $w_0^{-1} = w_0$ and $w_0 S w_0 =S$. This implies that $\Delta \cdot \Sigma 
\cdot\Delta^{-1} = \Sigma$, and, consequently, there exists a permutation $\tau: \SS \to 
\SS$ such that $\Delta \alpha = \tau(\alpha) \Delta$ for all $\alpha \in G_\Gamma^+$.

Let $\alpha \in G_\Gamma^+$. Set $\alpha = a_1 a_2 \cdots a_r$, where $a_i = \delta(a_i a_{i+1} \cdots 
a_r) \in \SS$ for all $1 \le i\le r$. Using the above observation, it is easily shown that $\alpha 
\le_L \Delta^r$.

Let $\alpha, \beta \in G_\Gamma^+$. Set $\EE = \{ \gamma \in G_\Gamma^+; \alpha \le_L \gamma \text{ and 
} \beta \le_L \gamma\}$. We have $\EE \neq \emptyset$ since, by the above, it contains an element of 
the form $\Delta^r$. Let $\gamma_o$ be the smallest element of $\EE$ (this element exists by (1)). 
Then $\gamma_o = \alpha \vee_L \beta$.
\end{proof}

\subsection{Artin groups}

We turn now to present a geometrical interpretation of the Artin groups which extends the 
interpretation of the braid groups in term of configuration spaces. We focus our 
presentation on the spherical type Artin groups, but many of the results stated in this subsection can be 
extended in some sense to the other Artin groups.

Let $\Gamma$ be a spherical type Coxeter graph, let $(W,S)$ be the Coxeter system of type $\Gamma$, and 
let $(G, \Sigma)$ be the Artin system of type $\Gamma$. Recall the set $\Pi= \{ e_s; s \in S\}$ of 
simple roots, the vector space $V= \oplus_{s \in S} \R e_s$, and the canonical bilinear form 
$\langle, \rangle: V \times V \to \R$, which, by Theorem 3.8, is positive definite. We assume that $W$ 
is embedded in $GL(V)$ via the canonical representation.

Let $\RR$ be the set of reflections in $W$. For each $r \in \RR$, let $H_r$ be the hyperplane of $V$ 
fixed by $r$. Then $W$ acts freely on the complement of $\cup_{r \in \RR} H_r$ (see \cite{Bourb1}). 
Complexifying the action, we get an action of $W$ on $V_\C = \C \otimes V$ which is free on the 
complement of $\cup_{r \in \RR} \C \otimes H_r$. Set
\[
M_\Gamma = V_\C \setminus \left( \bigcup_{r \in \RR} \C \otimes H_r \right)\,, \quad N_\Gamma = 
M_\Gamma /W\,.
\]
By a theorem of Chevalley \cite{Cheva1}, Shephard, and Todd \cite{SheTod1},
$V_\C/W$ is isomorphic to $\C^n$, thus $N_\Gamma$ is the 
complement in $\C^n$ of an algebraic set, $( \cup_{r \in \RR} \C \otimes H_r)/W$, called the {\it 
discriminant of type $\Gamma$}\index{Discriminant}.

\begin{theorem}[Brieskorn \cite{Bries1}]
$\pi_1 (N_\Gamma) \simeq G_\Gamma$.
\end{theorem}

\begin{note}
Infinite Coxeter groups also act as reflection groups on $\R^n$. However, to extend Theorem 
3.14 to these groups we should replace $V$ by the Tits cone $U \subset V$ (see \cite{Bourb1}), and 
$V_\C$ by $(U+iV) \subset V_\C$. Then $W$ acts freely on $(U+iV) \setminus (\cup_{r \in \RR} \C \otimes 
H_r) = M_\Gamma$, and it was shown by Van der Lek \cite{VdLek1} that $\pi_1 (N_\Gamma) \simeq G_\Gamma$, where 
$N_\Gamma = M_\Gamma/W$.
\end{note}

An extension of Corollary 2.12 to the spherical type Artin groups is:

\begin{theorem}[Deligne \cite{Delig1}]
Let $\Gamma$ be a spherical type Coxeter graph. Then 
$N_\Gamma$ and $M_\Gamma$ are $K(\pi,1)$.
\end{theorem}

\begin{note}
It is an open problem to know whether $N_\Gamma$ is $K(\pi,1)$ if $\Gamma$ is not of 
spherical type. The answer is yes for the so-called FC-type Artin groups and 2-dimensional Artin groups 
\cite{ChaDav1}, and also for few affine type Artin groups (see \cite{ChaPei1}, \cite{CaMoSa1}).
\end{note}

\begin{note}
We may replace $W$ by a finite complex reflection group acting on $\C^n$, and $M_\Gamma$ 
by $M(W)= \C^n \setminus (\cup_{r \in \RR} H_r)$, where $\RR$ is the set of reflections in $W$, and 
$H_r$ denotes the hyperplane fixed by $r$. Here again, the group $W$ acts freely on $M(W)$ and, by 
\cite{Cheva1} and \cite{SheTod1}, $N(W)=M(W)/W$ is isomorphic to the complement in $\C^n$ of an algebraic set. It was 
recently proved by Bessis \cite{Bessi1} that $N(W)$ is always $K(\pi,1)$. A classification of the 
finite complex reflection groups was obtained by Shephard and Todd \cite{SheTod1}, and a nice 
presentation of $\pi_1( N(W))$ is known for all these groups but four exceptional cases (see 
\cite{BrMaRo1}, \cite{BesMic1}).
\end{note}

%%%%%%%%%%%%%%%%%%%%%%%%%%%%%%%%%

\section{Garside groups}

\subsection{Garside monoids}

A monoid $M$ is called {\it atomic}\index{Atomic monoid} if there exists a function $\nu: M \to \N$ such that
\begin{itemize}
\item
$\nu (\alpha) = 0$ if and only if $\alpha = 1$;
\item
$\nu (\alpha \beta) \ge \nu (\alpha) + \nu (\beta)$ for all $\alpha, \beta \in M$.
\end{itemize}
Such a function $\nu$ is called a {\it norm}\index{Norm} on $M$. An element $\alpha \in M$ is called an 
{\it atom}\index{Atom} 
if it is indecomposable, that is, if $\alpha = \beta \gamma$, then either $\beta=1$ or $\gamma=1$.

The following is proved in \cite{DehPar1}.

\begin{lemma}
Let $M$ be an atomic monoid. A subset $S \subset M$ generates $M$ if and only if 
it contains all the atoms. In particular, $M$ is finitely generated if and only if it contains finitely 
many atoms.
\end{lemma}

Let $M$ be an atomic monoid. We define on $M$ two partial orders $\le_L$ and $\le_R$ as follows.
\begin{itemize}
\item
Set $\alpha \le_L \beta$ if there exists $\gamma \in M$ such that $\alpha \gamma= \beta$.
\item
Set $\alpha \le_R \beta$ if there exists $\gamma \in M$ such that $\gamma \alpha = \beta$.
\end{itemize}
The orders $\le_L$ and $\le_R$ are called the {\it left and right divisibily orders}, respectively.

A monoid $M$ is called a {\it Garside monoid}\index{Garside monoid} if
\begin{itemize}
\item
$M$ is atomic and finitely generated;
\item
$M$ is cancelative (that is, if $\alpha \beta \gamma = \alpha \beta' \gamma$, then $\beta = \beta'$, for 
all $\alpha, \beta, \beta', \gamma \in M$);
\item
$(M, \le_L)$ and $(M, \le_R)$ are lattices;
\item
there exists an element $\Delta \in M$, called a {\it Garside element}\index{Garside element}, 
such that the sets $L(\Delta) = \{ 
\alpha \in M; \alpha \le_L \Delta\}$ and $R(\Delta) = \{ \alpha \in M; \alpha \le_R \Delta \}$ are equal 
and generate $M$.
\end{itemize}
If $M$ is a Garside monoid, then the lattice operations of $(M, \le_L)$ (resp. of $(M,\le_R)$) are denoted 
by $\vee_L$ and $\wedge_L$ (resp. by $\vee_R$ and $\wedge_R$).

Let $M$ be a monoid. The {\it group of fractions}\index{Group of fractions}
 of $M$ is defined to be the group $G(M)$ presented 
with the generating set $M$ and the relations $\alpha \cdot \beta = \gamma$ if $\alpha \beta = \gamma$ 
in $M$. Such a group has the universal property that, if $\varphi: M \to H$ is a homomorphism and $H$ 
is a group, then there exists a unique homomorphism $\hat \varphi : G(M) \to H$ such that $\varphi = 
\hat \varphi \circ \iota$, where $\iota : M \to G(M)$ is the natural homomorphism. Note that the latter 
homomorphism $\iota: M \to G(M)$ is not injective in general.

A {\it Garside group}\index{Garside group} is defined to be the group of fractions of a Garside monoid.

\begin{note}
Garside monoids and groups were introduced in \cite{DehPar1} in a slightly restricted 
sense, and in \cite{Dehor1} in the larger sense which is now generally used. This notion was extended 
to the notion of quasi-Garside monoids \cite{Digne1}, \cite{Bessi2}, to study some non-spherical Artin 
groups. Quasi-Garside monoids have the same definition as the Garside monoids except they are not 
required to be finitely generated. Recently, this notion was extended to the notion of Garside 
categories \cite{Kramm1}, \cite{Kramm2}, \cite{DigMic1}, \cite{Bessi3}, which, in some sense, has to be 
considered as a geometric object more than as an algebraic one. Garside categories are a central concept 
in Bessis' solution to the $K(\pi,1)$ problem for complex reflection arrangements (see \cite{Bessi1}).
\end{note}

Motivating examples of Garside groups are the Artin groups of spherical type:

\begin{theorem}
Let $\Gamma$ be a spherical type Coxeter graph. Then $G_\Gamma^+$ is a Garside 
monoid. In particular, $G_\Gamma$ is a Garside group.
\end{theorem}

Note that Theorem 4.2 is essentially a restatement of Theorem 3.13. 

Other interesting examples of 
Garside groups include all torus link groups (see \cite{Pican1}) and some generalized braid groups 
associated to complex reflection groups (see \cite{Bessi1}).

\begin{note}
Two different Garside monoids can have the same group of fractions. In particular, the 
Artin groups of spherical type are groups of fractions of other Garside monoids, called {\it dual Artin 
monoids}, introduced by Birman, Ko, and Lee \cite{BiKoLe1} for the braid groups, and by Bessis 
\cite{Bessi4} for the other ones.
\end{note}

\begin{note}
A Garside element is not unique. For instance, if $\Delta$ is a Garside element, then 
$\Delta^k$ is a Garside element for all $k \ge 1$ (see \cite{Dehor1}).
\end{note}

We say that a monoid $M$ satisfies the {\it \"Ore conditions} if
\begin{itemize}
\item
$M$ is cancelative;
\item
for all $\alpha, \beta \in M$, there exist $\alpha',\beta' \in M$ such that $\alpha \alpha' = \beta 
\beta'$.
\end{itemize}
It is well-known that a monoid which satisfies the \"Ore conditions\index{\"Ore conditions} embeds in its group of fractions. On 
the other hand, A Garside monoid clearly satisfies the \"Ore conditions. Thus:

\begin{proposition}
Let $M$ be a Garside monoid. Then the natural homomorphism $\iota: M \to 
G(M)$ is injective.
\end{proposition}

Let $M$ be a Garside monoid and let $G = G(M)$ be the group of fractions of $M$. Then the partial 
orders $\le_L$ and $\le_R$ can be extended to $G$ as follows.
\begin{itemize}
\item
Set $\alpha \le_L \beta$ if $\alpha^{-1} \beta \in M$.
\item
Set $\alpha \le_R \beta$ if $\beta \alpha^{-1} \in M$.
\end{itemize}
One can easily verify that $(G, \le_L)$ and $(G, \le_R)$ are lattices. This can be used, for example, 
to prove the following.

\begin{proposition}
A Garside group is torsion free.
\end{proposition}

\begin{proof}
Let $\alpha \in G$ such that $\alpha^n=1$ for some $n \ge 1$. Set $\beta = 1 \vee_L \alpha 
\vee_L \cdots \vee_L \alpha^{n-1}$. It is easily seen that $\le_L$ is invariant by left multiplication. 
This implies that $\alpha \beta = \beta$, hence $\alpha =1$.
\end{proof}

\begin{note}
Let $G$ be a Garside group. Finite dimensional $K(G,1)$ (that is, $K(\pi,1)$ spaces having 
$G$ as fundamental group) were described in \cite{DehLaf1} and \cite{ChMeWh1}. This implies that $G$ is 
torsion free, but also more.
\end{note}

\subsection{Reversing processes and presentations}

Let $\Sigma$ be a finite set. Let $\Sigma^\ast$ be the free monoid on $\Sigma$. Recall that the 
elements of $\Sigma^\ast$ are the finite sequences of elements of $\Sigma$ that are called {\it words} 
on $\Sigma$. If $\equiv$ is a congruence on $\Sigma^\ast$ and $M= (\Sigma^\ast/\equiv)$, then we denote 
by $\Sigma^\ast \to M$, $\omega \mapsto \bar \omega$ the natural epimorphism.

Define a {\it complement} \index{Complement} on $\Sigma$ to be a map $f: \Sigma \times \Sigma \to \Sigma^\ast$ such that 
$f(x,x)= \epsilon$ for all $x \in \Sigma$, where $\epsilon$ denotes the empty word. To a complement 
$f$ we associate two monoids:
\[
\begin{array}{rcl}
M_L^f &=& \langle \Sigma \ |\ x f(x,y) = y f(y,x) \text{ for all } x,y \in \Sigma \rangle^+\,;\\
M_R^f &=& \langle \Sigma \ |\ f(y,x)x = f(x,y)y \text{ for all } x,y \in \Sigma \rangle^+\,.
\end{array}
\]
For $u,v \in \Sigma^\ast$, we use the notation $u \equiv_L^f v$ (resp. $u \equiv_R^f v$) to mean that 
$\bar u = \bar v$ in $M_L^f$ (resp. in $M_R^f$).

\bigskip\noindent
{\bf Example.} Let $\Gamma$ be a Coxeter graph, and let $M= (m_{s\,t})_{s,t\in S}$ be the Coxeter 
matrix of $\Gamma$. Suppose that $m_{s\,t} \neq + \infty$ for all $s,t \in S$, $s \neq t$. Let $\Sigma 
= \{ \sigma_s; s \in S\}$. Let $f: \Sigma \times \Sigma \to \Sigma^\ast$ be the complement defined by
\[
f( \sigma_s, \sigma_t) = \pprod (\sigma_t, \sigma_s: m_{s\,t}-1)\,.
\]
Then $G_\Gamma^+ = M_L^f$.

\bigskip\noindent
Suppose given a complement $f :\Sigma \times \Sigma \to \Sigma^\ast$. Let $\Sigma^{-1} = \{ x^{-1}; x 
\in \Sigma\}$ be the set of {\it inverses} of elements of $\Sigma$. Let $\omega, \omega' \in (\Sigma 
\sqcup \Sigma^{-1} )^\ast$. We say that $\omega$ is {\it $f$-reversible on the left in one step} 
\index{Reversible word} to 
$\omega'$ if there exist $\omega_1, \omega_2 \in (\Sigma \sqcup \Sigma^{-1})^\ast$ and $x,y \in \Sigma$ 
such that
\[
\omega = \omega_1 x^{-1} y \omega_2 \quad \text{and} \quad \omega' = \omega_1 \cdot f(x,y) \cdot 
f(y,x)^{-1} \cdot \omega_2\,.
\]
Note that $y$ can be equal to $x$ in the above definition. In that case we have $\omega = \omega_1 
x^{-1} x \omega_2$ and $\omega' = \omega_1 \omega_2$. Note also that $\overline{\omega} = 
\overline{\omega'}$ in $G(M_L^f)$ if $\omega$ is $f$-reversible on the left in one step to $\omega'$.

Let $p \ge 0$. We say that $\omega$ is {\it $f$-reversible on the left in $p$ steps} to $\omega'$ if 
there exists a sequence $\omega = \omega_0, \omega_1, \dots, \omega_p=\omega'$ in $(\Sigma \sqcup 
\Sigma^{-1})^\ast$ such that $\omega_{i-1}$ is $f$-reversible on the left in one step to $\omega_i$ for 
all $1 \le i\le p$. The property that $\omega$ is $f$-reversible on the left to $\omega'$ is denoted by 
$\omega \mapsto_L^f \omega'$.

We define the {\it $f$-reversibility on the right} in the same way, replacing subwords of the form $y 
x^{-1}$ by their corresponding words $f(x,y)^{-1} \cdot f(y,x)$. The property that $\omega$ is 
$f$-reversible on the right to $\omega'$ is denoted by $\omega \mapsto_R^f \omega'$.

A word $\omega \in (\Sigma \sqcup \Sigma^{-1} )^\ast$ is said to be {\it $f$-reduced on the left}\index{Reduced word} 
(resp. {\it $f$-reduced on the right}) if it is of the form $\omega = v u^{-1}$ (resp. $\omega= u^{-1} 
v$) with $u,v \in \Sigma^\ast$.

It is shown in \cite{Dehor2} that a reversing process is confluent, namely:

\begin{proposition}[Dehornoy \cite{Dehor2}]
Let $f: \Sigma \times \Sigma \to \Sigma^\ast$ be 
a complement, and let $\omega \in (\Sigma \sqcup \Sigma^{-1})^\ast$. Suppose that there exist $p \ge 0$ 
and a $f$-reduced word $v u^{-1}$ on the left such that $\omega$ is $f$-reversible on the left in $p$ 
steps to $v u^{-1}$. Then any sequence of left $f$-reversing transformations 
starting from $\omega$ converges to $v u^{-1}$ in 
$p$ steps.
\end{proposition}

Let $u,v \in \Sigma^\ast$. Suppose there exist $u', v' \in \Sigma^\ast$ such that $u^{-1} v \mapsto_L^f 
v' (u')^{-1}$. By Proposition~4.5, the words $u'$ and $v'$ are unique. Moreover, it is easily checked 
that we also have $v^{-1} u \mapsto_L^f u' (v')^{-1}$. In this case we set
\[
u' = C_L^f (v,u) \quad \text{and} \quad v'= C_L^f (u,v)\,.
\]
Similarly, if there exist $u',v' \in \Sigma^\ast$ such that $v u^{-1} \mapsto_R^f (u')^{-1} (v')$, then  
$u v^{-1} \mapsto_R^f (v')^{-1} (u')$, $u'$ and $v'$ are unique, and we set
\[
u' = C_R^f (u,v) \quad \text{and} \quad v'= C_R^f (v,u)\,.
\]

\begin{lemma}[Dehornoy \cite{Dehor2}]
Let $f: \Sigma \times \Sigma \to \Sigma^\ast$ be a 
complement. Let $u,v \in \Sigma^\ast$. Suppose that $C_L^f (u,v)$ and $C_L^f (v,u)$ exist. Then
\[
u \cdot C_L^f (u,v) \equiv_L^f v \cdot C_L^f (v,u)\,.
\]
\end{lemma}

A complement $f: \Sigma \times \Sigma \to \Sigma^\ast$ is said to be {\it coherent on the left} 
\index{Coherent complement} if, for 
all $x,y,z \in \Sigma$, $C_L^f (f(x,y), f(x,z))$ and $C_L^f (f(y,x), f(y,z))$ exist and are 
$\equiv_L^f$-equivalent. Similarly, we say that $f$ is {\it coherent on the right} if, for all $x,y,z 
\in \Sigma$, $C_R^f (f(z,x), f(y,x))$ and $C_R^f (f(z,y), f(x,y))$ exist and are 
$\equiv_R^f$-equivalent.

\begin{theorem}[Dehornoy, Paris \cite{DehPar1}, \cite{Dehor1}]
Let $M$ be a finitely generated 
monoid, and let $\Sigma$ be a finite generating set of $M$. Then $M$ is a Garside monoid if and only if 
it satisfies the following three conditions.
\begin{itemize}
\item
$M$ is atomic.
\item
There exist a complement $f: \Sigma \times \Sigma \to \Sigma^\ast$ coherent on the left and a 
complement $g: \Sigma \times \Sigma \to \Sigma^\ast$ coherent on the right such that $M= M_L^f = 
M_R^g$.
\item
There exists an element $\Delta \in M$ such that the sets $L(\Delta) = \{ \alpha \in M; \alpha \le_L \Delta \}$ 
and $R(\Delta) = \{ \alpha \in M; \alpha \le_R \Delta \}$ are equal and generate $M$.
\end{itemize}
\end{theorem}

We refer to \cite{DehPar1} and \cite{Dehor1} for
more ``algorithmic'' conditions to detect a Garside monoid in terms of complements and 
presentations, and turn to explain some applications of the 
reversing processes.

Let $M$ be a Garside monoid, and let $G=G(M)$ be its group of fractions. Let $f: \Sigma \times \Sigma 
\to \Sigma^\ast$ and $g: \Sigma \times \Sigma \to \Sigma^\ast$ be complements such that $M=M_L^f = 
M_R^g$.

First, the complements $f$ and $g$ lead to algorithms:

\begin{proposition}[Dehornoy, Paris \cite{DehPar1}, \cite{Dehor1}]
\begin{enumerate}
\item
The complement $f$ is coherent on the left, and the complement $g$ is coherent on the right.
\item
Let $\omega \in (\Sigma \sqcup \Sigma^{-1})^\ast$. There exist a (unique) $f$-reduced word $v u^{-1}$ on the 
left, and a (unique) $g$-reduced word $(u')^{-1} (v')$ on the right, such that $\omega \mapsto_L^f v u^{-1}$ and 
$\omega \mapsto_R^g (u')^{-1} (v')$.
\end{enumerate}
\end{proposition}

They can be used to solve the word problem:

\begin{proposition}[Dehornoy, Paris \cite{DehPar1}, \cite{Dehor1}]
Let $\omega \in (\Sigma 
\sqcup \Sigma^{-1})^\ast$. Let $u,v \in \Sigma$ such that $\omega \mapsto_L^f vu^{-1}$ (see 
Proposition 4.8). Then $\bar \omega = 1$ in $G=G(M)$ if and only if $u^{-1} v \mapsto_L^f \epsilon$, 
where $\epsilon$ denotes the empty word.
\end{proposition}

They can  be also used to compute the lattice operations of $(M, \le_L)$ and $(M, \le_R)$.

\begin{proposition}[Dehornoy, Paris \cite{DehPar1}, \cite{Dehor1}]
Let $u,v \in \Sigma^\ast$. 
Set $u' = C_L^f (u,v)$ and $v' = C_L^f (v,u)$. Then $\bar u \vee_L \bar v$ is represented by
\[
uu' \equiv_L^f vv'\,,
\]
and $\bar u \wedge_L \bar v$ is represented by
\[
C_R^g (u, C_R^g (v',u')) \equiv_L^f C_R^g (v, C_R^g (u',v'))\,.
\]
\end{proposition}

\subsection{Normal forms and automatic structures}

Let $M$ be a Garside monoid, let $G=G(M)$ be the group of fractions of $M$, and let $\Delta$ be a fixed 
Garside element of $M$. Define the set of {\it simple elements} \index{Simple element} to be
\[
\SS = \{ a \in M\ ;\ a \le_L \Delta \} = \{ a \in M\ ;\ a \le_R \Delta \}\,.
\]
By definition, $\SS$ is finite and generates $M$.

Let $\alpha \in M$. Then $\alpha$ can be uniquely written in the form
\[
\alpha = a_1 a_2 \cdots a_l\,,
\]
where $a_1, a_2, \dots, a_l \in \SS$, and
\[
a_i = \Delta \wedge_L (a_i a_{i+1} \cdots a_l) \quad \text{for all } 1 \le i \le l\,.
\]
Such an expression of $\alpha$ is called the {\it normal form} \index{Normal form} of $\alpha$.

Let $\alpha \in G$. Then $\alpha$ can be written in the form $\alpha = \beta^{-1} \gamma$, where $\beta, 
\gamma \in M$ (see Proposition~4.8, for instance). Obviously, we can also assume that $\beta \wedge_L 
\gamma = 1$. In that case $\beta$ and $\gamma$ are unique. Let $\beta = b_1 b_2 \cdots b_p$ be the 
normal form of $\beta$ and let $\gamma = c_1 c_2 \cdots c_q$ be the normal form of $\gamma$. Then the 
expression
\[
\alpha = b_p^{-1} \cdots b_2^{-1} b_1^{-1} c_1c_2 \cdots c_q
\]
is called the {\it normal form} \index{Normal form} of $\alpha$.

There is another notion of normal forms for the elements of $G$, called $\Delta$-normal forms, that are 
used, in particular, in several solutions to the conjugacy problem for $G$. They are defined as follows.

It is easily seen that there exists a permutation $\tau : \SS \to \SS$ such that $\Delta a \Delta^{-1} = 
\tau(a)$ for all $a \in \SS$. Moreover, for all $a \in \SS$, there exists $a^\ast \in \SS$ such that 
$a^\ast a = \Delta$ (i.e. $a^{-1} = \Delta^{-1} a^\ast$). These two observations show that every 
$\alpha \in G$ can be written in the form $\alpha = \Delta^p \beta$, where $p \in \Z$ and $\beta \in 
M$. One can choose $p$ to be maximal, and, in that case, $\beta$ is unique. Let $b_1 b_2 \cdots b_r$ be the 
normal form of $\beta$. Then the expression
\[
\alpha = \Delta^p b_1 b_2 \cdots b_r
\]
is called the {\it $\Delta$-normal form} \index{Normal form} of $\alpha$.

A {\it finite state automaton} \index{Finite state automaton} is a quintuple $\AA = (Q, \SS, T, A, q_0)$, where
\begin{itemize}
\item
$Q$ is a finite set, called the set of {\it states}\index{State};
\item
$\SS$ is a finite set, called the {\it alphabet}\index{Alphabet};
\item
$T$ is a map $T: Q \times \SS \to Q$, called the {\it transition function}\index{Transition function};
\item
$A$ is a subset of $Q$, called the set of {\it accepted states}\index{Accepted state};
\item
$q_0$ is an element of $Q$, called the {\it initial state}\index{Initial state}.
\end{itemize}
The {\it iterated transition function} \index{Iterated transition function} is the map 
$T^\ast: Q \times \SS^\ast \to Q$ defined by 
induction on the length of the second component as follows.
\[
\begin{array}{rcl}
T^\ast (q, \epsilon) &=& q\,,\\
T^\ast (q, x_1 x_2 \cdots x_l) &=& T (T^\ast (q, x_1 \cdots x_{l-1}), x_l)\,.
\end{array}
\]
The set
\[
\LL_\AA = \{ \omega \in \SS^\ast\ ;\ T^\ast (q_0, \omega) \in A\}
\]
is called the {\it language recognized} by $\AA$. A {\it regular language} \index{Regular language}
is a language recognized by a 
finite state automaton.

Let $G$ be a group generated by a finite set $\SS$. Define the {\it word length} of an element $\alpha 
\in G$, denoted by $\lg_\SS (\alpha)$, to be the shortest length of a word in $(\SS \sqcup \SS^{-
1})^\ast$ which represents $\alpha$. The {\it distance} between two element $\alpha, \beta \in G$, 
denoted by $d_\SS (\alpha, \beta)$, is the length of $\alpha^{-1} \beta$.

Let $\LL \subset (\SS \sqcup \SS^{-1})^\ast$ be a language. We say that $\LL$ {\it represents} $G$ if 
every element of $G$ is represented by an element of $\LL$. We say, furthermore, that $\LL$ has the 
{\it uniqueness property} \index{Uniqueness property (language)}
if every element of $G$ is represented by a unique element of $\LL$. We say 
that $\LL$ is {\it symmetric} \index{Symmetric language} if $\LL^{-1} = \LL$, 
where $\LL^{-1} = \{ \omega^{-1}; \omega \in \LL\}$. 
We say that $\LL$ is {\it geodesic} \index{Geodesic language} if $\lg (\omega) = \lg_\SS (\bar \omega)$ for all 
$\omega \in \LL$. 
Let $\omega = x_1^{\varepsilon_1} \cdots x_l^{\varepsilon_l} \in (\SS \sqcup \SS^{-1})^\ast$. For $t 
\in \N$ we set
\[
\bar \omega(t) = \left\{
\begin{array}{ll}
1 &\quad \text{if } t=0\\
\overline{x_1^{\varepsilon_1} \cdots x_t^{\varepsilon_t}} & \quad \text{if } 1 \le t\le l\\
\bar \omega & \quad \text{if } t \ge l
\end{array}\right.
\]
Let $c$ be a positive integer. We say that $\LL$ has the {\it $c$-fellow traveler property} 
\index{Fellow traveler property} if
\[
d_\SS (\bar u(t), \bar v(t)) \le c \cdot d_\SS (\bar u, \bar v)
\]
for all $u,v \in \LL$ and all $t \in \N$.

A group $G$ is said to be {\it automatic} \index{Automatic group} if there exist a finite generating set $\SS \subset G$, a 
regular language $\LL \subset (\SS \sqcup \SS^{-1})^\ast$, and a constant $c>0$, such that $\LL$ 
represents $G$ and has the $c$-fellow traveler property. If, in addition, $\LL^{-1}$ has also the 
$c$-fellow traveler property, then $G$ is said to be {\it biautomatic}\index{Biautomatic group}. We say that $G$ is 
{\it fully 
biautomatic} \index{Fully biautomatic group}
if $\LL$ is symmetric, and that $G$ is {\it geodesically automatic} \index{Geodesically automatic group}
if $\LL$ is geodesic.

Biautomatic groups have many attractive properties. For instance, they have soluble word and 
conjugacy problems, and they have quadratic isoperimetric inequalities. We refer to \cite{ECHLPT1} for a 
general exposition on the subject.

\begin{theorem}[Charney \cite{Charn1}, Dehornoy, Paris \cite{DehPar1}]
Let $M$ be a Garside 
monoid, and let $G=G(M)$ be the group of fractions of $M$. Let $\LL \subset (\SS \sqcup \SS^{-1})^\ast$ 
be the language of normal forms. Then $\LL$ is regular, represents $G$, has the uniqueness property, 
has the 5-fellow traveler property, is symmetric, and is geodesic.
\end{theorem}

\begin{corollary}
Garside groups are fully geodesically biautomatic.
\end{corollary}

\begin{note}
The language of $\Delta$-normal forms is also regular and satisfies some fellow traveler 
property, and the language of inverses of $\Delta$-normal forms satisfies the same fellow traveler 
property. So, $\Delta$-normal forms determine another biautomatic structure on $G$. This was proved by 
Thurston \cite{ECHLPT1} for the braid groups and by Charney \cite{Charn2} for all the spherical type 
Artin groups, but the same proof works in general for all Garside groups.
\end{note}

\subsection{Conjugacy problem}

Let $G$ be a group and let $\SS$ be a finite generating set of $G$. A {\it solution to the conjugacy 
problem} \index{Conjugacy problem}
 for $G$ is an algorithm which, for given $u,v \in (\SS \sqcup \SS^{-1})^\ast$, decides whether 
$\bar u$ and $\bar v$ are conjugate or not.

The first solution to the conjugacy problem for the braid groups was obtained by Garside 
\cite{Garsi1}. Garside's algorithm was improved by El-Rifai and Morton \cite{ElrMor1}, and this 
improvement was extended to the Garside groups by Picantin \cite{Pican2}. Picantin's algorithm was 
improved by Franco and Gonz\'alez-Meneses \cite{FraGon1}, then by Gebhardt \cite{Gebha1}, and now by 
Gebhardt and Gonz\'alez-Meneses \cite{GebGon1}. The algorithm that we present here is not the optimal 
one, but is probably the simplest one. It is based on the algorithm of \cite{GebGon1}.

\begin{note}
In addition to the above mentioned papers, there are several recent papers where the 
algorithms are analyzed, in particular to obtain the best possible complexity (see \cite{BiGeGo1}, 
\cite{BiGeGo2}, \cite{BiGeGo3}, \cite{GebGon2}, \cite{Lee1}, \cite{LeeLee2}, \cite{LeeLee3}). These 
analysis often lead to new and unexpected results on the braid groups and, more generally, on the 
Garside groups.
\end{note}

Let $M$ be a Garside monoid, let $G=G(M)$ be its group of fractions, let $\Delta$ be a fixed Garside 
element, and let $\SS = \{ a \in M; a \le_L \Delta \}$ be the set of simple elements. Recall 
that, for every $a \in \SS$, there exists a unique $a^\ast \in \SS$ such that $a a^\ast = \Delta$. 
Recall also that there exists a permutation $\tau: \SS \to \SS$ such that $\Delta a \Delta^{-1} = 
\tau(a)$ for all $a \in \SS$.

Let $\alpha \in G$. Let $\alpha = \Delta^p a_1 a_2 \cdots a_r$ be the $\Delta$-normal form of $\alpha$. 
The number $p$ is called the {\it infimum} of $\alpha$ and is denoted by $\inf (\alpha)$, $p+r$ is 
called the {\it supremum} and is denoted by $\sup(\alpha)$, and $r$ is called the {\it canonical length 
} and is denoted by $\| \alpha \|$. The above terminology comes from the fact that $p$ is the 
greatest number $n$ such that $\Delta^n \le_L \alpha$, and $p+r$ is the smallest number $n$ such that 
$\alpha \le_L \Delta^n$. The (simple) element $\tau^p(a_1)$ is called the {\it initial factor} of 
$\alpha$ and is denoted by $i(\alpha)$, and $a_r$ is called the {\it terminal factor} and is denoted by 
$t(\alpha)$. It is easily checked that $i(\alpha^{-1}) = t(\alpha)^\ast$. Let
\[
\pi (\alpha) = i(\alpha) \wedge_L t(\alpha)^\ast = i(\alpha) \wedge_L i( \alpha^{-1})\,.
\]
Define the {\it sliding} of $\alpha$ to be
\[
S (\alpha) = \pi(\alpha)^{-1} \cdot \alpha \cdot \pi(\alpha)\,.
\]
Observe that $\| S(\alpha)\| \le \| \alpha \|$.

For $\alpha, \beta \in G$, we use the notation $\alpha \sim \beta$ to mean that $\alpha$ is conjugate 
to $\beta$. Let $\alpha \in G$. Define the {\it sliding circuits} \index{Sliding circuits} of $\alpha$ to be
\[
SC (\alpha) = \{ \beta \in G\ ;\ \beta \sim \alpha \text{ and } S^m (\beta) = \beta \text{ for some } m 
\ge 1\}\,.
\]
It is shown in \cite{GebGon1} that the elements of $SC(\alpha)$ have minimal canonical length in the 
conjugacy class of $\alpha$, but not all the elements of the conjugacy class of minimal canonical 
length belong to $SC(\alpha)$.

Clearly, if $\alpha \sim \beta$, then $SC(\alpha) = SC(\beta)$, and if $\alpha \not \sim \beta$, then 
$SC (\alpha) \cap SC (\beta) = \emptyset$. 
So, our solution to the conjugacy problem for $G$ follows the following stages.

\bigskip\noindent
{\bf Input.} Two elements $\alpha, \beta \in G$.

\bigskip\noindent
{\bf Stage 1.} Calculate an element $\alpha_0 \in SC (\alpha)$ and an element $\beta_0 \in SC (\beta)$.

\bigskip\noindent
{\bf Stage 2.} Calculate the whole set $SC (\alpha) = SC(\alpha_0)$ from $\alpha_0$.

\bigskip\noindent
{\bf Output.} YES if $\beta_0 \in SC(\alpha)$, and NO otherwise.

\bigskip\noindent
In order to find an element of $SC(\alpha)$ we use the following which is easy to prove.

\begin{lemma}
Let $\alpha \in G$. There exists $m,k \ge 1$ such that $S^{m+k}(\alpha) = S^k 
(\alpha)$. In particular, $S^k (\alpha) \in SC (\alpha)$.
\end{lemma}

The key result for Stage 2 is the following.

\begin{theorem}[Gebhardt, Gonz\'alez-Meneses \cite{GebGon1}]
Let $\alpha, \beta \in G$ and let 
$\gamma_1, \gamma_2 \in M$. If $\beta$, $\gamma_1^{-1} \beta \gamma_1$, and $\gamma_2^{-1} \beta 
\gamma_2$ are elements of $SC(\alpha)$, then $(\gamma_1 \wedge_L \gamma_2 )^{-1} \beta (\gamma_1 
\wedge_L \gamma_2)$ is also an element of $SC(\alpha)$.
\end{theorem}

\begin{corollary}
Let $\alpha, \beta, \gamma \in G$ such that $\beta$ and $\gamma^{-1} \beta 
\gamma$ are elements of $SC(\alpha)$. Let $\gamma = \Delta^p c_1 c_2 \cdots c_r$ be the $\Delta$-normal 
form of $\gamma$. Set $\beta_0= \Delta^{-p} \beta \Delta^p$, and $\beta_i = c_i^{-1} \beta_{i-1} c_i$ 
for $1 \le i \le r$. Then $\beta_i \in SC(\alpha)$ for all $0 \le i\le r$.
\end{corollary}

\begin{proof}
We prove that $\beta_i \in SC (\alpha)$ by induction on $i$. It is easily seen that, if 
$\beta \in SC (\alpha)$, then $\Delta^{-1} \beta \Delta \in SC (\alpha)$. In particular, we have 
$\beta_0 = \Delta^{-p} \beta \Delta^p \in SC(\alpha)$.

Let $i>0$. By induction, $\beta_{i-1} \in SC (\alpha)$. By the above observation, we have 
$\Delta^{-1} \beta_{i-1} \Delta \in SC (\alpha)$. On the other hand, we have $\gamma^{-1} \beta \gamma = (c_i c_{i+1} \cdots c_r)^{-1} 
\linebreak
\beta_{i-1} (c_i c_{i+1} \cdots c_r) 
\in SC(\alpha)$. By definition of 
a normal form, we have $\Delta \wedge_L (c_i c_{i+1} \cdots c_r) = c_i$. We conclude by Theorem 4.14 
that $\beta_i = c_i^{-1} \beta_{i-1} c_i \in SC (\alpha)$.
\end{proof}

From Corollary 4.15 we obtain the following which, together with 
\linebreak
Lemma~4.13, provides an algorithm to 
compute $SC (\alpha)$.

\begin{corollary}
Let $\alpha \in G$. Let $\Omega_\alpha$ be the graph defined by the 
following data.
\begin{itemize}
\item
The set of vertices of $\Omega_\alpha$ is $SC (\alpha)$.
\item
Two vertices $\beta, \beta' \in SC (\alpha)$ are joined by an edge if there exists $a \in \SS$ such 
that $\beta' = a^{-1} \beta a$.
\end{itemize}
Then $\Omega_\alpha$ is connected.
\end{corollary}

%%%%%%%%%%%%%%%%%%%%%%%%%%%%%%%%%

\section{Cohomology and Salvetti complex}

\subsection{Cohomology}

Let $\Gamma$ be a Coxeter graph, let $(W_\Gamma, S)$ be the Coxeter system of type $\Gamma$, and let 
$(G_\Gamma, \Sigma)$ be the Artin system of type $\Gamma$. Let $\Gamma_{ab}$ be the graph defined by 
the following data.
\begin{itemize}
\item
$S$ is the set of vertices of $\Gamma$;
\item
two vertices $s,t \in S$ are joined by an edge if $m_{s\,t} \neq + \infty$ and $m_{s\,t}$ is odd.
\end{itemize}
The following is easy to prove from the presentation of $G_\Gamma$.

\begin{proposition}
Let $d$ be the number of connected components of $\Gamma_{ab}$. Then the 
abelianization of $G_\Gamma$ is a free abelian group of rank $d$. In particular, $H^1 (G_\Gamma, \Z) 
\simeq \Z^d$.
\end{proposition}

Now, assume that $\Gamma$ is of spherical type, and recall the space $N_\Gamma$ defined in Subsection~3.4. 
Except Proposition 5.1, all know results on the cohomology of $G_\Gamma$ use the fact that $\pi_1 
(N_\Gamma)= G_\Gamma$ (see Theorem 3.14), and $N_\Gamma$ is a $K(\pi,1)$ space (see Theorem 3.15). 
Recall that these two results imply that $H^\ast (G_\Gamma, A)= H^\ast (N_\Gamma,A)$ for any 
$G_\Gamma$-module $A$.

In \cite{Arnol1} Arnol'd established the following properties on the cohomology of the braid groups.

\begin{theorem}[Arnol'd \cite{Arnol1}]
Let $n \ge 2$.
\begin{enumerate}
\item
$H^0 (\BB_n, \Z) = H^1 (\BB_n, \Z) = \Z$, $H^q (\BB_n, \Z)$ is finite for all $q \ge 2$, and $H^q 
(\BB_n, \Z) = 0$ for all $q \ge n$.
\item
If $n$ is even, then $H^q (\BB_n, \Z) = H^q (\BB_{n+1}, \Z)$ for all $q \ge 0$.
\item
$H^q (\BB_n, \Z) = H^q (\BB_{2q-2} , \Z)$ for all $q \le \frac{1}{2} n+1$.
\end{enumerate}
\end{theorem}

The study of the cohomology of the braid groups was continued by Fuks \cite{Fuks1} who calculated the 
cohomology of $\BB_n$ with coefficients in $\F_2 = \Z /2\Z$. Let $\BB_\infty = \varinjlim 
\BB_n$, where the limit is taken relative to the natural embeddings $\BB_n \hookrightarrow 
\BB_{n+1}$, $n \ge 2$.

\begin{theorem}[Fuks \cite{Fuks1}]
\begin{enumerate}
\item
$H^\ast (\BB_\infty, \F_2)$ is the exterior $\F_2$-algebra generated by $\{ a_{m,k}; m\ge 1 \text{ and 
} k \ge 0\}$ where $\deg a_{m,k} = 2^k (2^m-1)$.
\item
The natural embedding $\BB_n \hookrightarrow \BB_\infty$ induces a surjective homomorphism $H^\ast 
(\BB_\infty, \F_2) \to H^\ast (\BB_n, \F_2)$ whose kernel is generated by the monomials
\[
a_{m_1,k_1} a_{m_2,k_2} \cdots a_{m_t,k_t}
\]
such that
\[
2^{m_1+ \cdots m_t+k_1+ \cdots +k_t} >n\,.
\]
\end{enumerate}
\end{theorem}

Later, the cohomology with coefficients in $\F_p = \Z/p\Z$ (where $p$ is an odd prime number) and the 
cohomology with coefficients in $\Z$ were calculate by Cohen \cite{Cohen1}, Segal \cite{Segal1}, and 
Va\v{\i}n\v{s}te\v{\i}n \cite{Vains1}.

\begin{theorem}[Cohen \cite{Cohen1}, Segal \cite{Segal1}, Va\v{\i}n\v{s}te\v{\i}n \cite{Vains1}]
\begin{enumerate}
\item
$H^\ast (\BB_\infty, \F_p)$ is the tensor product of a polynomial algebra generated by $\{ x_i; i\ge 
0\}$, where $\deg x_i = 2 p^{i+1} -2$, and an exterior algebra generated by $\{ y_j; j \ge 0\}$, where 
$\deg y_j = 2 p^j -1$.
\item
The natural embedding $\BB_n \hookrightarrow \BB_\infty$ induces a surjective homomorphism $H^\ast 
(\BB_\infty, \F_p) \to H^\ast (\BB_n, \F_p)$, whose kernel is generated by the monomials
\[
x_{i_1} x_{i_2} \cdots x_{i_s} y_{j_1} y_{j_2} \cdots y_{j_t}
\]
such that
\[
2 (p^{i_1+1} + \cdots + p^{i_s+1} + p^{j_1} + \cdots p^{j_t}) >n\,.
\]
\end{enumerate}
\end{theorem}

Let $\beta_2: H^\ast (\BB_n, \F_2) \to H^\ast (\BB_n, \F_2)$ be the homomorphism defined by
\[
\beta_2 (a_{m,k}) = a_{m+1,0} a_{m,1} \cdots a_{m,k-1}\,.
\]
For an odd prime number $p$, let $\beta_p: H^\ast (\BB_n, \F_p) \to H^\ast (\BB_n, \F_p)$ be the 
homomorphism defined by
\[
\beta_p (x_i) = y_{i+1}\,, \quad \beta_p (y_j) = 0\,.
\]

\begin{theorem}[Cohen \cite{Cohen1}, Va\v{\i}n\v{s}te\v{\i}n \cite{Vains1}]
Let $q \ge 2$. Then
\[
H^q (\BB_n, \Z)= \bigoplus_p \beta_p (H^{q-1} (\BB_n, \F_p))\,,
\]
where the sum is over all primes $p$.
\end{theorem}

The integral cohomology of the Artin groups of type $B$ and $D$ were calculated by Goryunov 
\cite{Goryu1} in terms of the cohomology groups of the braid groups.

\begin{theorem}[Goryunov \cite{Goryu1}]
\begin{enumerate}
\item
Let $n \ge 2$, and let $q \ge 2$. Then
\[ 
H^q (G_{B_n}, \Z) = \bigoplus_{i=0}^n H^{q-i} (\BB_{n-i}, \Z)\,.
\]
\item
Let $n \ge 4$, and let $q \ge 2$. Then
\begin{multline*}
H^q (G_{D_n}, \Z) = H^q (\BB_n, \Z) \oplus \left( \bigoplus_{i=0}^{+\infty} \Ker\, \gamma_{n-2i}^{q-2i} 
\right) \oplus\\ 
\left( \bigoplus_{j=0}^{+\infty} H^{q-2j-3} (\BB_{n-3j-3}, \F_2)\right)\,,
\end{multline*}
where, for $k \ge 2$ and $j \ge 0$, $\gamma_k^j : H^j (\BB_k, \Z) \to H^j (\BB_{k-1}, \Z)$ denotes the 
homomorphism induced by the inclusion $\BB_{k-1} \hookrightarrow \BB_k$.
\end{enumerate}
\end{theorem}

Finally, the integral cohomology of the remainder irreducible Artin groups of spherical type were calculate by 
Salvetti in \cite{Salve1}.

\begin{theorem}[Salvetti \cite{Salve1}]
The integral cohomology of the Artin groups of type 
$I_2(p)$ ($p \ge 5$), $E_6$, $E_7$, $E_8$, $F_4$, $H_3$, and $H_4$ is given in Table 5.1.
\end{theorem}

%%%%%%%%%%%%%%%%%%%%%%%%
\begin{table}[htb]
\begin{center}
{\small\begin{tabular}{|c||c|c|c|c|c|}
\hline
&&&&&\\
\noalign{\vskip-8pt}
&$H^0$&$H^1$&$H^2$&$H^3$&$H^4$\\
\hline\hline
&&&&&\\
\noalign{\vskip-9pt}
$I_2(2q)$&$\Z$&$\Z^2$&$\Z$&$0$&$0$\\
\hline
$I_2(2q+1)$&$\Z$&$\Z$&$0$&$0$&$0$\\
\hline
$H_3$&$\Z$&$\Z$&$\Z$&$\Z$&$0$\\
\hline
$H_4$&$\Z$&$\Z$&$0$&$\Z \times \Z_2$&$\Z$\\
\hline
&&&&&\\
\noalign{\vskip-9pt}
$F_4$&$\Z$&$\Z^2$&$\Z^2$&$\Z^2$&$\Z$\\
\hline
$E_6$&$\Z$&$\Z$&$0$&$\Z_2$&$\Z_2 $\\
\hline
$E_7$&$\Z$&$\Z$&$0$&$\Z_2$&$\Z_2 \times \Z_2$\\
\hline
$E_8$&$\Z$&$\Z$&$0$&$\Z_2$&$\Z_2$\\
\hline
\end{tabular}}
\end{center}
\medskip
\centerline{{\bf Table 5.1.(a).} Cohomology of the spherical type Artin groups.}
\end{table}

\begin{table}[htb]
\begin{center}
{\small\begin{tabular}{|c||c|c|c|c|}
\hline
&&&&\\
\noalign{\vskip-8pt}
&$H^5$&$H^6$&$H^7$&$H^8$\\
\hline\hline
$E_6$&$\Z_6$&$\Z_3$&$0$&$0$\\
\hline
$E_7$&$\Z_6 \times \Z_6$&$\Z_3 \times \Z_6 \times \Z$&$\Z$&$0$\\
\hline
$E_8$&$\Z_2 \times \Z_6$&$\Z_3 \times \Z_6$&$\Z_2 \times \Z_6\times \Z$&$\Z$\\
\hline
\end{tabular}}
\end{center}
\medskip
\centerline{{\bf Table 5.1.(b).} Cohomology of the spherical type Artin groups.}
\end{table}
%%%%%%%%%%%%%%%%%%%%%%%%%%%%
\begin{note}
It is a direct consequence of \cite{Salve2} that $N_\Gamma$ has the same homotopy type as 
a CW-complex of dimension $n$, where $n= |S|$. This implies that the cohomological dimension of 
$G_\Gamma$ is $\le n$, and, therefore, that $H^q (G_\Gamma, \Z) = 0$ for all $q >n$.
\end{note}

\begin{note}
Recall the space $M_\Gamma$ of Subsection 3.4. The cohomology $H^\ast (M_\Gamma, \Z)$ was 
calculate by Brieskorn in \cite{Bries2}. In particular, $H^\ast (M_\Gamma, \Z)$ is torsion free and 
$H^n (M_\Gamma, \Z) \neq 0$. Let $CG_\Gamma$ be the kernel of the canonical epimorphism $\theta: 
G_\Gamma \to W_\Gamma$. By \cite{Delig1} we have $H^n (M_\Gamma, \Z)= H^n (CG_\Gamma,\Z)$, thus, by the 
above, $\cd (G_\Gamma) = \cd (CG_\Gamma) \ge n$, where $\cd (G_\Gamma)$ denotes the cohomological 
dimension of $G_\Gamma$. We already know that $\cd (G_\Gamma) \le n$, thus $\cd 
(G_\Gamma)=n$.
\end{note}

\begin{note}
The ring structure of $H^\ast (G_\Gamma, \Z)$, where $\Gamma$ is a Coxeter graph in the 
list of Theorem 5.7, was calculated in \cite{Landi1}. Some cohomologies with twisted coefficients were 
also considered.
An interesting case is the cohomology over the module of Laurent polynomials $\Q [q^{\pm 1}]$ (resp. 
$\Z [q^{\pm 1}]$), because it determines the rational (resp. integral) cohomology of the Milnor fiber 
of the discriminant of type $\Gamma$ (see \cite{Calle1}).
For the case 
$\Gamma = A_n$ ({\it i.e.} $G_\Gamma$ is the braid group $\BB_{n+1}$), the $\Q [q^{\pm 1}]$-cohomology 
was calculated by several people in several ways (see \cite{Frenk1}, \cite{Marka1}, \cite{CohSuc1}, 
\cite{DePrSa1}), and the $\Z [q^{\pm 1}]$-cohomology was calculated by Callegaro in \cite{Calle2}. The 
$\Q [q^{\pm 1}]$-cohomology for the other spherical type Artin groups was calculated in \cite{DePrSaSt1}. The 
$\Z[q^{\pm 1}]$-cohomology for the exceptional cases was calculated in \cite{CalSal1}, and the top $\Z 
[q^{\pm 1}]$-cohomology for all cases was calculated in \cite{DeSaSt1}.
\end{note}

\begin{note}
The cohomology of the non-spherical Artin groups is badly understood. Some 
calculations for the type $\tilde A_n$ were done in \cite{CaMoSa1}.
\end{note}

We refer to \cite{Versh1} for a more detailed exposition on the cohomology of the braid groups and the 
Artin groups of spherical type, and turn to present the Salvetti complex (of a real hyperplane 
arrangement). This is the main tool in Salvetti's calculations of the cohomology of Artin groups (see 
\cite{Salve1}), but it can be used for other purposes. For instance, it can be also used to prove Theorems 
3.14 and 3.15 (see \cite{Salve3} and \cite{Paris2}), and to produce a free resolution of $\Z$
by $\Z[G_\Gamma]$-modules (see Theorem 5.15).

\subsection{Salvetti complex}

Define a {\it (real) hyperplane arrangement} \index{Hyperplane arrangement} 
to be a finite family $\AA$ of linear hyperplanes of 
$\R^n$. For every $H \in \AA$ we denote by $H_\C$ the hyperplane of $\C^n$ having the same equation as 
$H$ ({\it i.e.} $H_\C = \C \otimes H$), and we set
\[
M(\AA) = \C^n \setminus \left( \bigcup_{H \in \AA} H_\C \right)\,.
\]
Note that $M(\AA)$ is an open connected subvariety of $\C^n$.

The arrangement $\AA$ subdivides $\R^n$ into {\it facets}\index{Facet}. We denote by $\FF (\AA)$ the set of all 
facets. The {\it support} of a facet $F \in \FF (\AA)$ is the linear subspace $\langle F \rangle$ 
spanned by $F$. We denote by $\bar F$ the closure of a facet $F$. We order $\FF (\AA)$ by $F \le G$ if 
$F \subset \bar G$. The set $\FF (\AA)$ has a unique minimal element: $\cap_{H \in \AA}H$. The maximal 
elements of $\FF(\AA)$ are the facets of codimension $0$, and they are called {\it chambers}\index{Chamber}. The set 
of all chambers is denoted by $\CC (\AA)$.

Set
\[
\XX = \{ (F,C) \in \FF(\AA) \times \CC(\AA)\ ;\ F \le C\}\,.
\]
We partially order $\XX$ as follows. For $F \in \FF(\AA)$ we set $\AA_F = \{ H \in \AA; H \supset F\}$. 
For $F \in \FF (\AA)$ and $C \in \CC (\AA)$ we denote by $C_F$ the chamber of $\AA_F$ which contains 
$C$. We set
\[
(F_1,C_1) \le (F_2,C_2)\quad \text{if} \quad F_1 \le F_2 \text{ and } (C_1)_{F_2} = (C_2)_{F_2}\,.
\]
(See Figure 5.1.)

%%%%%%%%%%%
\begin{figure}[htb]
\centerline{
\setlength{\unitlength}{.4cm}
\begin{picture}(25,8)
\put(2,0.5){\includegraphics[width=8cm]{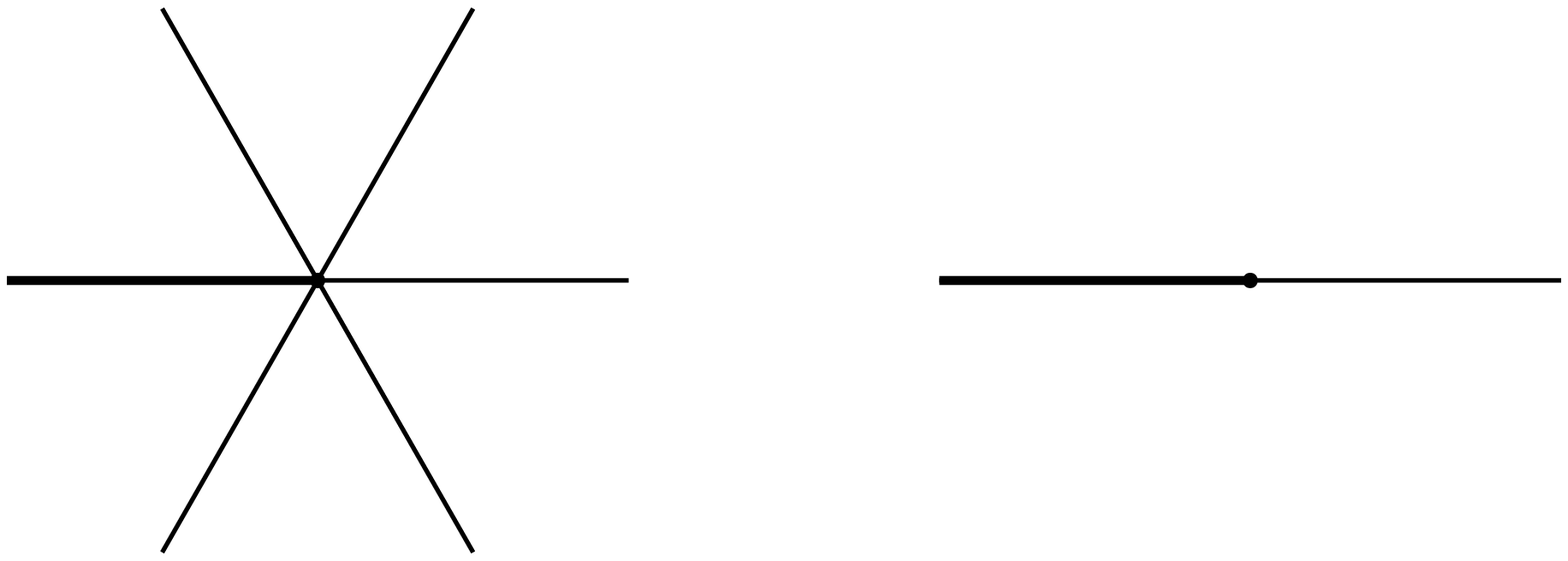}}
\put(0.5,3.9){\small $F_2$}
\put(3.5,2.5){\small $C_2$}
\put(5.5,1.5){\small $C_1$}
\put(6.5,3.4){\small $F_1$}
\put(15,3){\small $(C_1)_{F_2} = (C_2)_{F_2}$}
\put(22.5,3.8){\small $\langle F_2 \rangle$}
\end{picture}}
\medskip
\centerline{{\bf Figure 5.1.} Order in $\XX$.}
\end{figure}
%%%%%%%%%

Define the {\it Salvetti complex} \index{Salvetti complex}
$\Sal (\AA)$ of $\AA$ to be the (geometric realization of the)
flag complex of $(\XX, \le)$. That 
is, to every chain $X_0 <X_1< \cdots <X_d$ in $\XX$ corresponds a simplex $\Delta (X_0, X_1, \dots, 
X_d)$ of $\Sal (\AA)$, and every simplex of $\Sal (\AA)$ is of this form.

\begin{theorem}[Salvetti \cite{Salve2}]
The simplicial complex $\Sal (\AA)$ is homotopy 
equivalent to $M(\AA)$.
\end{theorem}

We turn now to describe a cellular decomposition of $\Sal (\AA)$ which is the version which is usually 
used in the literature.

Without loss of generality, we can and do assume that $\AA$ is 
{\it essential}\index{Essential arrangement}, that is, $\cap_{H \in 
\AA} H = \{ 0\}$. Consider the unit sphere $\S^{n-1}= \{ {\bf x} \in \R^n; \| {\bf x}\| =1\}$. The 
arrangement $\AA$ determines a cellular decomposition of $\S^{n-1}$: to each facet $F \in \FF(\AA) 
\setminus \{0\}$ corresponds the open cell $F \cap \S^{n-1}$, and each cell is of this form. This 
cellular decomposition is regular in the sense that the closure of a cell is a closed disk. Hence, one 
can consider the barycentric subdivision. For each facet $F \in \FF (\AA) \setminus \{ 0\}$ we fix a 
point ${\bf x} (F) \in F \cap \S^{n-1}$. To each chain $\{0\} \neq F_0<F_1< \cdots <F_d$ in $\FF (\AA) \setminus \{0\}$ 
corresponds a simplex $\Delta (F_0, F_1, \dots, F_d)$ whose vertices are ${\bf x} (F_0), {\bf x} (F_1), 
\dots, {\bf x} (F_d)$, and every simplex of $\S^{n-1}$ is of this form. So, the simplicial 
decomposition of $\S^{n-1}$ is the flag complex of $( \FF (\AA) \setminus \{ 0\}, \le)$.

We extend the above simplicial decomposition of $\S^{n-1}$ to a simplicial decomposition of the 
$n$-disk $\B^n= \{ {\bf x} \in \R^n; \| {\bf x}\| \le 1\}$, adding a single vertex ${\bf x}(0)=0$. That is, 
we view $\B^n$ as the cone of $\S^{n-1}$. Now, to any chain $F_0<F_1< \cdots <F_d$ in $\FF (\AA)$ 
corresponds a simplex $\Delta (F_0, F_1, \dots, F_d)$ of $\B^n$ (here we may have $F_0=0$), and every 
simplex of $\B^n$ is of this form. Note that this simplicial decomposition of $\B^n$ is the flag 
complex of $(\FF (\AA), \le)$.

Let $F_b \in \FF (\AA)$ be a facet. It can be easily checked that the union of the simplices of the form 
$\Delta (F_0, F_1, \dots, F_d)$ with $F_b=F_0<F_1< \cdots <F_d$ is a closed disk whose dimension is 
equal to $\codim\, F_b$. Its interior is denoted by $U(F_b)$. So, the set $\{ U(F); F \in \FF (\AA)\}$ 
forms a cellular decomposition of $\B^n$ called the {\it dual decomposition}.

\bigskip\noindent
{\bf Example.} Let $\AA$ be a collection of 3 lines in $\R^2$ (see Figure 5.2). The poset $\FF(\AA)$ 
contains 6 chambers, 6 facets of dimension 1 (half-lines), and $0$. The dual decomposition of $\B^2=\D$ 
has 6 vertices, 6 edges, and one 2-cell.

%%%%%%%%%%%
\begin{figure}[htb]
\centerline{
\setlength{\unitlength}{.4cm}
\begin{picture}(8,8)
\put(0,0.5){\includegraphics[width=3.2cm]{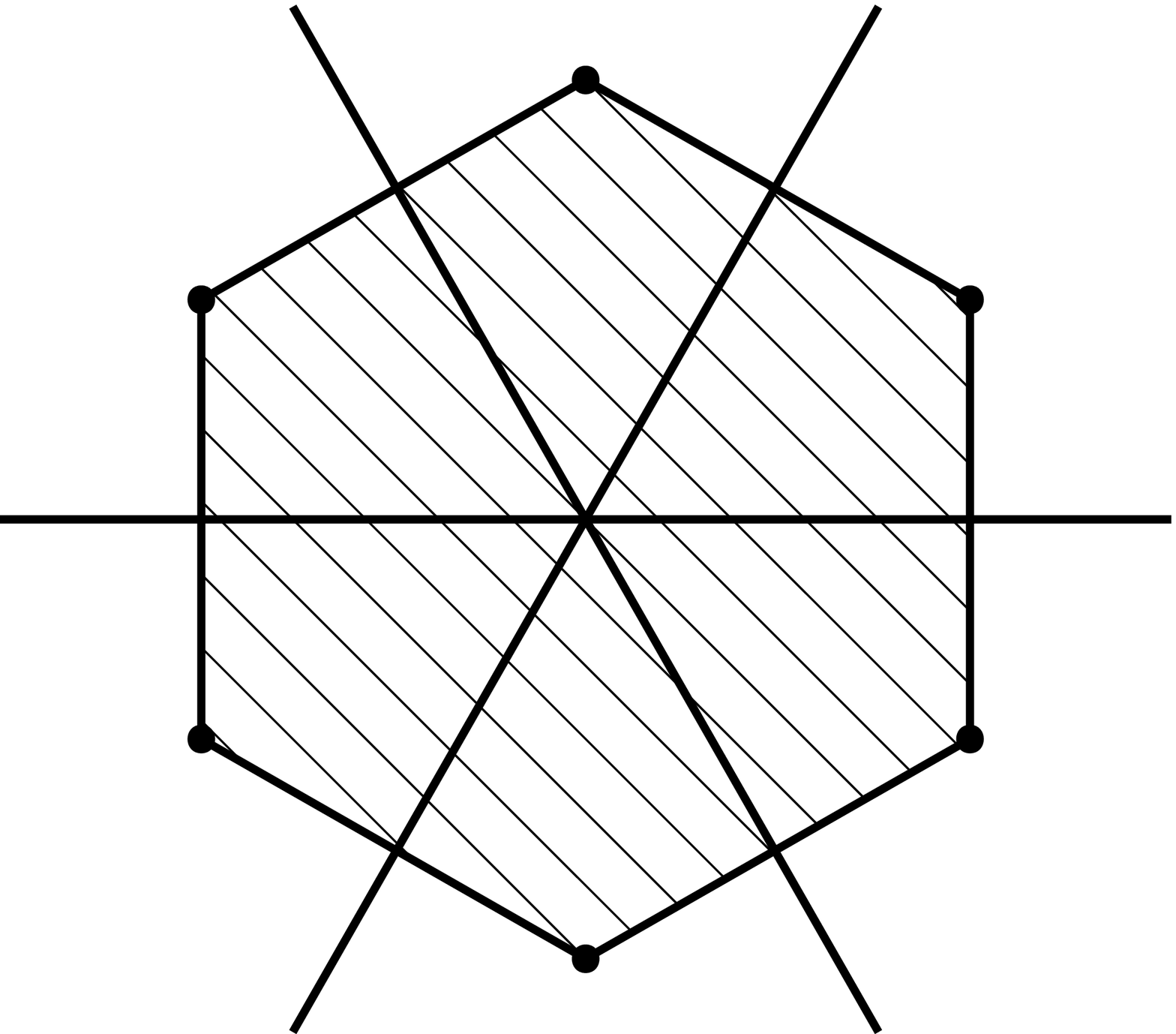}}
\end{picture}}
\medskip
\centerline{{\bf Figure 5.2.} A dual decomposition.}
\end{figure}
%%%%%%%%%

\bigskip\noindent
Let $X_b=(F_b,C_b) \in \XX$. We denote by $\bar U (X_b)$ the union of the simplices $\Delta (X_0, X_1, 
\dots, X_d)$ of $\Sal (\AA)$ such that $X_b=X_0<X_1< \cdots <X_d$. One can show (with some effort)
that, for every $F \ge 
F_b$, there exists a unique chamber $C \in \CC (\AA)$ such that $F \le C$ and $(F_b,C_b) \le (F,C)$. 
This implies that $\bar U(X_b)$ is homeomorphic to $\bar U(F_b)$ via the map $(F,C) \mapsto {\bf 
x}(F)$, $F_b \le F$. Hence, $\bar U(X_b)$ is a closed disk whose dimension is equal to 
$\codim\, F_b$. We denote by $U(X_b)$ the interior of $\bar U(X_b)$. So, $\{ U(X); X \in \XX \}$ forms a 
(regular) cell decomposition of $\Sal (\AA)$.

\bigskip\noindent
{\bf 0-skeleton.} For $C \in \CC(\AA)$, we set $\omega (C)= U(C,C)= \bar U(C,C)$. Then the 0-skeleton 
of $\Sal (\AA)$ is
\[
\Sal_0 (\AA) = \{ \omega (C)\ ;\ C \in \CC (\AA)\}\,.
\]

\bigskip\noindent
{\bf 1-skeleton.} Let $F \in \FF (\AA)$ be a facet of codimension 1. There are exactly two chambers 
$C,D \in \CC (\AA)$ such that $F \le C$ and $F \le D$. Then there are two edges, $U(F,C)$ and $U(F,D)$, 
joining $\omega(C)$ and $\omega(D)$ in the 1-skeleton of $\Sal (\AA)$ (see Figure 5.3). We use the 
convention that $U(F,C)$ is endowed with an orientation which goes from $\omega (C)$ to $\omega (D)$.

%%%%%%%%%%%
\begin{figure}[htb]
\centerline{
\setlength{\unitlength}{.4cm}
\begin{picture}(16,5)
\put(2,0){\includegraphics[width=4.4cm]{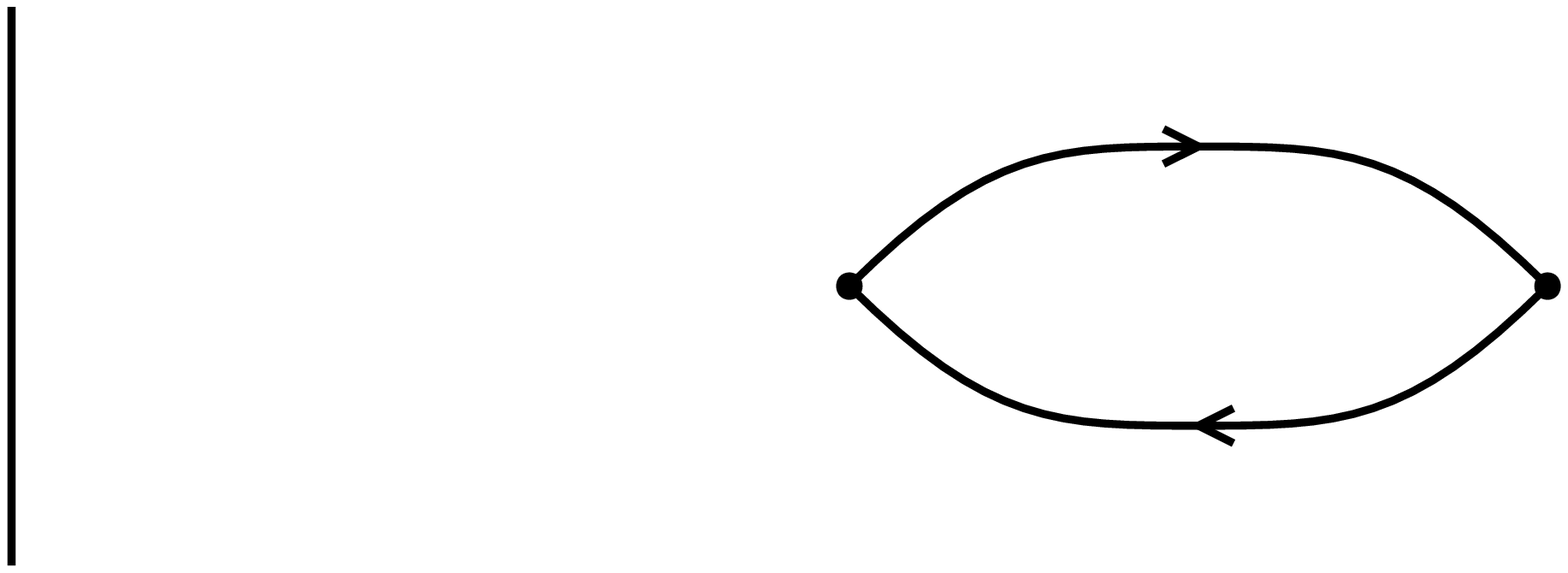}}
\put(0.4,1.8){\small $C$}
\put(1.7,4.5){\small $F$}
\put(3,1.8){\small $D$}
\put(6,1){\small $\omega (C)$}
\put(9,0){\small $U(F,D)$}
\put(9,3.3){\small $U(F,C)$}
\put(13,1){\small $\omega(D)$}
\end{picture}}
\medskip
\centerline{{\bf Figure 5.3.} 1-skeleton of $\Sal (\AA)$.}
\end{figure}
%%%%%%%%%

\bigskip\noindent
{\bf 2-skeleton.} Let $F_b \in \FF(\AA)$ be a facet of codimension 2, and let $C_b \in \CC (\AA)$ such 
that $F_b \le C_b$. Let $C_0=D_0=C_b, C_1, \dots, C_l=D_l, \dots, D_1$ be the chambers $C \in \CC 
(\AA)$ such that $F_b\le C$, arranged like in Figure 5.4. Let $F_1, \dots, F_l, G_1, \dots, G_l$ be 
the facets $F \in \FF (\AA)$ of codimension 1 such that $F_b \le F$, arranged like in Figure 5.4. Set 
$a_i=U(F_i,C_{i-1})$ and $b_i= U(G_i,D_{i-1})$ for $1 \le i\le l$. Then $U(F_b,C_b)$ is a 2-disk whose 
boundary is $(a_1a_2 \cdots a_l)(b_1b_2 \cdots b_l)^{-1}$.

%%%%%%%%%%%
\begin{figure}[htb]
\centerline{
\setlength{\unitlength}{.4cm}
\begin{picture}(29,9)
\put(2,0.5){\includegraphics[width=8.8cm]{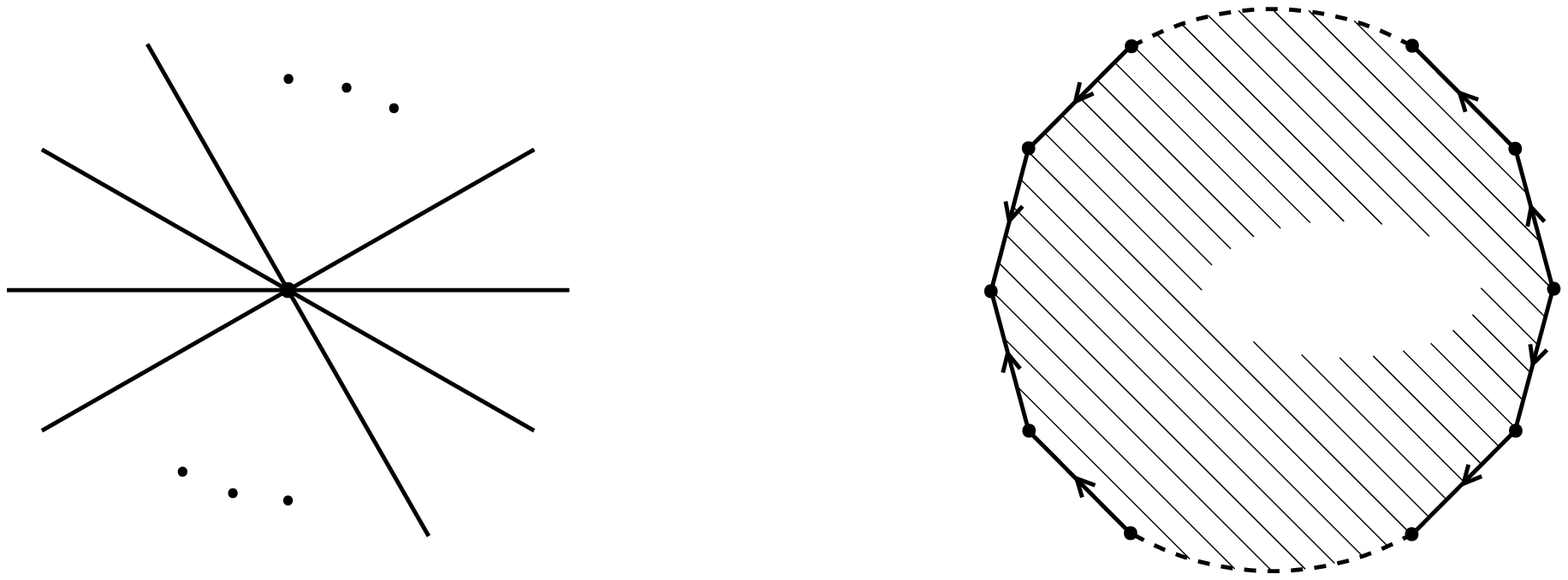}}
\put(5.8,5){\small $F_b$}
\put(9.6,1.9){\small $F_1$}
\put(8.2,0.4){\small $F_2$}
\put(1.5,1.7){\small $F_{l-1}$}
\put(0.7,4.3){\small $F_l$}
\put(10.5,4.3){\small $G_1$}
\put(10,6.6){\small $G_2$}
\put(3.2,8.2){\small $G_{l-1}$}
\put(1.3,6.3){\small $G_l$}
\put(8.5,3.6){\small $C_b$}
\put(7.6,2){\small $C_1$}
\put(2.4,3.5){\small $C_{l-1}$}
\put(2.5,4.9){\small $C_l$}
\put(8.5,4.9){\small $D_1$}
\put(3.1,6.3){\small $D_{l-1}$}
\put(24.3,4.5){\small $\omega(C_b)$}
\put(13.5,4.4){\small $\omega(C_l)$}
\put(24.1,3.5){\small $a_1$}
\put(23,1.5){\small $a_2$}
\put(16,1.1){\small $a_{l-1}$}
\put(15.1,3.3){\small $a_l$}
\put(24,5.4){\small $b_1$}
\put(22.8,7.2){\small $b_2$}
\put(16,7.4){\small $b_{l-1}$}
\put(15.1,5.2){\small $b_l$}
\put(19.2,4.3){\small $U(F_b,C_b)$}
\end{picture}}
\medskip
\centerline{{\bf Figure 5.4.} 2-skeleton of $\Sal (\AA)$.}
\end{figure}
%%%%%%%%%

\bigskip\noindent
Let $\Gamma$ be a Coxeter graph of spherical type, let $(W_\Gamma, S)$ be the Coxeter system of type 
$\Gamma$, and let $(G_\Gamma, \Sigma)$ be the Artin system of type $\Gamma$. Recall the set $\Pi = \{ 
e_s; s \in S\}$ of simple roots, the linear space $V= \oplus_{s \in S} \R e_s$, and the canonical 
bilinear form $\langle, \rangle : V \times V \to \R$. Recall also from Theorem 3.8 that $\langle, 
\rangle$ is positive definite, and that $W = W_\Gamma$ can be viewed as a finite subgroup of $O(V)= 
O(V, \langle, \rangle)$ generated by reflections.

Let $\AA_\Gamma$ denote the set of reflecting hyperplanes of $W$. Then $M_\Gamma = M(\AA_\Gamma)$, the 
group $W_\Gamma$ acts freely on $M_\Gamma$, $N_\Gamma = M_\Gamma /W_\Gamma$, and $\pi_1 (N_\Gamma) = 
G_\Gamma$ (see Subsection 3.4).

Fix a (base) chamber $C_b \in \CC (\AA_\Gamma)$. A hyperplane $H \in \AA_\Gamma$ is called a {\it wall} 
\index{Wall} of $C_b$ if $\codim (\bar C_b \cap H) = 1$. The following is proved in \cite{Bourb1}.

\begin{proposition}
\begin{enumerate}
\item
$C_b$ is a simplicial cone.
\item
Let $H_1, \dots, H_n$ be the walls of $C_b$, and, for $1 \le i \le n$, let $s_i$ be the orthogonal 
reflection with respect to $H_i$. Then, up to conjugation, $S = \{ s_1, \dots, s_n \}$ is the Coxeter 
generating set of $W$.
\end{enumerate}
\end{proposition}

For $T \subset S$ we denote by $W_T$ the subgroup of $W$ generated by $T$, and by $\Gamma_T$ the full 
subgraph of $\Gamma$ spanned by $T$. It is a well-know fact (see \cite{Bourb1}, for example) that 
$(W_T,T)$ is the Coxeter system of type $\Gamma_T$. The {\it Coxeter complex}\index{Coxeter complex}
of $(W,S)$ is defined to 
be the set
\[
\Cox_\Gamma = \{ w W_T \ ;\ T \subset S \text{ and } w \in W\}
\]
ordered by the reverse inclusion ({\it i.e.} $w_1 W_{T_1} \le w_2 W_{T_2}$ if $w_1 W_{T_1} \supset w_2 
W_{T_2}$).

We fix a base chamber $C_b$ and we take $S = \{ s_1, \dots, s_n\}$ like in Proposition~5.9. For each $s 
\in S$ we denote by $H_s$ the hyperplane fixed by $s$. So, $\{ H_s; s \in S\}$ is the set of walls of 
$C_b$. Since $C_b$ is a simplicial cone, for every $T \subset S$ there exists a unique facet $F(T) \in 
\FF (\AA_\Gamma)$ such that $F(T) \le C_b$ and $\langle F(T) \rangle =\cap_{s \in T} H_s$. The proof of 
the following can be found in \cite{Bourb1}.

\begin{proposition}
The map
\[
\begin{array}{rccc}
\psi:& \Cox_\Gamma &\to &\FF (\AA_\Gamma)\\
&wW_T &\mapsto &w F(T)
\end{array}
\]
is well-defined and is an isomorphism of ordered sets.
\end{proposition}

Now, the following lemmas 5.11 and 5.12 are used to describe the poset $\XX$ in terms of Coxeter 
complexes.

\begin{lemma}[Bourbaki \cite{Bourb1}]
Let $T \subset S$ and $w \in W$. Then $wW_T$ has a 
smallest element $u$ for the order $\le_L$ (defined in Subsection 3.2). That is, for all $w' \in wW_T$ 
there exists a unique $v' \in W_T$ such that $w' = uv'$ and $\lg_S (w') = \lg_S (u) + \lg_S (v')$.
\end{lemma}

The smallest element of $w W_T$ is denoted by $u=\min_T(w)$, and such an element is called {\it $T$-
minimal}. The set of $T$-minimal elements is denoted by $\Min (T)$. For $w \in W$, we denote by 
$\pi_T(w)$ the element $v \in W_T$ such that $w = \min_T(w) \cdot v$.

The proof of the following is left to the reader.

\begin{lemma}
Let $C_b$ be a base chamber, let $T \subset S$, and let $F=F(T)$. Let $w_1, w_2 
\in W$. We have $(w_1 C_b)_F = (w_2C_b)_F$ if and only if $\pi_T (w_1) = \pi_T (w_2)$.
\end{lemma}

Set
\[
\widehat{\Cox}_\Gamma = \{ (T,w)\ ;\ w \in W \text{ and } T \subset S\}\,.
\]
Let $\le$ be the partial order on $\widehat{\Cox}_\Gamma$ defined by 
\begin{multline*}
(T_1,w_1) \le (T_2,w_2) \quad \text{if} \quad T_1\supset T_2\,,\ \min_{T_1} (w_1) = \min_{T_1} (w_2)\,,\\ 
\text{ and } \pi_{T_2} (w_1) = \pi_{T_2} (w_2)\,.
\end{multline*}
Note that the conditions ``$T_1 \supset T_2$ and $\min_{T_1} (w_1) = \min_{T_1} (w_2)$'' are equivalent to 
the condition $w_1 W_{T_1} \supset w_2 W_{T_2}$, and, by Lemma 5.12, the condition $\pi_{T_2} (w_1) = 
\pi_{T_2} (w_2)$ is equivalent to the condition $(w_1 C_b)_{F(T_2)} = (w_2 C_b)_{F(T_2)}$. So:

\begin{theorem}
The map
\[
\begin{array}{rccc}
\hat \psi : & \widehat{\Cox}_\Gamma & \to & \XX (\AA_\Gamma)\\
& (T,w) & \mapsto & (w F(T), wC_b)
\end{array}
\]
is well-defined and is an isomorphism of posets.
\end{theorem}

For $(T, w) \in \widehat{\Cox}_\Gamma$ we set $U (T,w) = U (\hat \psi (T,w))$. So, $\{ U(T,w); (T,w) 
\in \widehat{\Cox}_\Gamma\}$ is a cellular decomposition of $\Sal (\AA_\Gamma)$. Moreover, the 
dimension of $U(T,w)$ is $|T|$ for all $(T,w) \in \widehat{\Cox}_\Gamma$.

The Coxeter group $W$ acts on $\widehat{\Cox}_\Gamma$ by
\[
u \cdot (T,w) = (T, uw)\,, \quad \text{for } (T,w) \in \widehat{\Cox}_\Gamma \text{ and } u \in W\,.
\]
It turns out that this action preserves the order of $\widehat{\Cox}_\Gamma$ and induces a cellular 
action on $\Sal (\AA_\Gamma)$ defined by
\[
u \cdot U (T,w) = U(T,uw) \quad \text{for } (T,w) \in \widehat{\Cox}_\Gamma \text{ and } u \in W\,.
\]

\begin{theorem}[Salvetti \cite{Salve1}]
There exists an embedding $\Sal (\AA_\Gamma) 
\hookrightarrow M_\Gamma$ and a 
(strong) retracting deformation of $M_\Gamma$ onto $\Sal (\AA_\Gamma)$ that are equivariant under the 
action of $W$. In particular, there exists an embedding $\Sal(\AA_\Gamma) /W \hookrightarrow M_\Gamma/W 
= N_\Gamma$ and a (strong) retracting deformation of $N_\Gamma$ onto $\Sal (\AA_\Gamma)/W$.
\end{theorem}

To each $T \subset S$ corresponds a unique cell $U_N(T)$ of $\Sal (\AA_\Gamma)/W$ of dimension $|T|$. 
This cell is the orbit of $U(T,w)$ for all $w \in W$. Every cell of $\Sal (\AA_\Gamma)/W$ is of this 
form.

The 0-skeleton of $\Sal (\AA_\Gamma)/W$ contains a unique vertex, $\omega_N = U_N(\emptyset)$. For every 
$s \in S$ there is an edge $U_N(s)$ in $\Sal(\AA_\Gamma)/W$ and each edge is of this form. For every 
pair $\{s,t\} \subset S$ there is a 2-cell $U_N (s,t)$ in $\Sal (\AA_\Gamma)/W$ whose boundary is
\[
\pprod( U_N(s), U_N(t): m_{s\,t}) \cdot \pprod (U_N(t), U_N(s): m_{s\,t})^{-1}\,,
\]
and every 2-cell is of this form. Note that the 2-skeleton of $\Sal(\AA_\Gamma)/W$ is equal to the 
2-cell complex associated to the standard presentation of $G_\Gamma$. This gives an alternative proof to 
Theorems 2.2 and 3.14.

For $0 \le q \le |S|$, set
\[
C_q (G_\Gamma) = \bigoplus_{\begin{subarray}{c}
T\subset S\\
|T| =q
\end{subarray}}
\Z [G_\Gamma] \cdot E_T\,,
\]
the free $\Z [G_\Gamma]$-module freely spanned by $\{ E_T; T \subset S \text{ and } |T|=q\}$. We fix a 
total order $S= \{ s_1, \dots, s_n \}$ on $S$ and we define $d: C_q (G_\Gamma) \to C_{q-1} (G_\Gamma)$ 
as follows. Let $T= \{ s_{i_1}, \dots, s_{i_q} \} \subset S$, $i_1 < \cdots < i_q$. Then
\[
dE_T = \sum_{j=1}^q (-1)^{j-1} \left( \sum_{\begin{subarray}{c}
u \in W_T\\
u \in \Min (T \setminus \{ s_{i_j}\})
\end{subarray}}
(-1)^{\lg_S(u)} \kappa (u) \right) \cdot E_{T \setminus \{s_{i_j}\}} \,,
\]
where $\kappa : W \to G_\Gamma$ is the set-section of the canonical epimorphism $\theta: G_\Gamma \to W$ 
defined in Subsection 3.3.

\begin{theorem}[De Concini, Salvetti \cite{DeCSal1}, Squier \cite{Squie1}]
The complex 
\linebreak
$(C_\ast (G_\Gamma), d)$ is a free 
resolution of $\Z$ by $\Z [G_\Gamma]$-modules.
\end{theorem}

\begin{note}
Squier's proof of Theorem 5.15 does not use the Salvetti complexes at all 
and is independent from the proof of De Concini and Salvetti.
\end{note}

%%%%%%%%%%%%%%%%%%%%%%%%%%%%%%%%%

\section{Linear representations}

The existence (or non-existence) of faithful linear representations of the braid groups was one of the 
major problem in the field. This problem was solved by Bigelow \cite{Bigel1} and Krammer \cite{Kramm3} 
in 2000. They representation, which is known now as the LKB representation\index{LKB representation}, 
was right afterwards 
extended to the Artin groups of type $D_n$ ($n \ge 4)$ and $E_k$ ($k=6,7,8$) by Digne \cite{Digne2}, 
Cohen, and Wales \cite{CohWal1}, and to all Artin groups of small type in \cite{Paris1}. The 
representations of Digne, Cohen and Wales are proved to be faithful. Hence, since any spherical type 
Artin group embeds in a direct product of Artin groups of type $A_n$ ($n \ge 1$), $D_n$ ($n \ge 4$), 
and $E_k$ ($k=6,7,8$) (see \cite{Crisp1}), any Artin group of spherical type is linear. The extension 
to the non-spherical type Artin groups gives rise to a linear representation over an infinite 
dimensional vector space, so it cannot be used for proving that these groups are linear. However, these 
representations are useful tools to study the non-spherical type Artin groups. In particular, they are 
the main tool in the proof of Theorem 3.9.

In Subsection 6.1 we present the algebraic approach to the LKB representations as constructed in 
\cite{Paris1} for the Artin groups of small type. Subsection~6.2 is dedicated to the topological 
construction of the LKB representations. Curiously, this topological point of view is known only for the braid 
groups.

\subsection{Algebraic approach}

Let $\Gamma$ be a Coxeter graph, let $M=(m_{s\,t})_{s,t\in S}$ be the Coxeter matrix of $\Gamma$, let 
$(W_\Gamma,S)$ be the Coxeter system of type $\Gamma$, let $(G_\Gamma, \Sigma)$ be the Artin system of 
type $\Gamma$, and let $G_\Gamma^+$ be the Artin monoid of type $\Gamma$.

We say that $\Gamma$ is of {\it small type}\index{Small type Coxeter graph}
if $m_{s\,t} \le 3$ for all $s,t \in S$, $s \neq t$, and we 
say that $\Gamma$ is {\it without triangle}\index{Coxeter graph without triangle}
if there is no triple $\{ s,t,r\}$ in $S$ such that 
$m_{s\,t} = m_{t\,r} = m_{r\,s} = 3$. We assume from now on that $\Gamma$ is of small type and without 
triangle.

Recall from Subsection 3.2 the set $\Pi= \{ e_s; s \in S\}$ of simple roots, the space $V = \oplus_{s 
\in S} \R e_s$, the canonical bilinear form $\langle, \rangle : V \times V \to \R$, and the root system 
$\Phi= \{ we_s; s \in S \text{ and } w \in W\}$. Recall also that we have the disjoint union $\Phi = 
\Phi_+ \sqcup \Phi_-$, where $\Phi_+$ is the set of positive roots and $\Phi_-$ is the set of negative 
roots (see Proposition 3.4).

Set $\EE = \{ u_f; f \in \Phi_+ \}$ an abstract set in one-to-one correspondence with $\Phi_+$, and $\K 
= \Q (x,y)$. Note that $\EE$ is finite if and only if $\Gamma$ is of spherical type. We denote by $\VV$ 
the $\K$-vector space having $\EE$ as a basis.

For all $s \in S$ we define a linear transformation $\varphi_s : \VV \to \VV$ by
\[
\varphi_s (u_f) = \left\{
\begin{array}{cl}
0&\quad\text{if } f=e_s\\
u_f&\quad\text{if }\langle e_s,f \rangle =0\\
y \cdot u_{f-ae_s}&\quad\text{if } \langle e_s,f \rangle = a >0 \text{ and } f \neq e_s\\
(1-y) \cdot u_f + u_{f+ae_s} &\quad\text{if } \langle e_s,f \rangle = -a<0\\
\end{array}\right.
\] 
The following is easy to prove.

\begin{lemma}
The mapping $\sigma_s \mapsto \varphi_s$, $s \in S$, induces a homomorphism of 
monoids $\varphi: G_\Gamma^+ \to \End (\VV)$.
\end{lemma}

For all $s \in S$ and all $f \in \Phi_+$ we choose a polynomial $T(s,f) \in \Q [y]$ and we define 
$\Phi_s: \VV \to \VV$ by
\[
\Phi_s (u_f) = \varphi_s (u_f) + x \cdot T(s,f) \cdot u_{e_s}\,.
\]
Now, we have:

\begin{theorem}[Paris \cite{Paris1}]
There exists a choice of polynomials $T(s,f)$, $s \in S$ 
and $f \in \Phi_+$, such that the mapping $\sigma_s \mapsto \Phi_s$, $s \in S$, induces a homomorphism 
$\Phi: G_\Gamma^+ \to \GL (\VV)$.
\end{theorem}

\begin{theorem}[Paris \cite{Paris1}]
The above defined homomorphism $\Phi: G_\Gamma^+ \to \GL 
(\VV)$ is injective.
\end{theorem}

\begin{corollary}[Paris \cite{Paris1}]
The natural homomorphism $\iota: G_\Gamma^+ \to 
G_\Gamma$ is injective.
\end{corollary}

\begin{proof}
Since $G_\Gamma$ is the group of fractions of $G_\Gamma^+$, there exists a unique 
homomorphism $\hat \Phi: G_\Gamma \to \GL (\VV)$ such that $\Phi= \hat \Phi \circ \iota$. Since $\Phi$ 
is injective, we conclude that $\iota$ is also injective.
\end{proof}

\begin{corollary}[Bigelow \cite{Bigel1}, Krammer \cite{Kramm3}, Digne \cite{Digne2}, Cohen, Wales 
\cite{CohWal1}]
Suppose that $\Gamma$ is of spherical type. Let $\hat \Phi: G_\Gamma \to \GL 
(\VV)$ be the homomorphism induced by $\Phi$. Then $\hat \Phi$ is injective.
\end{corollary}

\begin{proof}
Let $\alpha \in \Ker\, \hat\Phi$. By Proposition 4.8, $\alpha$ can be written in the form 
$\alpha = \beta^{-1} \gamma$, with $\beta, \gamma \in G_\Gamma^+$. We have $1 = \hat \Phi (\alpha) = 
\Phi( \beta)^{-1} \Phi(\gamma)$, thus $\Phi (\beta) = \Phi(\gamma)$. Since $\Phi$ is injective, it 
follows that $\beta = \gamma$, thus $\alpha = \beta^{-1} \gamma =1$.
\end{proof}

\begin{note}
It is shown in \cite{Paris1} that any Artin monoid $G_\Gamma^+$ can be embedded in an 
Artin monoid $G_\Omega^+$, where $\Omega$ is of small type without triangle. Moreover, if $\Gamma$ is 
of spherical type, then $\Omega$ can be chosen to be of spherical type (see also \cite{Crisp1}, 
\cite{Godel1}, \cite{CriPar3}, \cite{Caste1}). So, Corollary 6.4 implies that $\iota: G_\Gamma^+ \to 
G_\Gamma$ is injective for all Coxeter graphs $\Gamma$, and Corollary 6.5 implies that all the Artin 
groups of spherical type are linear.
\end{note}

\begin{note}
It is shown in \cite{Marin1} that: if $\Gamma$ is of type $A_n$, $D_n$, $E_k$ ($k=6,7,8$), 
then the image of $\hat \Phi$ is Zariski dense in $\GL(\VV)$. In particular, this shows that $\hat 
\Phi$ is irreducible (see also \cite{Zinno1}, \cite{Marin2}, \cite{CoGiWa1}).
\end{note}

\begin{note}
The proof of Theorem 6.3 given in \cite{Paris1} is largely inspired by Krammer's proof of 
the same theorem for the braid groups \cite{Kramm3}. A new, short, and elegant proof can be found now
in \cite{Hee1}.
\end{note}

\subsection{Topological approach}

Now, we give a topological interpretation of the representation $\hat \Phi: G_\Gamma \to \GL (\VV)$ in 
the case $\Gamma= A_{n-1}$, that is, when $G_\Gamma = \BB_n$ is the braid group on $n$ strands. Such an 
interpretation is unknown for the other Artin groups.

Let $M$ be a connected CW-complex, let $G= \pi_1 (M)$, and let $R$ be a (right) $\Z [G]$-module. Let 
$\tilde M$ be the universal cover of $M$. The action of $G$ on $\tilde M$ induces an action of $G$ on 
the group $C_q (\tilde M)$ of (cellular) $q$-chains of $\tilde M$, and this action makes $C_q (\tilde 
M)$ a module over the group ring $\Z [G]$. It is also easily seen that the boundary maps $\partial : 
C_q (\tilde M) \to C_{q-1} (\tilde M)$ are $\Z[G]$-module homomorphisms. We define $C_q(M,R)$ to be $R 
\otimes_{\Z[G]} C_q (\tilde M)$. These groups form a chain complex with boundary map $\Id \otimes 
\partial$. The homology groups $H_q (M,R)$ of this chain complex are the {\it homology 
groups of $M$ with local coefficients $R$}\index{Homology with local coefficients}.

Now, for $n \ge 1$, $M_n$ denotes the space of ordered configurations of $n$ points in $\C$, and $N_n= 
M_n/ \Sym_n$ denotes the space of (unordered) configurations of $n$ points in $\C$ (see Section 2). Let 
$n,m \ge 2$. By \cite{FadNeu1}, the map
\[
\begin{array}{rccc}
p_{n,m}: & M_{n+m} &\to& M_n\\
&(z_1, \dots z_n,z_{n+1}, \dots, z_{n+m}) & \mapsto &(z_1, \dots, z_n)
\end{array}
\]
is a locally trivial fiber bundle which admits a cross-section. The fiber of $p_{n,m}$ is as follows. 
Set
\[
\begin{array}{rcll}
H_{i\,j}&=&\{ w \in \C^m\ ;\ w_i=w_j\}&\quad\text{for } 1 \le i<j\le m\\
K_{i\,k}&=&\{  w \in \C^m\ ;\ w_i=k\} &\quad\text{for }1 \le i\le m \text{ and } 1 \le k\le n
\end{array}
\]
Set
\[
X_{n,m} = \C^m \setminus \left( \left( \bigcup_{i<j} H_{i\,j} \right) \cup \left( \bigcup_{\substack{1 
\le i\le m\\ 1 \le k\le n}} K_{i\,k} \right) \right)\,.
\]
Then
\[
p_{n,m}^{-1} (1,2, \dots, n) = \{ (1,2, \dots, n) \} \times X_{n,m}\,.
\]
Let $\Sym_n \times \Sym_m$ act on $M_{n+m}$, $\Sym_n$ acting by permutations on the $n$ first 
coordinates, and $\Sym_m$ acting on the $m$ last ones. Set
\begin{gather*}
N_{n,m} = M_{n+m}/(\Sym_n \times \Sym_m)\,,\\
Y_{n,m} = X_{n,m}/\Sym_m\,.
\end{gather*}
Then $p_{n,m}$ induces a locally trivial fiber bundle $\bar p_{n,m} : N_{n,m} \to N_n$ whose fiber is 
$Y_{n,m}$.

For $ z \in \C^n$ we set
\[
\|  z \|_\infty = \max \{ |z_i|\ ;\ 1 \le i\le n\}\,.
\]
It is easily checked that the map
\[
\begin{array}{rccc}
\kappa: &M_n &\to& M_{n+m}\\
& z &\mapsto &( z, \|  z\|_\infty +1,\| z\|_\infty +2, \dots, \|  z\|_\infty +m)
\end{array}
\]
is a well-defined cross-section of $p_{n,m}$ which is equivariant by the action of $\Sym_n$, thus it 
induces a cross-section $\bar \kappa: N_n \to N_{n,m}$ of $\bar p_{n,m}$. By the homotopy long exact 
sequence of a fiber bundle (see Theorem 2.9), we conclude that $\pi_1 (N_{n,m})$ can be written as a 
semi-direct product $\pi_1 (N_{n,m}) = \pi_1 (Y_{n,m}) \rtimes \BB_n$.

Set $G_{n,m} = \pi_1 (Y_{n,m})$. We consider $G_{n,m}$ as a subgroup of $\pi_1( N_{n,m})$ which, in its turn, is  
viewed as a subgroup of $\pi_1(N_{n+m}) = \BB_{n+m}$. It is easily seen that $G_{n,m}$ 
is generated by the set
\[
\{ \sigma_k\ ;\ n+1 \le k\le n+m\} \cup \{ \delta_{i\,k}\ ;\ 1 \le i\le n \text{ and } n+1 \le k\le 
n+m\}\,,
\]
where $\delta_{i\,k}$ is the pure braid defined in Theorem 2.3. Let $b$ be the homology class of 
$\sigma_{n+1}$ in $H_1 (G_{n,m}) = H_1 (Y_{n,m})$, and let $a_i$ be the homology class of 
$\delta_{i,n+1}$, $1 \le i\le n$. The proof of the following is left to the reader.

\begin{proposition}
$H_1 (Y_{n,m}) = H_1 (G_{n,m})$ is a free abelian group freely generated by 
$\{b,a_1, a_2, \dots, a_n \}$.
\end{proposition}

Let $\bar \rho : H_1 (G_{n,m}) \to \Q (x,y)^\ast$ be the homomorphism which sends $a_i$ to $x$ for all 
$1 \le i\le n$, and sends $b$ to $y$. Let $\rho: G_{n,m} \to \Q (x,y)^\ast$ be the composition of the 
natural projection $G_{n,m} \to H_1 (G_{n,m})$ with $\bar \rho$. This homomorphism makes $\Q (x,y)$ a 
$\Z [G_{n,m}]$-module that we denote by $\Gamma_\rho$.

The proof of the following is also left to the reader.

\begin{proposition}
The kernel of $\rho$ is invariant under the action of $\BB_n$, and $\BB_n$ 
acts trivially on the quotient $G_{n,m}/ \Ker\, \rho \simeq \Z \times \Z$.
\end{proposition}

From Proposition 6.7 follows that the fibration $\bar p_{n,m} : N_{n,m} \to N_n$ induces a monodromy 
representation $\Phi_{n,m}: \BB_n \to \Aut_{\Q (x,y)} (H_\ast (Y_{n,m}, \Gamma_\rho))$.

The following was announced by Krammer \cite{Kramm4}, \cite{Kramm3}, and proved in \cite{Bigel1} (see 
also \cite{PaoPar1}).

\begin{theorem}[Bigelow \cite{Bigel1}]
The homomorphism $\Phi_{n,2}: \BB_n \to 
\linebreak
\Aut_{\Q(x,y)} 
(H_2 (Y_{n,2},\Gamma_\rho))$ coincides with the representation $\hat \Phi : G_{A_{n-1}} \to \GL 
(\VV)$ defined in Subsection 6.1.
\end{theorem}

\begin{note}
The representation $\hat \Phi: G_{A_{n-1}} \to \GL (\VV)$ also coincides with the 
representation studied by Lawrence in \cite{Lawre1}. Lawrence's construction is also geometric. It 
slightly differs from the one presented above, but I do not know exactly how to relate them without the 
formulas.
\end{note}

\begin{note}
It is announced in \cite{Zheng1} that $\Phi_{n,m}$ is faithful for all $m\ge 2$, and it is 
announced in \cite{Chen1} that $\Phi_{n,m}: \BB_n \to \Aut_{\Q (x,y)} (H_m (Y_{n,m}, \Gamma_\rho))$ 
is irreducible for all $m \ge 2$.
\end{note}

%%%%%%%%%%%%%%%%%%%%%%%%%%%%%%%%%

\section{Geometric representations}

\subsection{Definitions and examples}

Let $\Sigma$ be an oriented compact surface, possibly with boundary, and let $\PP$ be a finite 
collection of punctures in the interior of $\Sigma$. Let $\MM (\Sigma, \PP)$ denote the mapping class 
group of the pair $(\Sigma, \PP)$, as defined in Subsection 2.3. Let $\Gamma$ be a Coxeter graph, and 
let $G_\Gamma$ be the Artin group of type $\Gamma$. Define a 
{\it geometric representation}\index{Geometric representation} of 
$G_\Gamma$ in $\MM (\Sigma, \PP)$ to be a homomorphism from $G_\Gamma$ to $\MM (\Sigma, \PP)$.

The main tools for constructing geometric representations of Artin groups are the Dehn twists and the 
braid twists. The 
braid twists are defined in Subsection 2.3, and the Dehn twists are defined as follows.

An {\it essential circle}\index{Essential circle} is an embedding $a: \S^1 \hookrightarrow \Sigma \setminus \PP$ of 
the circle whose image is contained in the interior of $\Sigma$ and does not bound any disk in $\Sigma$ 
containing $0$ or $1$ puncture. Two essential circles $a,a'$ are {\it isotopic} if there exists a 
continuous family $\{ a_t\}_{t \in [0,1]}$ of essential circles such that $a=a_0$ and $a'=a_1$. Isotopy 
of essential circles is an equivalence relation that we denote by $a \sim a'$.

Let $a: \S^1 \to \Sigma \setminus \PP$ be an essential circle. Take an embedding $A: [0,1] \times \S^1 
\to \Sigma \setminus \PP$ of the annulus such that $A(\frac{1}{2}, z)= a(z)$ for all $z \in \S^1$, and 
define $T \in \Homeo^+(\Sigma, \PP)$ by
\[
(T \circ A) (t,z) = A(t, e^{2i\pi t} z)\,,
\]
and $T$ is the identity outside the image of $A$ (see Figure 7.1). The {\it Dehn twist}\index{Dehn twist} along $a$, 
denoted by $\sigma_a$, is defined to be the element of $\MM (\Sigma, \PP)$ represented by $T$. Note 
that
\begin{itemize}
\item
the definition of $\sigma_a$ does not depend on the choice of the map $A$;
\item
if $a$ is isotopic to $a'$, then $\sigma_a = \sigma_{a'}$.
\end{itemize}

%%%%%%%%%%%%%%%%%%%
\begin{figure}[htb]
\centerline{
\setlength{\unitlength}{.4cm}
\begin{picture}(20,8)
\put(0,0){\includegraphics[width=8cm]{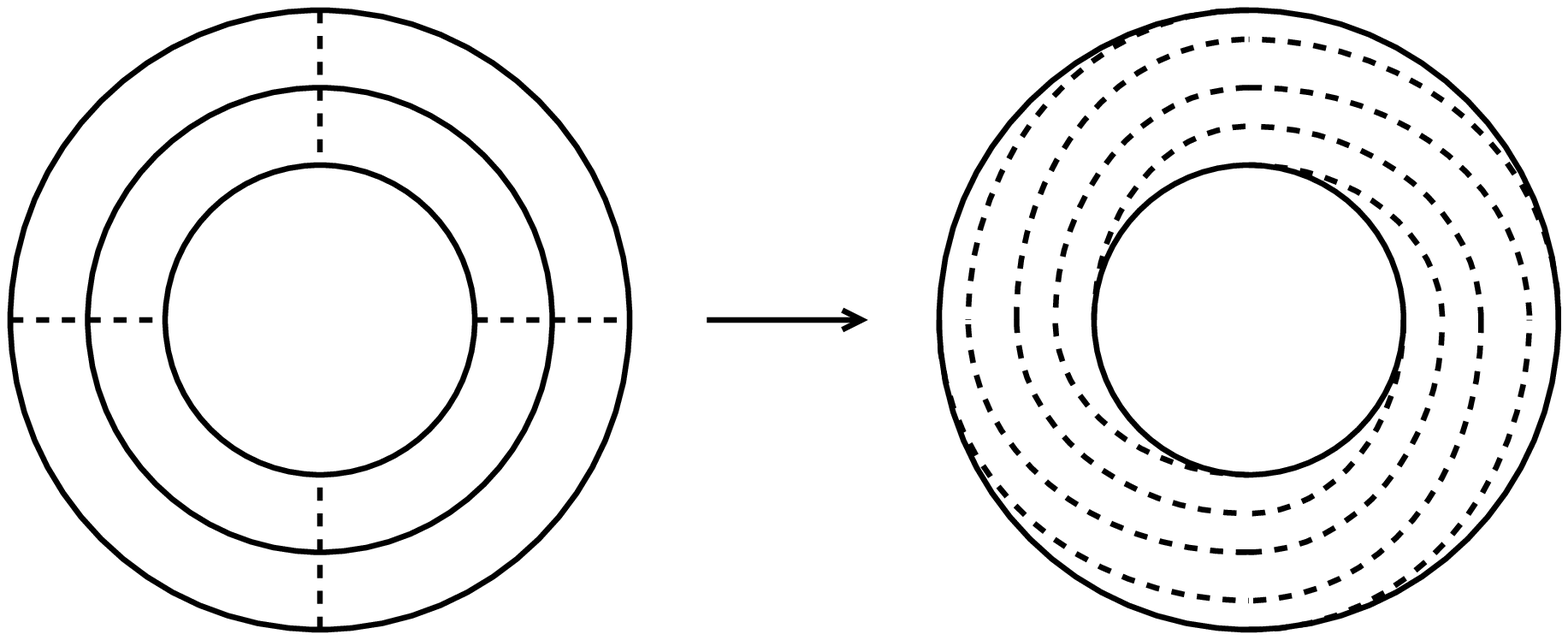}}
\put(1.5,6.4){\small $a$}
\put(9.8,4.3){$T$}
\end{picture}}
\medskip
\centerline{{\bf Figure 7.1.} Dehn twist.}
\end{figure}
%%%%%%%%%%%%%%

Recall that, for an essential arc $a$ of $(\Sigma, \PP)$, $\tau_a$ denotes the braid twist along $a$. The Dehn 
twists and the braid twists satisfy the following relations (see \cite{Birma3}, \cite{LabPar1}).

\begin{proposition}
\begin{enumerate}
\item
Let $a,b$ be two essential circles that intersect transversely. Then
\[
\begin{array}{cl}
\sigma_a \sigma_b = \sigma_b \sigma_a &\quad \text{if } a \cap b = \emptyset\\
\sigma_a \sigma_b \sigma_a = \sigma_b \sigma_a \sigma_b &\quad \text{if } | a \cap b| = 1
\end{array}
\]
\item
Let $a,b$ be two essential arcs of $(\Sigma, \PP)$. Then
\[
\begin{array}{cl}
\tau_a \tau_b = \tau_b \tau_a &\quad \text{if }a \cap b = \emptyset\\
\tau_a \tau_b \tau_a = \tau_b \tau_a \tau_b &\quad \text{if } a(0)= b(1) \text{ and } a \cap b = \{ 
a(0)\}
\end{array}
\]
\item
let $a$ be an essential arc, and let $b$ be an essential circle which intersects $a$ transversely. Then
\[
\begin{array}{cl}
\tau_a \sigma_b = \sigma_b \tau_a &\quad\text{if } a \cap b = \emptyset\\
\tau_a \sigma_b \tau_a \sigma_b = \sigma_b \tau_a \sigma_b \tau_a &\quad \text{if } |a \cap b|=1
\end{array}
\]
\end{enumerate}
\end{proposition}

\bigskip\noindent
{\bf Example 1.} Suppose $\Sigma=\D$ is a disk, and $\PP_n = \{ P_1, \dots, P_n \}$ is a collection of 
$n$ punctures in the interior of $\Sigma$. Then the Artin isomorphism $\Phi: \BB_n \to \MM (\D, \PP_n)$ of 
Theorem 2.16 is a geometric representation of $G_{A_{n-1}} = \BB_n$.

\bigskip\noindent
{\bf Example 2.} Let $n \ge 3$. Suppose that, if $n$ is odd, then $\Sigma$ is a surface of genus 
$\frac{n-1}{2}$ with one boundary component, and if $n$ is even, then $\Sigma$ is a surface of genus 
$\frac{n-2}{2}$ with two boundary components. Let $a_1, \dots, a_{n-1}$ be the essential circles of 
$\Sigma$ pictured in Figure 7.2. By Proposition 7.1, the mapping $\sigma_i \mapsto \sigma_{a_i}$, $1 
\le i\le n-1$, induces a representation $\rho_M: \BB_n \to \MM (\Sigma)$ called the {\it monodromy 
representation}\index{Monodromy representation}
 of $\BB_n$. This geometric representation was introduced by Birman and Hilden in 
\cite{BirHil2}, where it is proved that $\rho_M$ is faithful and its image consists on mapping classes arising 
from homeomorphisms symmetric with respect to a hyperelliptic involution (see also \cite{BirHil1}, 
\cite{Ziesc1}, and \cite{MacHar1}). It is also the geometric monodromy of the simple singularity of 
type $A_{n-1}$ (see \cite{PerVan1}). Let $P_0 \in \partial \Sigma$ be a base-point. Then $\rho_M$  
induces a homomorphism $\rho_{M\, \ast}: \BB_n \to \Aut (\pi_1 (\Sigma, P_0))$ which 
turns out to coincide with the homomorphism $\rho_D: \BB_n \to \Aut (F_{n-1})$ defined in Subsection 
3.1 (see \cite{CriPar2}).

%%%%%%%%%%%%%%%%%%%%%%%
\begin{figure}[htb]
\centerline{
\setlength{\unitlength}{.4cm}
\begin{picture}(21,12)
\put(0,0){\includegraphics[width=7cm]{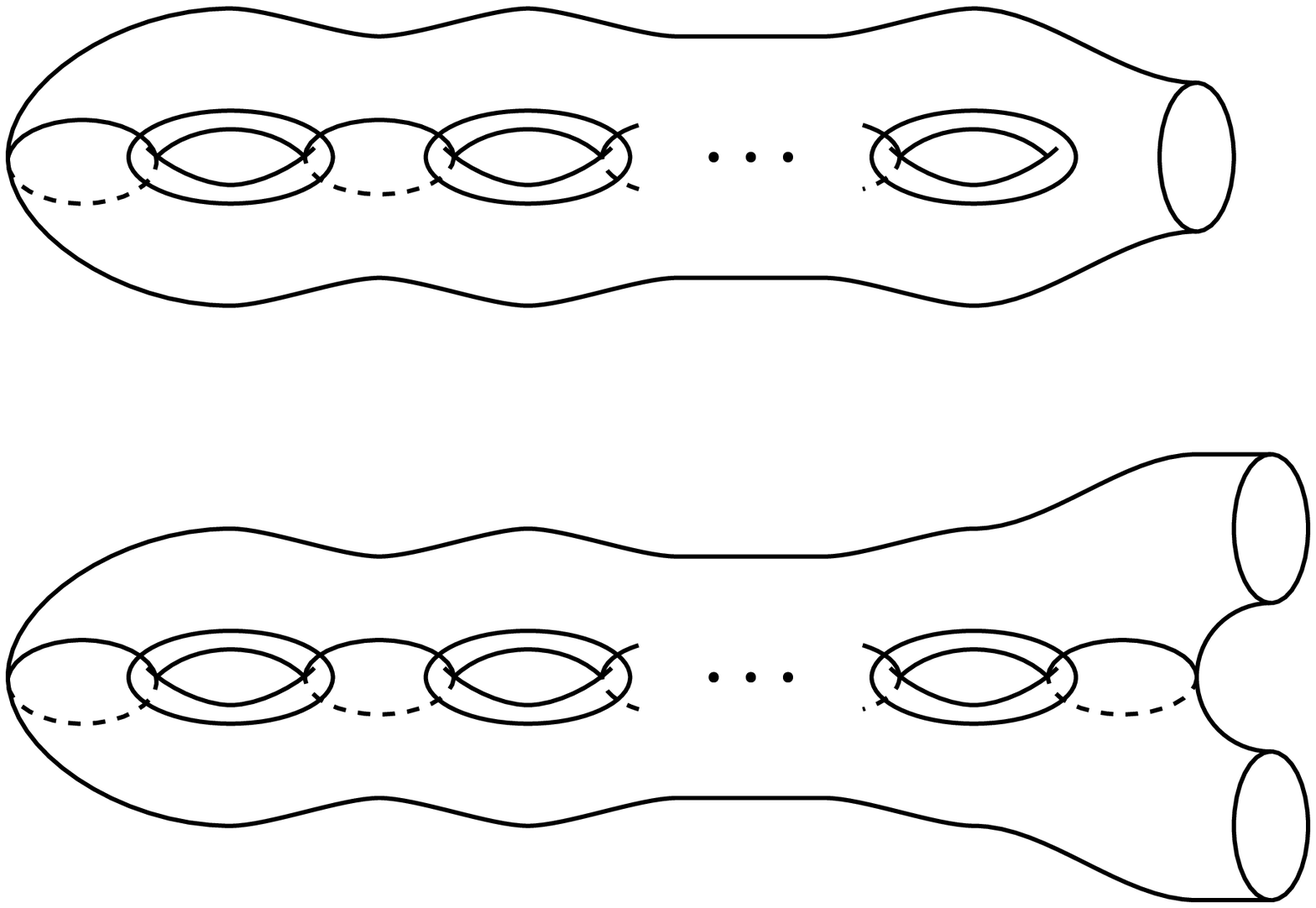}}
\put(0.6,10.6){\small $a_1$}
\put(2.5,10.9){\small $a_2$}
\put(4.5,10.6){\small $a_3$}
\put(6.5,10.9){\small $a_4$}
\put(12.5,10.9){\small $a_{n-1}$}
\put(17,10){\small $n$ odd}
\put(0.6,3.7){\small $a_1$}
\put(2.5,4){\small $a_2$}
\put(4.5,3.7){\small $a_3$}
\put(6.5,4){\small $a_4$}
\put(12.4,4){\small $a_{n-2}$}
\put(14.5,3.7){\small $a_{n-1}$}
\put(18,3){\small $n$ even}
\end{picture}}
\medskip
\centerline{{\bf Figure 7.2.} Monodromy representation of $\BB_n$.}
\end{figure}
%%%%%%%%%%%%%%%%%%%%%%ù

\bigskip\noindent
{\bf Example 3.} Let $\D^2 = \{z \in \C; |z| \le 1 \}$ be the standard disk. A {\it chord diagram}\index{Chord diagram} in 
$\D^2$ is defined to be a collection $\{S_1, \dots, S_n\}$ of segments in $\D^2$ such that
\begin{itemize}
\item
the extremities of $S_i$ belong to $\partial \D^2$ and its interior is contained in the interior of 
$\D^2$, for all $1 \le i\le n$;
\item
either $S_i$ and $S_j$ are disjoint, or they intersect transversely in a unique point in the interior 
of $\D^2$, for all $1 \le i \neq j \le n$.
\end{itemize}
From this data one can define a Coxeter matrix $M=(m_{i\,j})_{1 \le i,j \le n}$ setting $m_{i\,j}=2$ if 
$S_i$ and $S_j$ are disjoint, and $m_{i\,j}=3$ if they intersect. The Coxeter graph $\Gamma$ of $M$ is 
called the {\it intersection diagram}\index{Intersection diagram} of the chord diagram.

From this data one can also define a surface $\Sigma$ by attaching to $\D^2$ a handle $H_i$ which joins 
both extremities of $S_i$, for all $1 \le i \le n$ (see Figure 7.3). Let $a_i$ be the essential circle 
of $\Sigma$ made with $S_i$ and the central arc of $H_i$. Then, by Proposition 7.1, the mapping 
$\sigma_i \mapsto \sigma_{a_i}$, $1 \le i\le n$, induces a geometric representation $\rho_\PV: 
G_\Gamma \to \MM (\Sigma)$, called {\it Perron-Vannier representation}\index{Perron-Vannier representation}.

%%%%%%%%%%%%%%%%%%%%%%%%%%%%
\begin{figure}[htb]
\centerline{
\setlength{\unitlength}{.4cm}
\begin{picture}(18,14)
\put(0,0){\includegraphics[width=7.2cm]{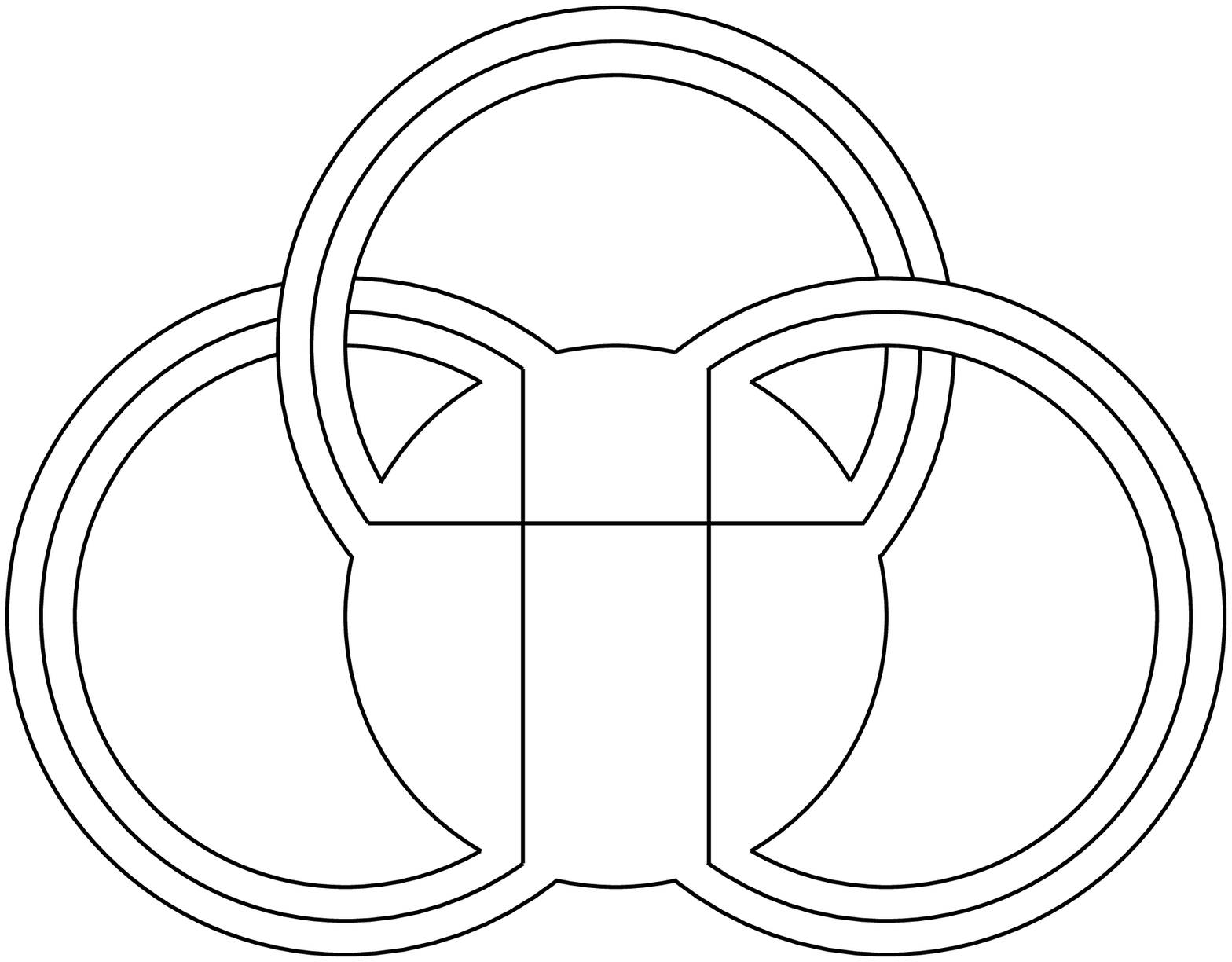}}
\put(1.5,4.5){\small $H_1$}
\put(8.5,12.2){\small $H_2$}
\put(15.5,4.5){\small $H_3$}
\put(8,4){\small $S_1$}
\put(8.5,5.5){\small $S_2$}
\put(9.5,4){\small $S_3$}
\end{picture}}
\medskip
\centerline{{\bf Figure 7.3.} Chord diagram and associated surface.}
\end{figure}
%%%%%%%%%%%%%%%%%%%%%%%

The Perron-Vannier representations were introduced in \cite{PerVan1}. If $\Gamma=A_{n-1}$, then 
$\rho_{\PV}$ is equal to the monodromy representation $\rho_M$ defined in Example 2. More generally, if 
$\Gamma$ is $A_n$ ($n \ge 1$), $D_n$ ($n \ge 4$), or $E_k$ ($k=6,7,8$), then $\rho_{\PV}$ is the 
geometric monodromy of the simple singularity of type $\Gamma$ (see \cite{PerVan1}). For a connected 
graph $\Gamma$, the representation $\rho_{\PV}$ is faithful if and only if either $\Gamma=A_n$ for some 
$n \ge 1$, or $\Gamma = D_n$ for some $n \ge 4$ (see \cite{PerVan1}, \cite{Labru1}, \cite{Wajnr1}).

\bigskip\noindent
{\bf Example 4.} This example comes from \cite{CriPar3}. Recall that a Coxeter graph $\Gamma$ is of 
{\it small type} if $m_{s\,t}\le 3$ for all $s,t \in S$, where $M=(m_{s\,t})_{s,t\in S}$ is the Coxeter 
matrix of $\Gamma$. Let $\Gamma$ be a small type Coxeter graph. We choose (arbitrarily) a total order 
$<$ on $S$. For $s \in S$, we set $\St_s= \{ t \in S; m_{s\,t}=3\} \cup\{s\}$. Write $\St_s= \{t_1, 
t_2, \dots, t_k\}$ such that $t_1<t_2< \cdots <t_k$, and suppose that $s=t_j$. For $1 \le i\le k$, the 
difference $i-j$ is called the {\it relative position} of $t_i$ with respect to $s$ and is denoted by 
$\pos (t_i:s)$. In particular, $\pos (s:s)=0$.

Let $s \in S$ and let $k = |\St_s|$. Let $\An_s$ denote the annulus $\An_s= (\R/2k \Z) \times [0,1]$. 
We define the surface $\Sigma= \Sigma_\Gamma$ by
\[
\Sigma= \left( \bigsqcup_{s \in S} \An_s \right) / \sim\,,
\]
where $\sim$ is the equivalence relation defined as follows. Let $s,t \in S$ such that $s<t$ and 
$m_{s\,t}=3$. Set $p=\pos (t:s) >0$ and $q= \pos (s:t) <0$. For all $(x,y) \in [0,1] \times [0,1]$ the 
relation $\sim$ identifies the point $(2p+x,y)$ of $\An_s$ with the point $(2q+1-y,x)$ of $\An_t$ (see 
Figure 7.4).

%%%%%%%%%%%%%%%%%%%%%%%%%%%
\begin{figure}[htb]
\centerline{
\setlength{\unitlength}{.4cm}
\begin{picture}(29,11)
\put(3,0){\includegraphics[width=10.4cm]{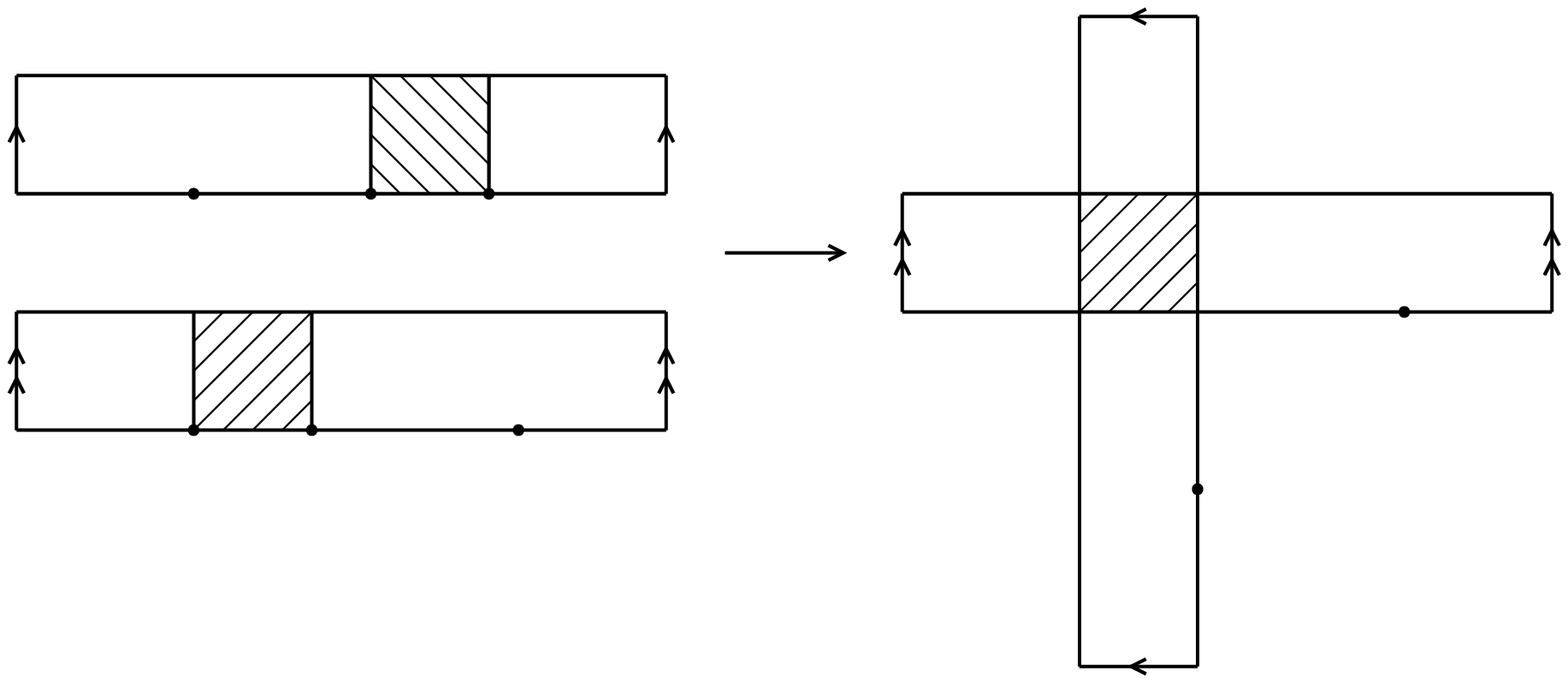}}
\put(0.5,5){\small $\An_t$}
\put(0.5,9){\small $\An_s$}
\put(15.5,7.2){\small $\sim$}
\put(5,3.4){\scriptsize $(2q,0)$}
\put(7,3.4){\scriptsize $(2q+1,0)$}
\put(11,3.4){\scriptsize $(0,0)$}
\put(5.4,7.5){\scriptsize $(0,0)$}
\put(8,7.5){\scriptsize $(2p,0)$}
\put(10,7.5){\scriptsize $(2p+1,0)$}
\end{picture}}
\medskip
\centerline{{\bf Figure 7.4.} Identification of annuli.}
\end{figure}
%%%%%%%%%%%%%%%%%

We identify each annulus $\An_s$ with its image in $\Sigma$, and we denote by $a_s$ its central curve. 
Note that $a_s$ is an essential circle, $a_s \cap a_t = \emptyset$ if $m_{s\,t} =2$, and $| a_s \cap a_t 
| = 1$ if $m_{s\,t}=3$. So, by Proposition 7.1, the mapping $\sigma_s \mapsto \sigma_{a_s}$, $s \in S$, 
induces a geometric representation $\rho_{\CP}: G_\Gamma \to \MM (\Sigma)$.

We have $\rho_{CP} = \rho_{\PV}$ if $\Gamma$ is a tree. (Note that it may happen that 
$\rho_{\PV}$ is not defined if $\Gamma$ is not a tree.) If $\Gamma= \tilde A_n$, then 
$\rho_\CP$ is 
faithful  (while, by \cite{Labru1}, $\rho_\PV$ is not faithful in this case).

\subsection{Presentations}

Let $\Sigma_{g,r}$ be a surface of genus $g \ge 1$ with $r \ge 0$ boundary components, and let $\PP_n$ 
be a collection of $n$ punctures in the interior of $\Sigma_{g,r}$, where $n \ge 0$.

Assume first that $r \ge 1$. Consider the essential circles $a_0, a_1, \dots, a_r, b_1, b_2, 
\linebreak
\dots, 
b_{2g-1}, c, 
d_1, \dots, d_{r-1}$, and the essential arcs $e_1, e_2, \dots, e_{n-1}$ drawn in Figure 7.5. Note 
that there is no $c$ if $g=1$, there is no $d_i$ if $r=1$, there is no $a_r$ if $n=0$, and there is no 
$e_i$ if $n=0$ or $1$. Let $\Gamma (g,r,n)$ be the Coxeter graph drawn in Figure 7.6. One can show 
that the set
\[
\{\sigma_{a_0}, \sigma_{a_1}, \dots, \sigma_{a_r}, \sigma_{b_1}, \sigma_{b_2}, \dots, \sigma_{b_{2g-
1}}, \sigma_c, \sigma_{d_1}, \sigma_{d_2}, \dots, \sigma_{d_{r-1}}, \tau_{e_1}, \dots \tau_{e_{n-1}}\}
\]
generates $\MM (\Sigma_{g,r}, \PP_n)$. On the other hand, by Proposition 7.1, the mapping
\begin{gather*}
x_i \mapsto \sigma_{a_i}\ (0 \le i\le r)\,, \quad y_i \mapsto \sigma_{b_i}\ (1 \le i\le 2g-1)\,, \quad 
z \mapsto \sigma_c\\
u_i \mapsto \sigma_{d_i}\ (1 \le i\le r-1)\,, \quad v_j \mapsto \tau_{e_j}\ (1 \le j\le n-1)\,,
\end{gather*}
induces a homomorphism $\rho: G_{\Gamma (g,r,n)} \to \MM (\Sigma_{g,r}, \PP_n)$. So, in order to obtain 
a presentation for $\MM (\Sigma_{g,r}, \PP_n)$, it suffices to find normal generators for $\Ker\, 
\rho$. This was done in \cite{Matsu1} for $r=1$ and $n=0$, and in \cite{LabPar1} for the other cases.

%%%%%%%%%%%%%%%%%%
\begin{figure}[htb]
\centerline{
\setlength{\unitlength}{.4cm}
\begin{picture}(28,12)
\put(0.5,0){\includegraphics[width=10.8cm]{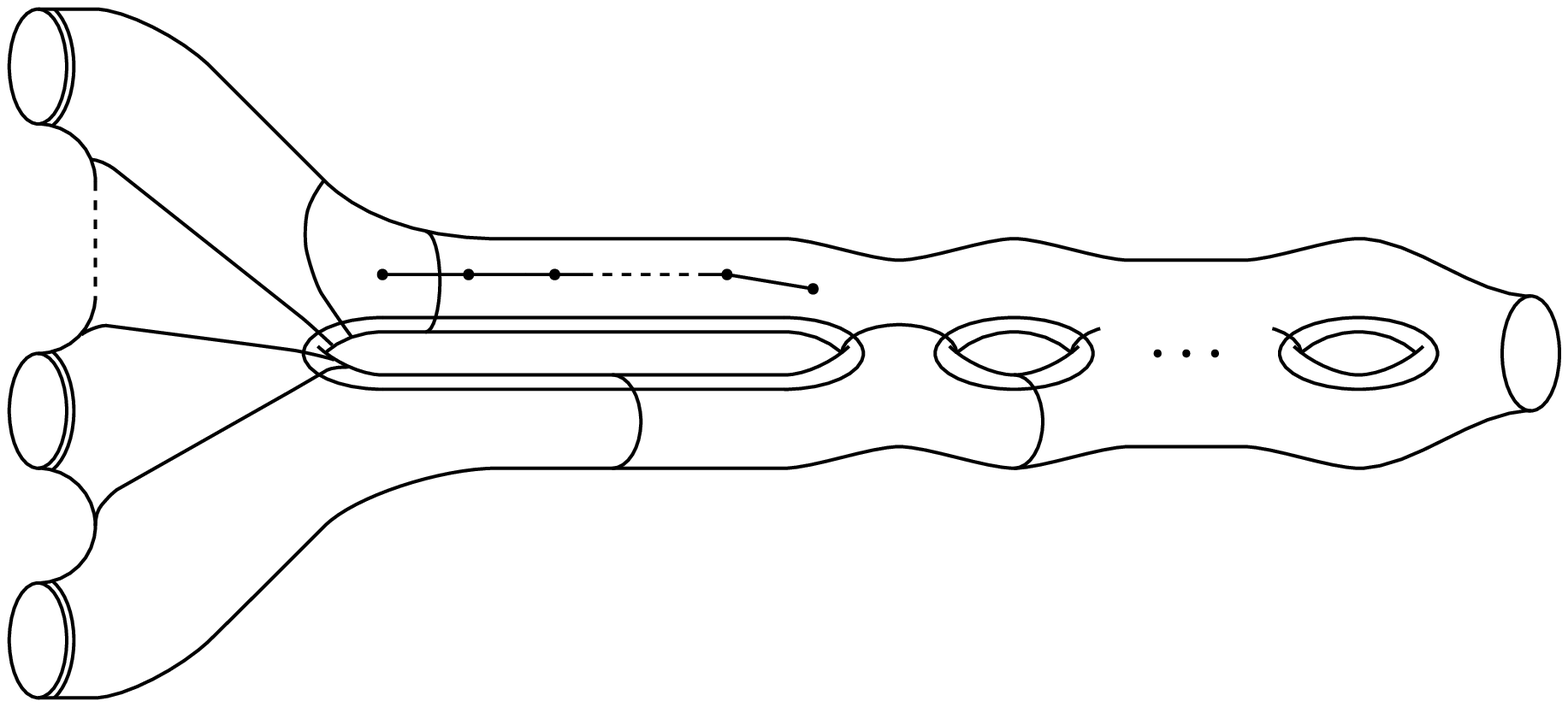}}
\put(10.5,3.5){\small $a_0$}
\put(3.5,3.5){\small $a_1$}
\put(3,5.5){\small $a_2$}
\put(3,9.2){\small $a_{r-2}$}
\put(6,9.5){\small $a_{r-1}$}
\put(7.5,8.5){\small $a_r$}
\put(13,4.7){\small $b_1$}
\put(15.5,6.8){\small $b_2$}
\put(17.5,7){\small $b_3$}
\put(23.5,7){\small $b_{2g-1}$}
\put(18,3.5){\small $c$}
\put(0.5,2.5){\small $d_1$}
\put(0.5,6.5){\small $d_2$}
\put(0.5,9){\small $d_{r-1}$}
\put(7.1,7){\small $e_1$}
\put(9,7){\small $e_2$}
\put(13.3,6.8){\small $e_{n-1}$}
\end{picture}}
\medskip
\centerline{{\bf Figure 7.5.} Generators of $\MM (\Sigma_{g,r}, \PP_n)$.}
\end{figure}
%%%%%%%%%%%%%%

%%%%%%%%%%%%%%%
\begin{figure}[htb]
\centerline{
\setlength{\unitlength}{.4cm}
\begin{picture}(21,6)
\put(0.5,1){\includegraphics[width=7.8cm]{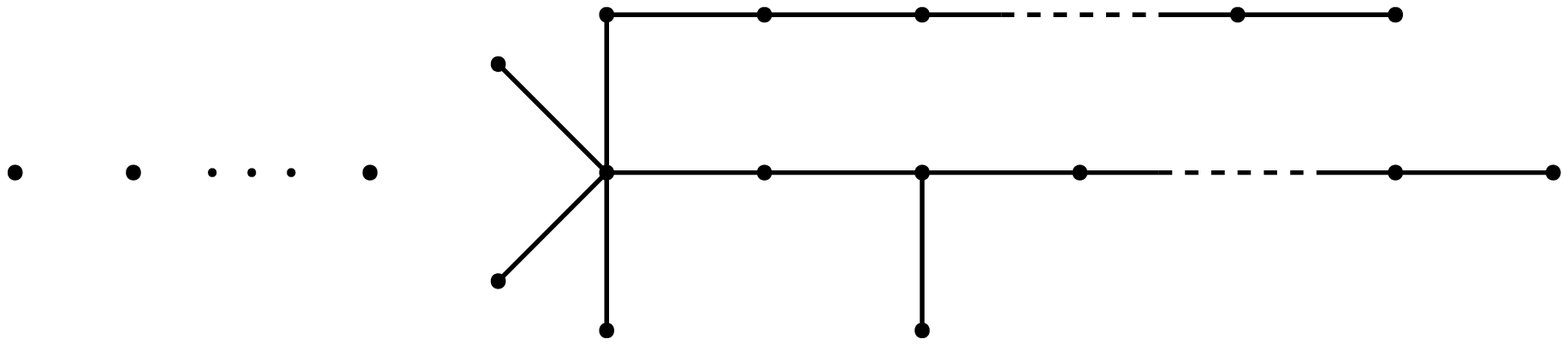}}
\put(7.6,0.5){\small $x_0$}
\put(6,1){\small $x_1$}
\put(5.5,4.8){\small $x_{r-1}$}
\put(8.1,4.5){\small $x_r$}
\put(8.3,2.5){\small $y_1$}
\put(9.7,2.5){\small $y_2$}
\put(12.2,2.5){\small $y_3$}
\put(13.5,2.5){\small $y_4$}
\put(19,2.5){\small $y_{2g-1}$}
\put(11.9,0.5){\small $z$}
\put(0,2.5){\small $u_1$}
\put(1.7,2.5){\small $u_2$}
\put(4.5,2.5){\small $u_{r-1}$}
\put(9.6,4.5){\small $v_1$}
\put(11.7,4.5){\small $v_2$}
\put(17.5,4.5){\small $v_{n-1}$}
\put(8.7,5.4){\small $4$}
\end{picture}}
\medskip
\centerline{{\bf Figure 7.6.} The Coxeter graph $\Gamma (g,r,n)$.}
\end{figure}
%%%%%%%%%%%%%%%%%%%%%%%

One can use the same kind of arguments for the case $r=0$. Consider the essential circles $a_0, a_1, 
b_1, b_2, \dots, b_{2g-1}, c$, and the essential arcs $e_1, e_2, \dots, e_{n-1}$ drawn in Figure 7.7. Then the 
set
\[
\{ \sigma_{a_0}, \sigma_{a_1}, \sigma_{b_1}, \sigma_{b_2}, \dots, \sigma_{b_{2g-1}}, \sigma_c, 
\tau_{e_1}, \tau_{e_2}, \dots, \tau_{e_{n-1}}\}
\]
generates $\MM (\Sigma_{g,0}, \PP_n)$, and the mapping
\begin{gather*}
x_i \mapsto \sigma_{a_i}\ (i=0,1)\,, \quad y_i \mapsto \sigma_{b_i}\ (1 \le i\le 2g-1)\,,\\
 z \mapsto \sigma_c\,, \quad v_j \mapsto \tau_{e_j}\ (1 \le j\le n-1)\,,
\end{gather*}
induces a homomorphism $\rho: G_{\Gamma (g,1,n)} \to \MM (\Sigma_{g,0}, \PP_n)$. Here again, the kernel 
of $\rho$ was calculated in \cite{Matsu1} for $n=0$, and in \cite{LabPar1} for $n \ge 1$.

%%%%%%%%%%
\begin{figure}[htb]
\centerline{
\setlength{\unitlength}{.4cm}
\begin{picture}(23,6)
\put(0,1){\includegraphics[width=9.2cm]{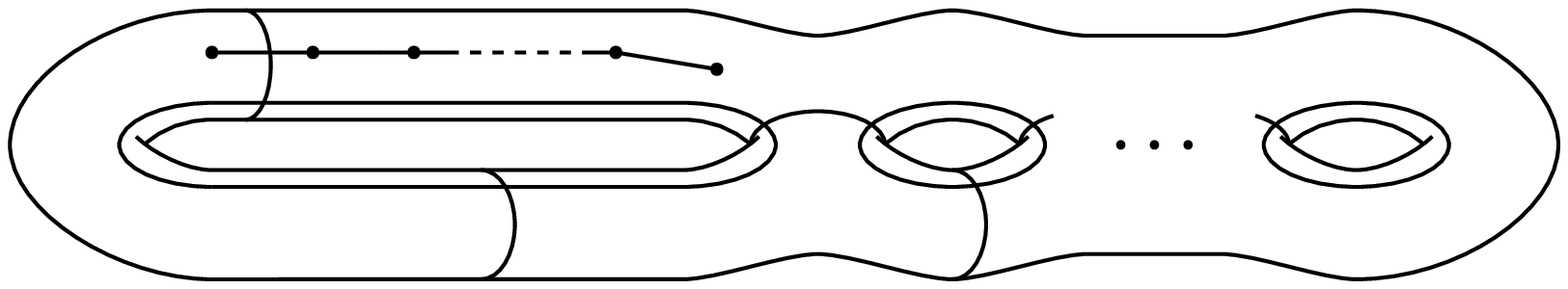}}
\put(7,0.5){\small $a_0$}
\put(3,5.4){\small $a_1$}
\put(9,1.7){\small $b_1$}
\put(11.5,3.8){\small $b_2$}
\put(13.5,4){\small $b_3$}
\put(19.5,4){\small $b_{2g-1}$}
\put(14,0.5){\small $c$}
\put(3,3.9){\small $e_1$}
\put(5,3.9){\small $e_2$}
\put(9,3.9){\small $e_{n-1}$}
\end{picture}}
\medskip
\centerline{{\bf Figure 7.7.} Generators of $\MM (\Sigma_{g,0}, \PP_n)$.}
\end{figure}
%%%%%%%%%%

In order to state the results of \cite{Matsu1} and \cite{LabPar1}, we need the following notations. Let 
$\Gamma$ be a Coxeter graph, let $M=(m_{s\,t})_{s,t \in S}$ be the Coxeter matrix of $\Gamma$, and let 
$(G, \Sigma)$ be the Artin system of type $\Gamma$. For $X \subset S$, we denote by $\Gamma_X$ the full 
subgraph of $\Gamma$ generated by $X$, we set $\Sigma_X= \{ \sigma_s; s \in X\}$, and we denote by 
$G_X$ the subgroup of $G$ generated by $\Sigma_X$. By \cite{VdLek1}, $(G_X, \Sigma_X)$ is the Artin 
system of type $\Gamma_X$ (see also \cite{Paris3}). If $\Gamma_X$ is of spherical type, then we denote 
by $\Delta (X)$ the Garside element of $(G_X, \Sigma_X)$, viewed as an element of $G$.

\begin{theorem}[Matsumoto \cite{Matsu1}]
\begin{enumerate}
\item
$\MM (\Sigma_{g,1})$ is isomorphic with the quotient of $G_{\Gamma (g,1,0)}$ by the following relations
\[
\begin{array}{rcl}
(R1)\quad & \Delta (y_1,y_2,y_3,z)^4 = \Delta (x_0, y_1,y_2,y_3,z)^2 &\ \text{if } g\ge 2\\
\noalign{\smallskip}
(R2)\quad & \Delta (y_1,y_2,y_3,y_4,y_5,z)^2 = \Delta (x_0, y_1,y_2,y_3,y_4,y_5, z) &\ \text{if } 
g\ge 3
\end{array}
\]
\item
$\MM (\Sigma_{g,0})$ is isomorphic with the quotient of $G_{\Gamma (g,1,0)}$ by the above relations 
(R1) and (R2) together with
\[
(R3)\quad\begin{array}{cl}
(x_0y_1)^6=1 &\quad\text{if } g=1\\
\noalign{\smallskip}
x_0^{2g-2} = \Delta (y_2, y_3, z, y_4, \dots, y_{2g-1} ) &\quad\text{if } g\ge 2
\end{array}
\]
\end{enumerate}
\end{theorem}

\begin{theorem}[Labru\`ere, Paris \cite{LabPar1}]
Let $g \ge 1$, $r\ge 1$, and $n \ge 0$. Then 
$\MM (\Sigma_{g,r}, \PP_n)$ is isomorphic with the quotient of $G_{\Gamma (g,r,n)}$ by the following 
relations.
\begin{itemize}
\item
Relations from $\MM (\Sigma_{g,1})$.
\[
\begin{array}{rcl}
(R1)\quad & \Delta (y_1,y_2,y_3,z)^4 = \Delta (x_0, y_1,y_2,y_3,z)^2 &\ \text{if } g\ge 2\\
\noalign{\smallskip}
(R2)\quad & \Delta (y_1,y_2,y_3,y_4,y_5,z)^2 = \Delta (x_0, y_1,y_2,y_3,y_4,y_5, z) &\ \text{if } 
g\ge 3
\end{array}
\]
\item
Relations of commutation.
\[
\begin{array}{rp{3cm}l}
(R3)\quad& \multicolumn{2}{l}{\quad x_k \cdot \Delta (x_{i+1}, x_j, y_1)^{-1} x_i \Delta (x_{i+1}, x_j ,y_1)}\\ 
&\multicolumn{2}{l}{= \Delta (x_{i+1}, x_j, y_1)^{-1} x_i \Delta (x_{i+1}, x_j, y_1) \cdot x_k}\\ 
&&\text{if } 0 \le k<j<i\le r-1\\
\noalign{\smallskip}
(R4)\quad& 
\multicolumn{2}{l}{\quad y_2 \cdot \Delta (x_{i+1}, x_j, y_1)^{-1} x_i \Delta (x_{i+1}, x_j ,y_1)}\\
 &\multicolumn{2}{l}{= \Delta (x_{i+1}, x_j, y_1)^{-1} x_i \Delta (x_{i+1}, x_j, y_1) \cdot y_2}\\
 &&\text{if } 0 \le j<i\le r-1 \text{ and } g \ge 2
\end{array}
\]
\item
Expressions of the $u_i$'s.
\[
\begin{array}{rp{3cm}l}
(R5)\quad &\multicolumn{2}{l}{u_1 = \Delta (x_0, x_1, y_1, y_2,y_3,z) \cdot \Delta (x_1, y_1,y_2,y_3,z)^{-2}}\\
&&\text{if }g \ge 2\\
\noalign{\smallskip}
(R6)\quad &\multicolumn{2}{l}{u_{i+1} = \Delta (x_i, x_{i+1}, y_1,y_2,y_3,z) \cdot \Delta (x_{i+1}, y_1,y_2,y_3,z)^{-2}}\\
&\multicolumn{2}{l}{\hskip1.5cm \cdot \Delta (x_0, x_{i+1}, y_1)^2 \cdot \Delta (x_0, x_i, x_{i+1}, y_1)^{-1}}\\
&&\text{if } 1 \le i\le r-2 \text{ and } g \ge 2
\end{array}
\]
\item
Other relations.
\[
\begin{array}{rp{3cm}l}
(R7)\quad&\multicolumn{2}{c}{ \Delta (x_{r-1}, x_{r}, y_1,v_1) = \Delta (x_{r}, y_1,v_1)^2}\\ &&\text{if } n \ge 2\\
\noalign{\smallskip}
(R8a)\quad &\multicolumn{2}{c}{\Delta (x_0, x_1, y_1,y_2,y_3,z) = \Delta (x_1, y_1,y_2,y_3,z)^2}\\
&&\text{if }n \ge 1,\ g\ge 2,\text{ and } r=1\\
\noalign{\smallskip}
(R8b)\quad &\multicolumn{2}{c}{\Delta (x_{r-1}, x_r, y_1,y_2,y_3,z) \cdot \Delta (x_r, y_1,y_2,y_3,z)^{-2}}\\
&\multicolumn{2}{c}{= \Delta (x_0, x_{r-1}, x_r, y_1) \cdot \Delta (x_0, x_r, y_1)^{-2}}\\ 
&&\text{if } n\ge 1,\ g\ge 2, \text{ and } r\ge 2
\end{array}
\]
\end{itemize}
\end{theorem}

Note that only the relations (R1), (R2), (R7), and (R8a) remain in the presentation if $r=1$, and (R8a) 
must be replaced by (R8b) if $r\ge 2$. Note also that, if $g \ge 2$, then $u_1, \dots, u_{r-1}$ can be 
removed from the generating set. However, to do so, one must add new long relations.

\begin{theorem}[Labru\`ere, Paris \cite{LabPar1}]
Let $g\ge 1$ and $n \ge 1$. Then 
\linebreak
$\MM (\Sigma_{g,0}, \PP_n)$ 
is isomorphic with the quotient of $G_{\Gamma (g,1,n)}$ by the following 
relations.
\begin{itemize}
\item
Relations from $\MM (\Sigma_{g,1}, \PP_n)$.
\[
\begin{array}{rp{3cm}l}
(R1)\quad &\multicolumn{2}{c}{\Delta (y_1,y_2,y_3,z)^4 = \Delta (x_0, y_1,y_2,y_3,z)^2}\\ &&\text{if } g\ge 2\\
\noalign{\smallskip}
(R2)\quad &\multicolumn{2}{c}{ \Delta (y_1,y_2,y_3,y_4,y_5,z)^2 = \Delta (x_0, y_1,y_2,y_3,y_4,y_5, z)}\\
&&\text{if } g\ge 3\\
\noalign{\smallskip}
(R7)\quad&\multicolumn{2}{c}{\Delta (x_0, x_1, y_1,v_1) = \Delta (x_1, y_1,v_1)^2}\\
&&\text{if } n \ge 2\\
\noalign{\smallskip}
(R8a)\quad &\multicolumn{2}{c}{\Delta (x_0, x_1, y_1,y_2,y_3,z) = \Delta (x_1, y_1,y_2,y_3,z)^2}\\
&&\text{if }n \ge 1 \text{ and } g\ge 2
\end{array}
\]
\item
Other relations.
\[
\begin{array}{rp{3cm}l}
(R9a) \quad &\multicolumn{2}{c}{x_0^{2g-n-2} \cdot \Delta (x_1, v_1, \dots, v_{n-1}) = \Delta (z, y_2, \dots, y_{2g-1})^2}\\ 
&&\text{if }g\ge 2\\
\noalign{\smallskip}
(R9b) \quad &\multicolumn{2}{c}{x_0^n = \Delta (x_1, v_1, \dots, v_{n-1})}\\
&&\text{if }g=1\\
\noalign{\smallskip}
(R9c) \quad & \multicolumn{2}{c}{\Delta (x_0, y_1)^4 = \Delta (v_1, \dots, v_{n-1})^2}\\
&&\text{if } g=1
\end{array}
\]
\end{itemize}
\end{theorem}

\begin{note}
Presentations of $\MM (\Sigma_{g,r})$, also in terms of Artin groups,
with more generators but simpler relations, were 
obtained by Gervais in \cite{Gerva1}. On the other hand, a unified proof of all these presentations 
can be found in \cite{Benve1}.
\end{note}

\subsection{Classification}

This subsection is an account of Castel's results \cite{Caste2} on the geometric representations of the 
braid group $\BB_n$ on mapping class groups of surfaces of genus $g \le \frac{n-1}{2}$.

Suppose first that $n$ is odd, $n \ge 5$. Write $n=2k+1$, where $k \ge 2$. Let $r \ge 0$. We present 
the surface $\Sigma_{k,r}$ as the union of three subsurfaces, $\Omega_0$, $\A$, and $\Omega_1$, where 
$\Omega_0$ is a surface of genus $k$ with one boundary component, $c$, $\Omega_1$ is a surface of genus 
$0$ with $r+1$ boundary components, $c', d_1, \dots, d_r$, and $\A$ is an annulus bounded by $c$ and 
$c'$ (see Figure 7.8). Consider the essential circles $a_1, a_2, \dots, a_{2k}$ drawn in Figure 7.8. 
Then, by Proposition 7.1, there exists a homomorphism $\rho_M: \BB_n \to \MM (\Sigma_{k,r})$ which 
sends $\sigma_i$ to $\sigma_{a_i}$ for all $1 \le i\le n-1=2k$.

%%%%%%%%%%
\begin{figure}[htb]
\centerline{
\setlength{\unitlength}{.4cm}
\begin{picture}(22,12)
\put(-1,0){\includegraphics[width=9cm]{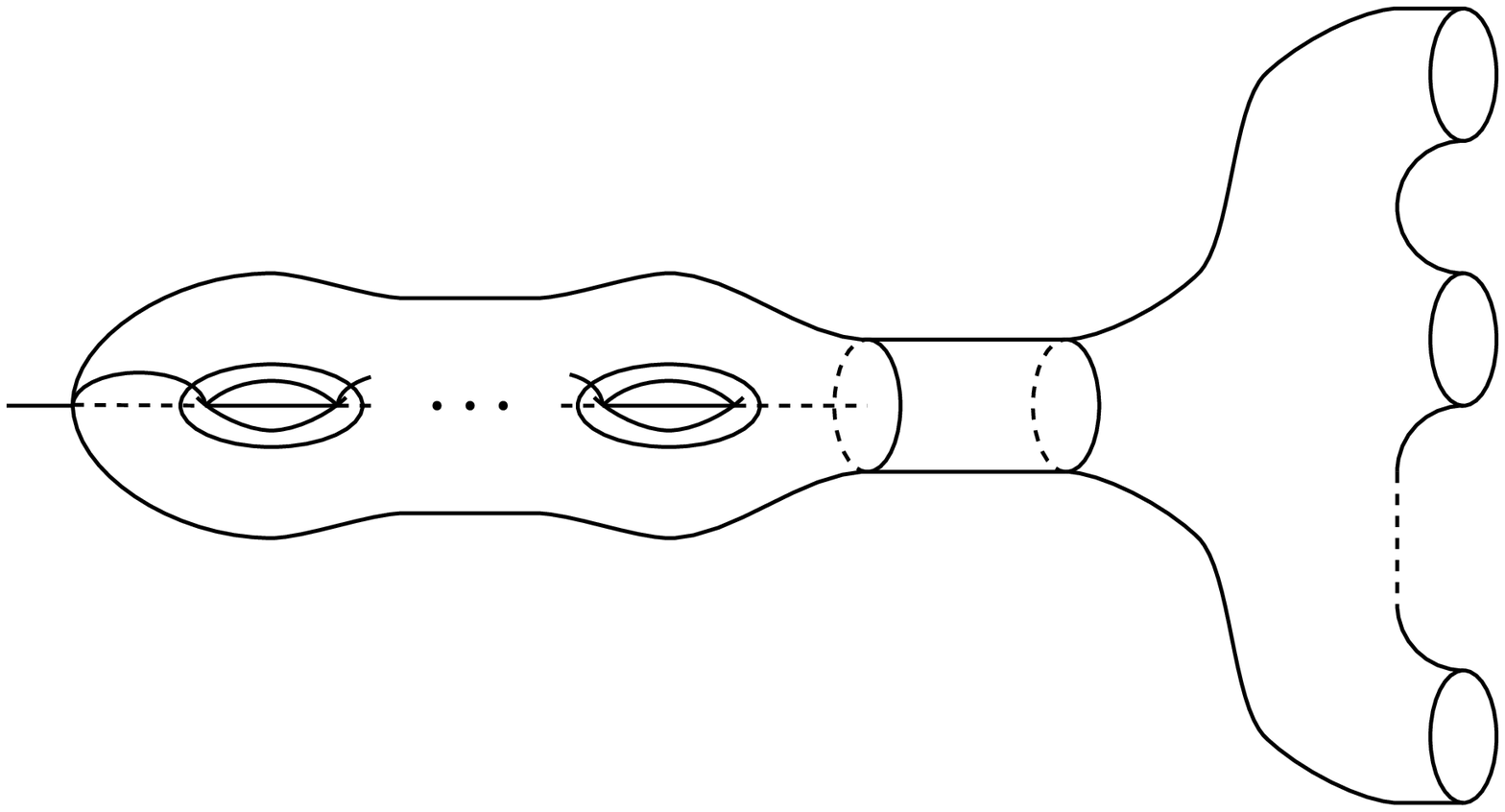}}
\put(0.7,6.7){\small $a_1$}
\put(2.5,7){\small $a_2$}
\put(8.5,7){\small $a_{2k}$}
\put(11.7,7.3){\small $c$}
\put(14.7,7.3){\small $c'$}
\put(21,9.2){\small $d_1$}
\put(21,5.2){\small $d_2$}
\put(21,2.3){\small $d_r$}
\put(5,3){\small $\Omega_0$}
\put(13,4){\small $\A$}
\put(16,2){\small $\Omega_1$}
\put(-1,6.3){\small $D$}
\end{picture}}
\medskip
\centerline{{\bf Figure 7.8.} Decomposition of $\Sigma_{k,r}$ ($n$ odd).}
\end{figure}
%%%%%%%%%%

The statement of Castel's classification of the geometric representations of $\BB_n$ in $\MM 
(\Sigma_{k,r})$ involves the centralizer of $\Im\, \rho_M$ in $\MM (\Sigma_{k,r})$. That is why we start 
with a description of the latter.

The inclusion of $\Omega_1$ in $\Sigma_{k,r}$ induces a homomorphism $\MM (\Omega_1) \to \MM 
(\Sigma_{k,r})$ which is injective (see \cite{ParRol1}). It is easily checked that the image of this 
homomorphism is contained in the centralizer of $\Im\, \rho_M$. Another element of the centralizer is 
the element $u \in \MM (\Sigma_{k,r})$ represented by the homeomorphism $U: \Sigma_{k,r} \to 
\Sigma_{k,r}$ which is the axial symmetry relative to the axis $D$ on $\Omega_0$, a half-twist which 
pointwise fixes $c'$ on the annulus $\A$, and the identity on $\Omega_1$.

\begin{proposition}[Castel \cite{Caste2}]
The centralizer of $\Im\, \rho_M$ in $\MM 
(\Sigma_{k,r})$ is generated by $\MM (\Omega_1) \cup \{ u\}$.
\end{proposition}

If $r=0$, then $\MM (\Omega_1) = \{1\}$, $u$ is of order $2$, and $Z_{\MM (\Sigma_{k,r})} (\Im\, \rho_M) 
= \langle u\rangle$ is cyclic of order 2. If $r=1$, then $\MM (\Omega_1) = \langle \tau_c \rangle$, 
$u^2 = \tau_c$, and $Z_{\MM (\Sigma_{k,r})} (\Im\, \rho_M)= \langle u\rangle$ is an infinite cyclic 
group. If $r=2$, then $Z_{\MM (\Sigma_{k,r})} (\Im\, \rho_M)$ is a free abelian group of rank 3 freely 
generated by $\{ u, \sigma_{d_1}, \sigma_{d_2}\}$. If $r \ge 3$, then $Z_{\MM (\Sigma_{k,r})} (\Im\, 
\rho_M)$ is more complicated.

For $\varepsilon \in \{ \pm 1 \}$ and $z \in Z_{\MM (\Sigma_{k,r})} (\Im\, \rho_M)$, 
the mapping $\sigma_i \mapsto \sigma_{a_i}^{\varepsilon} z$, $1 \le 
i\le n-1$, induces a homomorphism $\rho_{M} (\varepsilon, z): \BB_n \to \MM (\Sigma_{k,r})$ called the {\it 
transvection of $\rho_M$ by $(\varepsilon,z)$}\index{Transvection}. On the other hand, 
a homomorphism $\varphi: \BB_n \to G$, 
where $G$ is a group, is called {\it cyclic} if there exists $\alpha \in G$ such that $\varphi 
(\sigma_i) = \alpha$ for all $1 \le i \le n-1$.

\begin{theorem}[Castel \cite{Caste2}]
Suppose $n$ odd, $n \ge 5$, and set $n=2k+1$. Let $g \ge 
0$ and $r \ge 0$.
\begin{enumerate}
\item
If $g <k$, then all the homomorphisms $\varphi : \BB_n \to \MM (\Sigma_{g,r})$ are cyclic.
\item
All the non-cyclic homomorphisms $\varphi: \BB_n \to \MM (\Sigma_{k,r})$ are conjugate to transvections 
of $\rho_M$.
\item
The homomorphism $\rho_M : \BB_n \to \MM (\Sigma_{k,r})$ is injective if and only if $r \ge 1$.
\end{enumerate}
\end{theorem}

Now, we suppose that $n$ is even, $n \ge 6$, and we set $n=2k+2$. We choose $r_1,r_2 \ge 0$ such that 
$r_1+r_2 = r$ and we represent the surface $\Sigma_{k,r}$ as the union of three subsurfaces, a surface 
$\Omega_0$ of genus $k$ with two boundary components, $c_1$ and $c_2$, a surface $\Omega_1$ of genus 
$0$ with $r_1+1$ boundary components $c_1, d_1, \dots, d_{r_1}$, and a surface $\Omega_2$ of genus $0$ 
with $r_2+1$ boundary components $c_2, d_{r_1+1}, \dots, d_{r_1+r_2}$ (see Figure 7.9). Consider the 
essential circles $a_1, \dots, a_{n-1}$ drawn in Figure 7.9. Then, by Proposition 7.1, there exists a 
homomorphism $\rho_{M}(r_1,r_2): \BB_n \to \MM (\Sigma_{k,r})$ which sends $\sigma_i$ to $\sigma_{a_i}$ 
for all $1 \le i \le n-1$.

%%%%%%%%%%
\begin{figure}[htb]
\centerline{
\setlength{\unitlength}{.35cm}
\begin{picture}(30,12)
\put(0.5,0){\includegraphics[width=10.15cm]{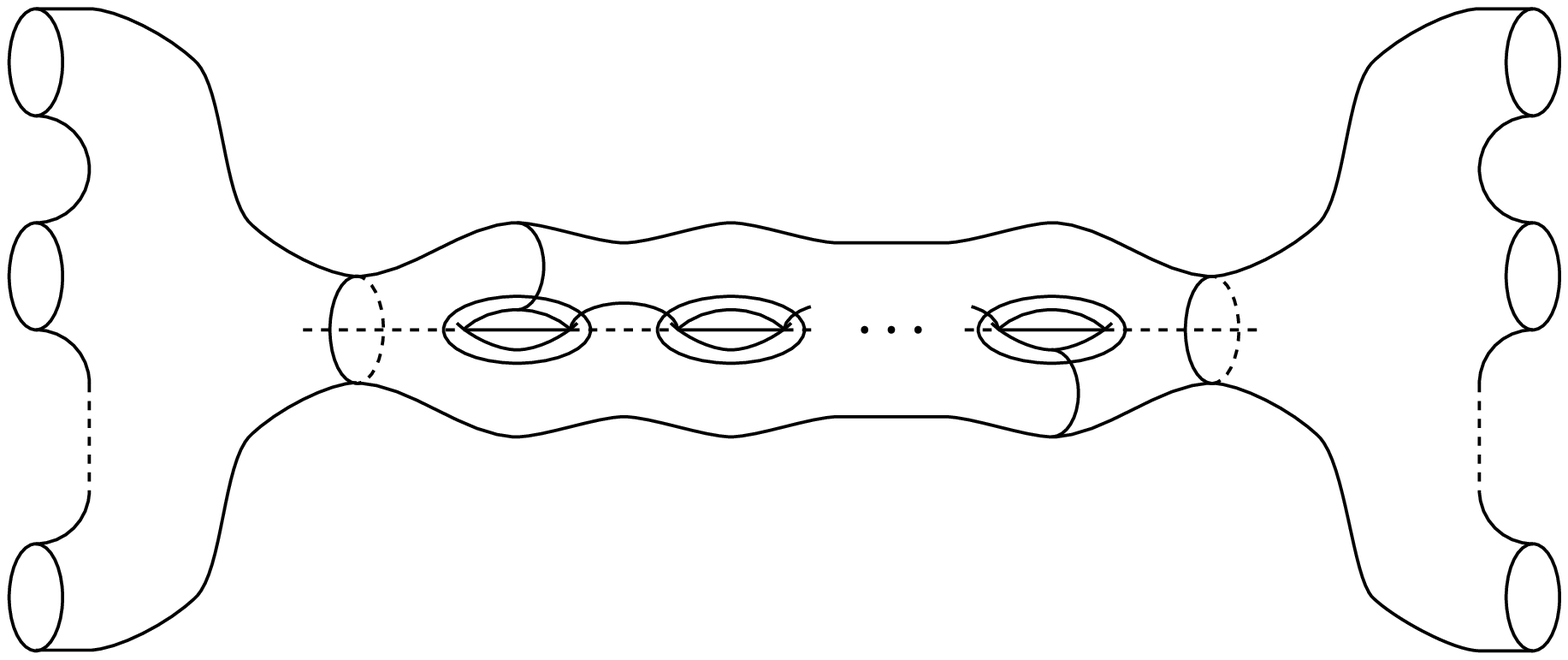}}
\put(9.7,8.3){\small $a_1$}
\put(9.7,4.8){\small $a_2$}
\put(11.5,6.8){\small $a_3$}
\put(13.5,7){\small $a_4$}
\put(19.5,7){\small $a_{2k}$}
\put(19.5,3.2){\small $a_{2k+1}$}
\put(6.7,4.2){\small $c_1$}
\put(22.7,4.2){\small $c_2$}
\put(0.5,9.2){\small $d_1$}
\put(0.5,5.2){\small $d_2$}
\put(0.5,2.5){\small $d_{r_1}$}
\put(29,9.1){\small $d_{r_1+1}$}
\put(29,2.4){\small $d_r$}
\put(3,6){\small $\Omega_1$}
\put(14,2.5){\small $\Omega_0$}
\put(25.5,6){\small $\Omega_2$}
\put(5.2,6){\small $D$}
\end{picture}}
\medskip
\centerline{{\bf Figure 7.9.} Decomposition of $\Sigma_{k,r}$ ($n$ even).}
\end{figure}
%%%%%%%%%%

The inclusions $\Omega_1, \Omega_2 \subset \Sigma_{k,r}$ induce a homomorphism $\MM (\Omega_1) \times 
\MM (\Omega_2) \to \MM (\Sigma_{k,r})$ which is injective (see \cite{ParRol1}), and we have:

\begin{proposition}[Castel \cite{Caste2}]
\begin{enumerate}
\item
If $r>0$, then the centralizer of 
\linebreak
$\Im\, \rho_{M}(r_1,r_2)$ in $\MM (\Sigma_{k,r})$ is $\MM (\Omega_1) 
\times \MM (\Omega_2)$.
\item
If $r=0$, then the centralizer of $\Im\, \rho_{M}(r_1,r_2)$ in $\MM (\Sigma_{k,r})$ is a cyclic group of 
order 2 generated by an element represented by the axial symmetry relative to the axis $D$ of Figure 
7.9.
\end{enumerate}
\end{proposition}

For $\varepsilon \in \{ \pm 1\}$ and $z \in Z_{\MM (\Sigma_{k,r})} (\Im\, \rho_{M}(r_1,r_2))$, 
the mapping $\sigma_i \mapsto \sigma_{a_i}^{\varepsilon} 
z$, $1 \le i\le n-1$, induces a homomorphism $\rho_{M}(r_1,r_2,\varepsilon,z): \BB_n \to \MM (\Sigma_{k,r})$ 
called 
the {\it transvection of $\rho_{M}(r_1,r_2)$ by $(\varepsilon,z)$}\index{Transvection}.

\begin{theorem}[Castel \cite{Caste2}]
Suppose $n$ even, $n\ge 6$, and set $n=2k+2$. Let $g \ge 
0$ and $r \ge 0$.
\begin{enumerate}
\item
If $g <k$, then all the homomorphisms $\varphi: \BB_n \to \MM (\Sigma_{g,r})$ are cyclic.
\item
If $\varphi: \BB_n \to \MM (\Sigma_{k,r})$ is a non-cyclic homomorphism, then there exist $r_1,r_2 \ge 
0$ such that $r_1+r_2=r$ and $\varphi$ is conjugate to a transvection of $\rho_{M}(r_1,r_2)$.
\item
Let $r_1,r_2 \ge 0$ such that $r_1+r_2=r$. The homomorphism $\rho_{M}(r_1,r_2): \BB_n \to 
\MM(\Sigma_{k,r})$ is injective if and only if $r_1 \ge 1$ and $r_2 \ge 1$.
\end{enumerate}
\end{theorem}

Recall that, for a group $G$, $\Out (G)$ denotes the group of outer automorphisms of $G$. Now, Theorems 
7.6 and 7.8 can be used for new proofs of the following two theorems.

\begin{theorem}[Dyer, Grossman \cite{DyeGro1}]
We have $\Out (\BB_n) = \Z/2 \Z$ if $n \ge 5$.
\end{theorem}

\begin{theorem}[Ivanov \cite{Ivano2}, McCarthy \cite{McCar1}]
Let $g \ge 2$ and $r \ge 0$. 
Then
\[
\Out (\MM (\Sigma_{g,r})) = \left\{ \begin{array}{ll}
\{1\} &\quad\text{if } r\ge 1\\
\Z/2\Z &\quad \text{if } r=0 \text{ and } g\ge 3\\
\Z/2\Z \times \Z/2\Z &\quad\text{if } r=0 \text{ and } g=2
\end{array}\right.
\]
\end{theorem}

%%%%%%%%%%%%%%%%%%%%%%%%

%%%%%%%%%%%%%%%%%%%%%%%%%%%
%\printindex

\end{document}